\pdfoutput=1
\documentclass[11pt,dvipsnames]{article}
\usepackage[T1]{fontenc}
\usepackage[utf8]{inputenc}
\usepackage{lmodern}
\usepackage[margin=1in]{geometry}
\usepackage{amsmath,amssymb,amsthm,mathtools}
\usepackage{booktabs}
\usepackage{graphicx}
\usepackage{xcolor}
\usepackage{microtype}
\usepackage[colorlinks=true,linkcolor=NavyBlue,citecolor=PineGreen,urlcolor=RoyalBlue,pdfusetitle]{hyperref}

\newtheorem{theorem}{Theorem}[section]
\newtheorem{proposition}[theorem]{Proposition}
\newtheorem{lemma}[theorem]{Lemma}
\newtheorem{corollary}[theorem]{Corollary}
\theoremstyle{definition}
\newtheorem{definition}[theorem]{Definition}
\theoremstyle{remark}
\newtheorem{remark}[theorem]{Remark}

\newcommand{\R}{\mathbb{R}}
\newcommand{\Z}{\mathbb{Z}}
\newcommand{\Wopt}{W_{\mathrm{opt}}}

\newcommand{\PaperAuthor}{Kevin M. Boergens}
\newcommand{\PaperDepartments}{Departments of Biological Sciences and Physics}
\newcommand{\PaperAffiliation}{University of Illinois Chicago}
\newcommand{\PaperEmail}{boergens@uic.edu}

\title{\textbf{Optimal Patrol Widths in the Infinite-Corridor \mbox{Barrier-Guarding} Game}}
\author{\PaperAuthor\thanks{Electronic address: \texttt{\href{mailto:\PaperEmail}{\PaperEmail}}.}\\[0.4ex]{\normalsize\PaperDepartments}\\{\normalsize\PaperAffiliation}}
\hypersetup{pdftitle={Optimal Patrol Widths in the Infinite-Corridor Barrier-Guarding Game},pdfauthor={Kevin M. Boergens}}
\date{November 2026}

\begin{document}
\maketitle

\begin{abstract}
A guard of capture radius $r$ commits to a patrol in an infinite parallel-walled corridor; an intruder knows the patrol and moves at most $q$ times as fast. Let $W_{\mathrm{opt}}(q,H)$ be the supremal width guarded by (possibly curved, nonperiodic) patrols with longitudinal excursion $|Y|\le H$ against crossings in a known direction, and $W_{\mathrm{opt}}(q)=\sup_{H\ge 0}W_{\mathrm{opt}}(q,H)$. Periodic polygonal patrols give the same suprema; a transverse sweep ($H=0$) guards $W_{\mathrm{bf}}(q)=r\sqrt{4+q^{-2}}$; and $W_{\mathrm{opt}}(q,H)\le W_{\mathrm{bf}}(q)+(1-q)_+H/q$ with initial $H$-slope $(1-q)/q$ on $(0,1)$, so the sweep is optimal for $q\ge 1$. For $1/3 < q < 1$, straight wall runners force $W_{\mathrm{opt}}(q)\le F(q)=2r\sqrt{1+(3q^2+2q-1)^{-1}}$. A four-leg bowtie attains $F(q)$ on $[q_a,1)$ ($0.457892036651 < q_a < 0.457892036652$), we find $W_{\mathrm{opt}}(q,H)$ for all $H\ge 0$ on $[1/2,1)$, and certifying $F(q_a)-\varepsilon r$ by a space-time barrier takes $\Theta(\varepsilon^{-1/4})$ triangles. On $(1/3,q_a)$, paths rounding the wall-parallel legs force $W_{\mathrm{opt}}(q) < F(q)$, with $F(q)-W_{\mathrm{opt}}(q)=\Theta(r(q_a-q)^2)$ as $q\uparrow q_a$. Eight legs beat every symmetric full-speed four-leg bowtie on $(0,q_a)$; near $q_a$, five legs are the fewest beating all periodic four-leg patrols (including asymmetric routes with waits), and twelve lower the quadratic deficit coefficient from $51.2835\ldots$ to below $48.4251$. As $q\downarrow 0$, $2r/q-W_{\mathrm{opt}}(q)$ is below $1.235353\,r$ for $q\le 1/40$ and at least $2^{-21}r$ for $q\le 2^{-24}$ (versus $2r$ in the limit for symmetric full-speed four-leg bowties), while $qW_{\mathrm{opt}}(q,H)/r\to\min\{2,1+H/r\}$. Blocking both directions makes the sweep optimal for $q\ge 1/2$, and as $q\downarrow 0$ unknown direction costs $\Theta(r/\sqrt{q})$ of width and requires $H(q)/r\to\infty$ to reach $2r/q$.

\end{abstract}

\tableofcontents
\section{Introduction}
\subsection{The channel patrol problem}
A guard with capture radius $r$ commits to a continuous patrol route inside an infinite parallel-walled corridor. An intruder knows the entire schedule in advance and can move in any direction, wait, or reverse at speeds up to $q$ times the guard's maximum speed, succeeding if it traverses the corridor longitudinally in a given direction without ever entering the closed capture disk of radius $r$ around the guard. How wide a corridor can a single route-committed guard protect, and how does the geometry of optimal patrols change across the speed ratio $q$?

Half speed illustrates why moving along the corridor helps. At $q=1/2$, a guard restricted to sweeping horizontally back and forth between the walls protects width $2\sqrt{2}\,r\approx 2.828\,r$. Tilting the route into a symmetric four-leg ``bowtie''---running parallel to each wall against the oncoming traffic and rising diagonally across the corridor, as on the left of Figure~\ref{fig:landscape}---protects the exact optimal width $2\sqrt{7/3}\,r\approx 3.055\,r$, an $8.01\%$ increase over horizontal sweeping.

Washburn posed this channel-patrol problem for a single searcher and compared it with two solvable restrictions: joining the walls into a cylindrical tube, and confining the guard to a transverse back-and-forth sweep~\cite[Sec.~2.C]{Washburn2009}, \cite[pp.~36--37]{Washburn2010}. In the transverse sweep the critical width and optimal half-amplitude are
\begin{equation}\label{eq:washburn-bf}
 W_{\mathrm{bf}}(q)=r\sqrt{4+q^{-2}},\qquad a_0=\frac{r}{2q\sqrt{1+4q^2}}.
\end{equation}
In the tube approximation the critical width for $q<1$ is $2r/q$. Even the transverse value in~\eqref{eq:washburn-bf} exceeds $2r/q$ when $\sqrt{3}/2<q<1$, refuting the natural guess that identifying the walls always helps the searcher~\cite[Sec.~2.C]{Washburn2009}.

Bowtie and oblique patrols have a longer operational history against restricted evaders. Koopman's crossover patrol chooses its lead angle against targets on a fixed straight course~\cite[Sec.~7.1 and Fig.~6]{Koopman1946}, \cite[pp.~200--203]{Koopman1980}; Balinski, Lucas, Owen, and Sand~\cite[App.~I--II]{Balinski1965}, Washburn~\cite{Washburn1976}, and Galyaev and Maslov~\cite{GalyaevMaslov2011} likewise analyze straight-line transitors. Hew~\cite{Hew2015,Hew2017} studies prescribed patrols against maneuvering evaders, and Rowling~\cite{Rowling2019} treats effective sweep width under counter-detection. Here every admissible path, at every starting phase and with arbitrary maneuvers inside the corridor, must be intercepted with certainty.

Testing intruders that move straight upward along vertical lines of constant transverse coordinate $x$ (which we call \emph{lanes}) already yields Washburn's universal area upper bound~\cite[eqs.~(2)--(3)]{Washburn1982}:
\begin{equation}\label{eq:area-bound}
 \Wopt(q)\le W_{\mathrm{area}}(q):=2r\sqrt{1+q^{-2}}.
\end{equation}
Indeed, normalizing the guard's maximum speed to one, the disk swept over $[0,L]$ in straight-runner coordinates $(x,qt-y)$ has area at most $\pi r^2+2r\int_0^L\sqrt{X'^2+(q-Y')^2}\,dt$ by polygonal approximation. Concavity of the square root and $X'^2+Y'^2\le 1$ bound the integrand by $\sqrt{1+q^2}-qY'/\sqrt{1+q^2}$. When $|Y|\le H$, every straight-runner phase in $[H+r,qL-H-r]$ on every lane $x\in[-b,b]$ must be covered during $[0,L]$; comparing areas, dividing by $L$, and letting $L\to\infty$ gives~\eqref{eq:area-bound}.

\subsection{Mathematical model and conventions}

Fix capture radius $r>0$, maximum guard speed $V>0$, maximum intruder speed $U>0$, and dimensionless speed ratio $q := U/V > 0$. Rescaling space by $r$ and time by $r/V$ normalizes units to $r=1$, $V=1$, and $U=q$; we restore explicit factors of $r$ in theorem statements. In normalized units a unit-speed patrol means a trajectory with speed bounded almost everywhere by one, so waiting and slowing down are permitted.

\begin{definition}[Corridor, patrols, and evasion]\label{def:model}
For half-width $b>0$ (corridor width $W=2b$), let $\mathcal{C}_b := [-b,b]\times\R$.
\begin{enumerate}
    \item A \emph{prescribed patrol} is an absolutely continuous trajectory $g=(X,Y):\R\to\mathcal{C}_b$ satisfying the speed bound $\operatorname{ess\,sup}_{t\in\R}\sqrt{X'(t)^2+Y'(t)^2} \le V$. It has \emph{bounded longitudinal excursion} if $H(g):=\sup_{t\in\R}|Y(t)|<\infty$, and is \emph{$T$-periodic} if $g(t+T)=g(t)$ for all $t\in\R$ (which implies $H(g)<\infty$).
    \item A \emph{(one-way upward) transit trajectory} is an absolutely continuous path $p=(x,y):\R\to\mathcal{C}_b$ satisfying $\operatorname{ess\,sup}_t \sqrt{x'(t)^2+y'(t)^2}\le U$ that traverses the corridor upward: $\lim_{t\to-\infty} y(t) = -\infty$ and $\lim_{t\to+\infty} y(t) = +\infty$. (When $H(g)\le H<\infty$, any finite-horizon segment $p:[t_0,t_1]\to\mathcal{C}_b$ with $y(t_0)\le -H-2r$ and $y(t_1)\ge +H+2r$ extends outside $[t_0,t_1]$ at velocity $(0,+U)$ without meeting $\overline{B}_r(g(t))$, so finite-slab and full-corridor transit are equivalent.)
    \item \emph{Closed-disk capture:} The patrol $g$ \emph{captures} $p$ if there exists $t\in\R$ such that
    \[
        \|p(t)-g(t)\|_2 = \sqrt{(x(t)-X(t))^2 + (y(t)-Y(t))^2} \le r.
    \]
    The patrol $g$ \emph{successfully guards} $\mathcal{C}_b$ if it captures every admissible upward transit trajectory $p$.
\end{enumerate}
\end{definition}

\begin{definition}[Minimax guarded corridor widths]\label{def:wopt}
For $q = U/V > 0$ and $H \ge 0$, define the excursion-constrained minimax width
\[
    \Wopt(q, H) := \sup\bigl\{\, 2b > 0 : \exists\text{ patrol } g\text{ on }\mathcal{C}_b\text{ with }\sup_{t\in\R}|Y(t)|\le H\text{ guarding }\mathcal{C}_b \,\bigr\},
\]
and the unrestricted bounded-excursion minimax width $\Wopt(q) := \sup_{H\ge 0} \Wopt(q,H)$, in which each competing patrol is bounded longitudinally ($H(g)<\infty$) but the bound may vary arbitrarily across patrols. Write $D_1(q):=2/q-\Wopt(q)/r$ for the dimensionless slow-speed deficit relative to $2r/q$. Because any speed-$q_2$ trajectory is admissible at speed $q_1 \ge q_2$, $\Wopt(q)$ is nonincreasing in $q$; $\Wopt(q,H)$ is nondecreasing in $H$; and feasible widths form an interval because horizontally clamping a patrol onto a narrower corridor preserves its speed bound and never increases distances within the strip.
\end{definition}

\begin{remark}[Model conventions]\label{rem:conventions}
Three modeling conventions merit explicit record:
\begin{enumerate}
    \item \emph{Stackelberg commitment:} The guard commits to $g(t)$ without real-time feedback on the intruder's position outside radius $r$, while the intruder knows $g$ completely and may curve, pause, or reverse longitudinally at any speed $\le U$.
    \item \emph{Closed vs.\ open capture disk:} Definition~\ref{def:model} uses the closed disk $\{\|p-g\|_2\le r\}$, so grazing contact at distance $r$ counts as capture. Under an open capture disk $\{\|p-g\|_2<r\}$, if $g$ guards $\mathcal{C}_b$ with $\overline{B}_r$, then for every $\varepsilon\in(0,1)$ the contracted trajectory $g_\varepsilon(t) := (1-\varepsilon)g(t/(1-\varepsilon))$ has speed $\le V$ and guards $\mathcal{C}_{(1-\varepsilon)b}$ with closed radius $(1-\varepsilon)r < r$ (hence with $B_r$). Letting $\varepsilon\downarrow 0$ shows that open-disk and closed-disk suprema coincide for every $q>0$.
    \item \emph{One-way and bidirectional variants:} Sections~\ref{sec:global-fast}--\ref{sec:asymptotics} analyze the known-direction (one-way) game across the full range of $q$; Section~\ref{sec:diagnostics} then analyzes the bidirectional (two-way) variant in which a single patrol must block both upward and downward transits simultaneously.
\end{enumerate}
\end{remark}

\subsection{Overview of results}
The main body of the paper establishes global structural properties valid across all speed ratios $q>0$, then characterizes the one-way minimax width $\Wopt(q)$ from fast intruders ($q\ge 1$) down to the slow-intruder limit ($q\downarrow 0$), and finally treats bidirectional guarding in a single self-contained chapter.

\paragraph{Global structure across all $q>0$ (Section~\ref{sec:global-fast}).}
Although Definition~\ref{def:wopt} admits arbitrary curved, variable-speed, nonperiodic trajectories of bounded excursion, periodic polygonal patrols (finitely many straight segments repeated on a fixed period) achieve the exact same width supremum $\Wopt(q,H)$ at every $q>0$ and every excursion cap $H\ge 0$, and hence achieve $\Wopt(q)$ as well (Theorem~\ref{thm:periodic-approximation}). An elementary two-level moving barrier for four-leg bowties (Theorem~\ref{thm:two-level}) pairs with a midline wall-phase inequality (Theorem~\ref{thm:midline}) to determine the exact zero-excursion baseline $\Wopt(q,0)=W_{\mathrm{bf}}(q)=r\sqrt{4+q^{-2}}$ at every $q>0$ and establish the universal excursion bound
\[
    \Wopt(q,H)\le W_{\mathrm{bf}}(q)+\frac{(1-q)_+}{q}H,
\]
with exact initial gain $\Wopt(q,H)=W_{\mathrm{bf}}(q)+\frac{1-q}{q}H+O(H^2/r)$ as $H\downarrow 0$ for every $0<q<1$ (Corollary~\ref{cor:two-level-conseq}).

\paragraph{Fast intruders, $q\ge 1$ (Section~\ref{sec:global-fast}).}
When $q\ge 1$, $(1-q)_+=0$ in the global excursion inequality, so longitudinal motion cannot widen the guarded corridor at all: horizontal back-and-forth sweeping is globally optimal, $\Wopt(q)=W_{\mathrm{bf}}(q)=r\sqrt{4+q^{-2}}$ (Theorem~\ref{thm:midline}(1)).

\paragraph{Above the diffraction threshold, $q_a\le q<1$ (Section~\ref{sec:above-qa}).}
For $1/3<q<1$, intruders running straight upward along the left and right walls force $\Wopt(q)\le F(q)$, where
\[
    F(q)=2r\sqrt{1+(3q^2+2q-1)^{-1}}
\]
is the \emph{wall-runner upper bound} (or \emph{straight-wall ceiling}, Theorem~\ref{thm:wall-upper}). The symmetric constant-speed four-leg bowtie that saturates this two-wall time budget descends parallel to each wall at an inset and crosses diagonally upward at inclination $\sin\alpha=(1-q)/(2q)$. Combining the flat two-level separator above $q_0\approx 0.471959$ with an attached eikonal diffraction fan below $q_0$ certifies this bowtie against all maneuvering intruders down to a unique diffraction threshold $q_a\in(0.457892036651,0.457892036652)$, proving $\Wopt(q)=F(q)$ for all $q_a\le q<1$ (Theorem~\ref{thm:qa-exact}). Weighting the wall-phase budget by vertical displacement further determines the exact constrained width $\Wopt(q,H)$ for all $H\ge 0$ when $1/2\le q<1$, and on an initial interval of caps at every $0<q<1$ (Theorem~\ref{thm:fixed-height}).

\begin{figure}[tb]
\centering
\includegraphics[width=\linewidth]{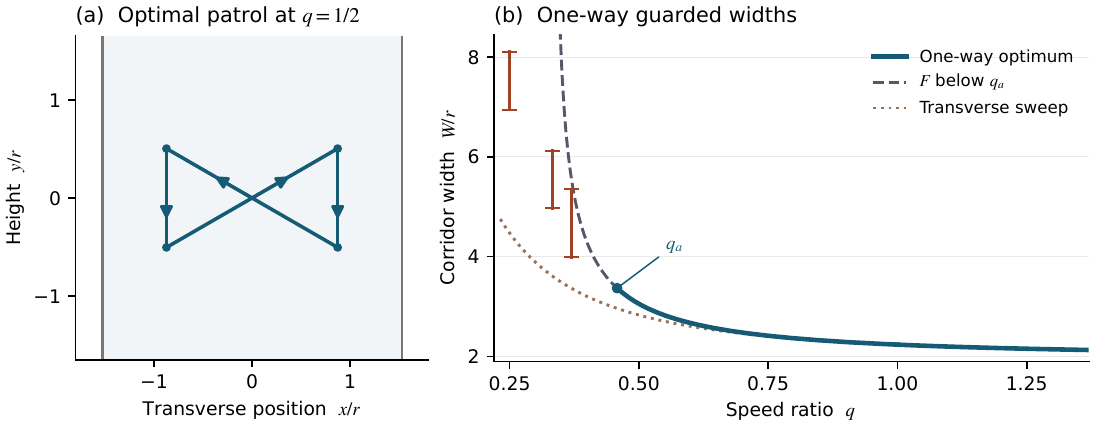}
\caption{Left: the four-leg optimal patrol at $q=1/2$. Both vertical legs point down and both diagonals point up. Right: one-way widths across $q$. The solid curve gives the exact optimum $\Wopt(q)$ for $q\ge q_a$ (attaining $F(q)$ on $[q_a,1)$ and $W_{\mathrm{bf}}(q)$ on $[1,\infty)$). The dashed gray curve gives the wall-runner upper bound $F(q)$ below $q_a$, where the strict gap theorem applies. Vertical bars show the certified lower and upper bounds at $q=1/4,1/3,37/100$.}
\label{fig:landscape}
\end{figure}

\paragraph{Below the diffraction threshold, $0<q<q_a$ (Section~\ref{sec:below-qa}).}
Below $q_a$, an intruder can peel off the wall, round the end of a downward leg, and escape the canonical bowtie with positive clearance (Theorem~\ref{thm:canonical-escape}). Synchronizing long crossings on the phase cylinder upgrades this single-patrol failure to a period- and excursion-uniform strict gap $\Wopt(q)<F(q)$ for all $1/3<q<q_a$, with quadratic scaling $F(q)-\Wopt(q)=\Theta\bigl(r(q_a-q)^2\bigr)$ as $q\uparrow q_a$ (Theorems~\ref{thm:global-gap} and~\ref{thm:global-quadratic}; Proposition~\ref{prop:tilted}). Across restricted patrol classes, $q_a$ is also the sharp threshold where four straight legs cease to be optimal:
\begin{itemize}
    \item Let $W_4(q)$ be the supremum over \emph{ordinary} bowties---symmetric constant-speed four-leg bowties with vertical wall legs, allowing either traversal orientation and transverse offsets. While $W_4(q)=\Wopt(q)$ on $[q_a,\infty)$, an eight-leg patrol strictly beats $W_4(q)$ at every $0<q<q_a$ (Theorem~\ref{thm:ordinary-transition}). As $q\uparrow q_a$, $F(q)-W_4(q)\sim C_4 r(q_a-q)^2$ with $C_4=51.2835202\ldots$, whereas a twelve-leg patrol bending both its wall legs and its crossings lowers the quadratic coefficient to at most $C_{12}=48.4250865\ldots$ (Theorem~\ref{thm:bent-deficit}).
    \item Immediately below $q_a$, allowing arbitrary vertex placement, variable leg speeds, and waits among periodic patrols with at most four straight legs yields no gain over $W_4(q)$ (Theorem~\ref{thm:four}), while adding a single inward vertex on one wall leg gives a five-leg patrol with $C_5=50.8421136\ldots<C_4$ that beats the entire four-leg class on a left neighborhood of $q_a$ (Theorem~\ref{thm:five}).
\end{itemize}
At midrange speeds $q\in\{1/4,1/3,37/100\}$, triangulated space-time surface-cut lower bounds and Galilean lane-flux / wall-parity upper bounds bracket $\Wopt(q)$ quantitatively (Figure~\ref{fig:landscape}, right), and Proposition~\ref{prop:q13-floor} identifies an intrinsic $2\lfloor 1/q\rfloor r$ floor for stateless flux-and-parity multipliers when $q\le 1/3$. The same surface-cut framework reveals an algebraic cost at $q_a$ itself: approximating the curved diffraction fan forces any triangulated space-time certificate within $\varepsilon r$ of $F(q_a)$ to use $N(\varepsilon)=\Theta(\varepsilon^{-1/4})$ triangles across all polygonal patrols (Theorem~\ref{thm:face-rate}).

\paragraph{Slow intruders, $q\downarrow 0$ (Section~\ref{sec:asymptotics}).}
As $q\downarrow 0$, $\Wopt(q)$ grows like Washburn's area bound $2r/q$, and under any fixed excursion cap $H=hr$ satisfies $\lim_{q\downarrow 0}q\Wopt(q,hr)/r=\min\{2,1+h\}$, reaching the full leading coefficient $2$ already at $H=r$ (Theorem~\ref{thm:slow-limit}). For unrestricted bounded excursion, the additive deficit $D_1(q)=2/q-\Wopt(q)/r$ has exact order $\Theta(1)$: ordinary four-leg bowties admit a unique analytic optimizer with $W_4(q)/r=2/q-2+3q-\frac{9}{4}q^3-\frac{19}{4}q^4+O(q^5)$ and limiting deficit $2$ (Theorem~\ref{thm:slow-ordinary}), an explicit ten-leg patrol keeps $D_1(q)<1.235353$ throughout $0<q\le 1/40$ (Theorem~\ref{thm:slow-hooks}), and a reachable-height turn estimate proves $D_1(q)\ge 2^{-21}$ for $0<q\le 2^{-24}$ (Theorem~\ref{thm:loss}).

\paragraph{Bidirectional guarding (Section~\ref{sec:diagnostics}).}
Section~\ref{sec:diagnostics} studies $W_{\mathrm{2way}}(q)$, the supremal width when a single patrol must block upward and downward transits simultaneously. Transverse sweeping is globally optimal already for all $q\ge 1/2$ (Theorem~\ref{thm:twoway}), so at $q=1/2$ knowing the intruder's direction widens the corridor by $\Wopt(1/2)/W_{\mathrm{2way}}(1/2)=\sqrt{7/6}\approx 1.0801$; just below $1/2$, an antipodal eight-leg patrol guards width at least $W_{\mathrm{bf}}(q)+(4\sqrt{2}-o(1))\,r(1/2-q)^2$ (Theorem~\ref{thm:twoway-slope}; Proposition~\ref{prop:twoway-half}). As $q\downarrow 0$, although $W_{\mathrm{2way}}(q)\sim 2r/q$, not knowing the transit direction forfeits $\Wopt(q)-W_{\mathrm{2way}}(q)=\Theta(r/\sqrt{q})$---an unbounded loss in capture radii with $\limsup_{q\downarrow 0}\sqrt{q}\,(\Wopt(q)-W_{\mathrm{2way}}(q))/r\le 2\sqrt{6}$ (Theorem~\ref{thm:direction-value})---and achieving $W_{\mathrm{2way}}(q,H(q))\sim 2r/q$ requires $H(q)/r\to\infty$ (Corollary~\ref{cor:two-way-height}).

\section{Global properties and fast intruders}
\label{sec:global-fast}\label{sec:structural-inputs}

Before examining individual speed intervals below $q=1$, we establish three properties that hold across the entire half-line $q>0$: arbitrary bounded-excursion patrols can be approximated in width by periodic polygonal patrols without increasing excursion; an elementary two-level moving barrier certifies symmetric four-leg bowties (and horizontal sweeps) whenever the downward leg is short enough to reset inside the capture disk; and a midline wall-phase bound determines both the exact $H=0$ width at all $q>0$ and the global optimum $\Wopt(q)=r\sqrt{4+q^{-2}}$ for all fast intruders $q\ge 1$.

\subsection{Periodic polygonal approximation across all speeds}
\label{subsec:periodic-approximation}
Work in units with guard speed one and keep the capture radius \(r>0\)
explicit. Height is measured from the same fixed horizontal axis as in
Definition~\ref{def:wopt}.

\begin{theorem}[Periodic polygonal approximation]\label{thm:periodic-approximation}
Suppose a prescribed patrol of longitudinal excursion at most \(H<\infty\) guards a corridor of
width \(W\) against every upward transit of speed at most \(q>0\).
For every \(0<\lambda<1\), there is a finite-piece periodic polygonal
patrol of excursion at most \(H\), with the same speed bounds and capture radius,
guarding a corridor of width greater than \(\lambda W\).
Consequently, periodic polygonal patrols achieve the same width supremum $\Wopt(q,H)$ at every fixed
\(H\ge 0\), and hence also achieve $\Wopt(q)=\sup_{H\ge 0}\Wopt(q,H)$.
The assertion holds for two-way guarding (Section~\ref{sec:diagnostics}) as well.
\end{theorem}
\begin{proof}
Let \(g\) be the given guard and choose
\(r<s<R_0<R<r/\lambda\). First, increasing the radius to \(s\)
allows a positive increase in intruder speed. To see this, fix \(L>H+R\).
Choose \(c<1\) sufficiently close to one that \(cL>H+r\) and
\((1-c)\sqrt{(W/2)^2+L^2}<s-r\). If a speed-\(q_1=q/c\)
intruder escaped the radius-\(s\) guard, a piece of its path would travel
from \(-L\) to \(L\) within the closed slab between these levels
(take the last visit to \(-L\) before a first subsequent visit to \(L\)).
Scale this piece about the origin by \(c\), keeping its time parameter.
It has speed at most \(q\), moves by less than \(s-r\) at each time, and
stays strictly beyond the radius-\(r\) disk. Its endpoints at \(\pm cL\)
admit safe vertical extensions to a full transit, contradicting guarding.

Choose \(1<A<q_1/q\). Take a time step \(\Delta>0\) and then a spatial
net radius \(\eta>0\) so that
\[
 \Delta/2+\eta<R_0-s,\qquad Aq\Delta/2+\eta<R-R_0,\qquad
 1+2\eta/\Delta\le A,\qquad Aq+2\eta/\Delta\le q_1 .
\]
Use rectangular product nets in the guard box and in
\([-W/2,W/2]\times[-L,L]\). Include the boundary rows and round each
coordinate to a nearest grid coordinate. Denote the finite nets by
\(G,E\) and the bottom and top rows of \(E\) by \(B,C\).
For \(i\ge0\) let \(z_i\in G\) be the rounded \(g(i\Delta)\).
The interpolant of these vertices has speed at most \(A\) and stays
within \(\Delta/2+\eta\) of \(g\). Indeed, at fractional time \(v\)
the deviation of a speed-\(K\) path from the chord between its true
endpoints is at most \(2v(1-v)K\Delta\).

Given consecutive guard vertices \(z,z'\), allow an evader step
\(e\to e'\) in \(E\) if \(|e-e'|\le q_1\Delta\) and the simultaneous
affine interpolants maintain distance strictly greater than \(R_0\)
throughout the step. Write \(\operatorname{Post}_{z,z'}(S)\) for
all destinations from \(S\). Define \(S_0=B\) and
\[
 S_{i+1}= B\cup\operatorname{Post}_{z_i,z_{i+1}}(S_i).
\]
A path reaching \(C\) in any \(S_i\) interpolates to a speed-\(q_1\)
transit avoiding the radius-\(s\) guard. The endpoints admit safe
extensions since \(L>H+s\). Thus every \(S_i\) avoids \(C\).
Some pair \((z_i,S_i)\) repeats because \(G\times 2^E\) is finite.
Repeat the polygon between the repeated states on the whole real line.
The intervening sets now give a periodic invariant containing \(B\)
at every step and disjoint from \(C\), including across the seam.

The periodic polygon \(g_P\) guards at radius \(R\) against speed
\(Aq\). Otherwise take a putative evader piece from \(-L\) to \(L\),
and extend its endpoints constantly to grid times. Round its samples
to \(E\). The rounded steps have length at most \(Aq\Delta+2\eta\),
and their interpolants differ from the evader by at most
\(Aq\Delta/2+\eta<R-R_0\). This produces a forbidden discrete path
from \(B\) to \(C\). Finally let \(\beta=r/R\) and put
\(g_*(t)=\beta g_P(t/(A\beta))\). This scales the two speed bounds
\(A,Aq\) to \(1,q\), the radius \(R\) to \(r\), and the height
to at most \(\beta H\). The guarded width is \(\beta W>\lambda W\).
For two-way transit, include in each finite state both the set
reachable from \(B\) and the set reachable from \(C\), with the
opposite row excluded in each case.
\end{proof}

\subsection{General two-level separator for bowtie and transverse patrols}

Consider a general symmetric four-leg bowtie in $[-b,b]\times\R$ ($b = a + \delta$, $\delta\ge 0$, $a>0$, $h\ge 0$) tracing
\[
    (-a,+h) \xrightarrow{\text{vertical } 2h} (-a,-h) \xrightarrow{\text{diagonal } 2\ell} (+a,+h) \xrightarrow{\text{vertical } 2h} (+a,-h) \xrightarrow{\text{diagonal } 2\ell} (-a,+h)
\]
at unit speed $V=1$, where $\ell := \sqrt{a^2+h^2}$ and period $T = 4h + 4\ell$ (when $h=0$, vertical legs have zero duration and the patrol reduces to horizontal back-and-forth oscillation between $\pm a$). Set
\[
    s := h(1+q), \qquad G := 2q\ell - (1-q)h.
\]

\begin{theorem}[General two-level separator]\label{thm:two-level}
Let $q > 0$, $r > 0$, $a > 0$, $h \ge 0$, and $b = a + \delta$ ($\delta\ge 0$). If
\[
    s = h(1+q) \le r \qquad\text{and}\qquad \delta^2 + G^2 \le r^2,
\]
then the four-leg bowtie $(a,h)$ captures every speed-$q$ upward transit trajectory (including curved, waiting, and reversing paths) in $[-b,b]\times\R$. No sign restriction on $G$ or lower bound on $h/\ell - q$ is required.
\end{theorem}
\begin{proof}
Construct a bounded $T$-periodic moving graph $y = f(t,x)$ spanning $x\in[-b,b]$ according to three barrier rules: (i) on every unlit (non-captured) horizontal flat piece, $f$ rises at speed $\partial_t f = +q$, so along any trajectory $p(t)=(x(t),y(t))$ with $\|p'\|\le q$, $\frac{d}{dt}(y(t)-f(t,x(t))) = y'(t) - q \le 0$, preventing below-to-above crossings; (ii) every connecting chord or vertical step lies entirely inside the current closed capture disk $\overline{B}_r(g(t))$ (by convexity of $\overline{B}_r(g(t))$, it suffices to check segment endpoints); and (iii) instantaneous resets are either upward (expanding the region below $f$) or have their entire swept polygonal region between old and new graphs contained inside $\overline{B}_r(g(t))$.

At the start of the right-hand downward wall leg $(+a,+h)\to(+a,-h)$, set the flat level on $[-b,a]$ to $Y - s$ and join $(a, Y-s)$ to the right-wall point $(b, Y-G)$ by a chord. During the downward leg of duration $2h$ ($Y'=-1$), the flat level rises at speed $+q$, so its vertical offset relative to $Y(t)$ increases linearly from $-s$ to $-s + 2h(1+q) = +s \in [-r,r]$; at every instant, $(a, Y(t)+\text{offset})\in\overline{B}_r(g(t))$ since $|\text{offset}|\le s\le r$, and $(b, Y(t)-G)\in\overline{B}_r(g(t))$ since $\delta^2+G^2\le r^2$.

At departure from $(+a,-h)$ ($Y_0 := -h$, $t=0$), reset the right strip $[a,b]$ downward to the flat level $Y_0 - G$: the swept cut is the triangle with vertices $(a, Y_0+s)$, $(a, Y_0-G)$, and $(b, Y_0-G)$, all three inside $\overline{B}_r(a,Y_0)$. During the leftward diagonal transfer $t\in[0,2\ell]$ ($X(t)=a - (a/\ell)t$, $Y(t)=Y_0 + (h/\ell)t$), place two rising flat levels separated at the guard column $X(t)$:
\[
    F_{\mathrm{ahead}}(t) = Y_0 + s + qt \quad\text{on }[-b, X(t)], \qquad F_{\mathrm{behind}}(t) = Y_0 - G + qt \quad\text{on }[X(t), b],
\]
joined by a vertical step at $x = X(t)$. Since $2\ell(q - h/\ell) = 2q\ell - 2h = G - s$, as $t$ runs from $0$ to $2\ell$, the vertical offset $F_{\mathrm{ahead}}(t) - Y(t)$ changes linearly from $+s$ to $+G$, and $F_{\mathrm{behind}}(t) - Y(t)$ changes linearly from $-G$ to $-s$. Because $|s|\le r$ and $|G|\le r$, both step endpoints at $x=X(t)$ stay inside $\overline{B}_r(g(t))$ throughout $[0,2\ell]$, with $F_{\mathrm{ahead}} - F_{\mathrm{behind}} = s + G = 2q(h+\ell) > 0$.

Upon arrival at $(-a,+h)$ ($Y_1 = +h$), $F_{\mathrm{behind}}$ has offset $-s$ on $[-a,b]$, while $F_{\mathrm{ahead}}$ on $[-b,-a]$ has offset $+G$. Reset $[-b,-a]$ to the chord joining $(-b, Y_1-G)$ to $(-a, Y_1-s)$: the swept region lies inside the convex hull of $(-b, Y_1\pm G)$ and $(-a, Y_1 + \{G, -s\})$, all four inside $\overline{B}_r(-a,Y_1)$ (if $G<0$, any non-illuminated part of the reset is upward and hence safe). Reflecting $x\leftrightarrow -x$ across the second half-period closes a bounded $T$-periodic barrier across $[-b,b]\times\R$ that no uncaptured transit can cross from $-\infty$ to $+\infty$.
\end{proof}

\begin{corollary}[All-speed transverse sufficiency and small-excursion gain]\label{cor:two-level-conseq}
Theorem~\ref{thm:two-level} implies:
\begin{enumerate}
    \item \emph{All-speed transverse sufficiency:} Setting $h=0$ ($s=0, G=2qa$) and maximizing $b = a + \sqrt{r^2 - 4q^2 a^2}$ at $a_0 = r/(2q\sqrt{1+4q^2})$ proves $\Wopt(q,0) \ge W_{\mathrm{bf}}(q) = r\sqrt{4+q^{-2}}$ against all maneuvering paths for every $q > 0$.
    \item \emph{Strict suboptimality of $W_{\mathrm{bf}}(q)$ for all $0<q<1$:} For any $0 < q < 1$ and all sufficiently small $H>0$, setting $h=H$ ($s\le r$) and $a = a_0$ gives $G = R_0 - (1-q)H + O(H^2/r)$ ($R_0 = r/\sqrt{1+4q^2}$, $\delta_0 = 2q R_0$), hence $\sqrt{r^2 - G^2} = \delta_0 + \frac{1-q}{2q}H - O(H^2/r)$ and $\Wopt(q,H) \ge W_{\mathrm{bf}}(q) + \frac{1-q}{q}H - O(H^2/r) > W_{\mathrm{bf}}(q)$. Every larger positive cap allows a sufficiently small successful stroke as well.
\end{enumerate}
\end{corollary}

\subsection{Global midline-excursion bound and the fast-intruder optimum}
\label{sec:fast}

Pairing Corollary~\ref{cor:two-level-conseq} with a global wall-phase bound across midline excursions determines $\Wopt(q,0)$ for all $q>0$, solves $\Wopt(q)$ completely for fast intruders $q\ge 1$, and pins down the exact first derivative of $\Wopt(q,H)$ at $H=0$.

\begin{theorem}[Global midline-excursion bound and fast-intruder optimum]\label{thm:midline}
Let $q > 0$, and let $g=(X,Y):\R\to[-b,b]\times\R$ be any prescribed unit-speed patrol (periodic or nonperiodic) satisfying $\sup_{t\in\R}|Y(t)|\le H < \infty$. If $g$ successfully guards $[-b,b]\times\R$, then
\begin{equation}\label{eq:midline-bound}
    W = 2b \le r\sqrt{4 + q^{-2}} + \frac{(1-q)_+}{q}\,H, \qquad (1-q)_+ := \max(0,\, 1-q).
\end{equation}
Consequently:
\begin{enumerate}
    \item For every $q>0$, $\Wopt(q,0)=r\sqrt{4+q^{-2}}$; and for all fast intruders $q \ge 1$, the exact unrestricted minimax corridor width is
    \begin{equation}\label{eq:fast-opt}
        \Wopt(q) = r\sqrt{4 + q^{-2}}, \qquad \forall\, q \ge 1.
    \end{equation}
    \item For every $0 < q < 1$, $\Wopt(q,H) = r\sqrt{4+q^{-2}} + \frac{1-q}{q}H + O(H^2/r)$ as $H\downarrow 0$.
\end{enumerate}
\end{theorem}
\begin{proof}
Normalize $r=1, V=1, U=q$, and define the transverse reference parameters
\[
    a_0 = \frac{1}{2q\sqrt{1+4q^2}}, \quad b_0 = \frac{\sqrt{1+4q^2}}{2q} = \frac{1}{2}\sqrt{4+q^{-2}}, \quad \delta_0 = b_0 - a_0 = \frac{2q}{\sqrt{1+4q^2}}, \quad R_0 = \frac{1}{\sqrt{1+4q^2}}.
\]
Note $R_0 = 2qa_0$ and $\delta_0 / R_0 = 2q$. Assume $b > b_0 > 1$ and set $\Delta := b - b_0 > 0$.

At the right wall $x=+b$, detection at time $t$ requires $b - X(t) \le 1$, hence $X(t) \ge b - 1 > b_0 - 1 > 0$. The detected straight-run wall phases at time $t$ form $[qt - Y(t) - R(b-X(t)),\, qt - Y(t) + R(b-X(t))]$ where $R(d)=\sqrt{1-d^2}$. By concavity of $R$ on $[0,1]$ and its supporting tangent at $\delta_0$ ($R'(\delta_0) = -\delta_0/R_0 = -2q$), every active right-wall contact satisfies
\[
    R(b - X) \le R_0 - 2q\bigl[(b - X) - \delta_0\bigr] = 2q b_0 - 2q(b - X) = 2q(X - \Delta).
\]
Consider any complete positive midline excursion of $X$ on $[0,D]$ ($X(0)=X(D)=0$ and $X(t)>0$ on $(0,D)$). Put $A(t) := qt - Y(t) + 2q X(t)$ and $B(t) := qt - Y(t) - 2q X(t)$. The total span length $I$ of right-wall phases detected during $[0,D]$ satisfies $I \le \max_{[0,D]} A - \min_{[0,D]} B - 4q\Delta$. Choose $u,v\in[0,D]$ attaining $\max A$ and $\min B$, and put $d := |u - v|$. Since $|X'|\le 1$ and $X(0)=X(D)=0$, $X(u)\le\min(u,D-u)$ and $X(v)\le\min(v,D-v)$ imply $X(u)+X(v)\le D - d$, while $|Y|\le H$ and $|Y'|\le 1$ give $Y(v)-Y(u)\le\min(2H,d)$. Therefore
\[
    A(u) - B(v) - 2qD = q(u-v) + \bigl(Y(v)-Y(u)\bigr) + 2q\bigl(X(u)+X(v)-D\bigr) \le
    \begin{cases}
        \min(2H,d) - qd, & u \ge v,\\
        \min(2H,d) - 3qd, & u < v.
    \end{cases}
\]
Maximizing over $d\ge 0$ yields $A(u)-B(v)-2qD \le 2H(1-q)_+ =: C$. By symmetry ($X\mapsto -X$), every complete active excursion on either side of the midline covers wall-phase span at most
\begin{equation}\label{eq:excursion-span}
    I \le 2q D + C - 4q\Delta.
\end{equation}
On any window $[0,L]$, at most two active excursions are truncated by $0$ and $L$, each contributing span at most $2q D_{\mathrm{tr}} + C + 4qb$. Every complete active excursion reaches $|X|\ge b-1>0$ from $0$ and returns in time $\ge 2(b-1)>0$, so their count $N(L)$ is finite with $\sum D_k \le L$. Covering $[H+2,\, qL - H - 2]$ on both walls requires $2qL - 4(H+2) \le 2qL + N(L)(C - 4q\Delta) + 2(C + 4qb)$. Because late wall phases $z\to+\infty$ on both walls require infinitely many sign switches, $N(L)\to\infty$ as $L\to\infty$, forcing $C - 4q\Delta \ge 0$, i.e., $2\Delta \le \frac{(1-q)_+}{q}H$, proving~\eqref{eq:midline-bound}. Setting $H=0$ in~\eqref{eq:midline-bound} and combining with Corollary~\ref{cor:two-level-conseq}(1) gives $\Wopt(q,0)=r\sqrt{4+q^{-2}}$ for all $q>0$; for $q\ge 1$, $(1-q)_+=0$ gives $\Wopt(q,H)\le r\sqrt{4+q^{-2}}$ for every $H$, proving~\eqref{eq:fast-opt}. Part~(2) follows by combining~\eqref{eq:midline-bound} with Corollary~\ref{cor:two-level-conseq}(2) as $H\downarrow0$.
\end{proof}

\section{Above the diffraction threshold (\texorpdfstring{$q_a \le q < 1$}{q\_a <= q < 1})}
\label{sec:above-qa}

When $q<1$, Theorem~\ref{thm:midline}(2) shows that longitudinal motion strictly improves on horizontal sweeping. For $1/3<q<1$, intruders running straight along the two walls impose an explicit upper bound $F(q)$ that merges continuously into $\sqrt{5}\,r=W_{\mathrm{bf}}(1)$ as $q\uparrow 1$. We establish this bound first, then certify that the symmetric four-leg bowtie saturating the wall budget blocks all maneuvering intruders down to $q_a\approx 0.457892036651$, and finally determine exact excursion-constrained widths $\Wopt(q,H)$.

\subsection{The wall-runner upper bound for \texorpdfstring{$1/3 < q < 1$}{1/3 < q < 1}}
\label{sec:wall-upper}

In this section normalize $V=1$, take $1/3<q<1$, and set
\begin{equation}\label{eq:F-def}
    F(q) := 2r\sqrt{1 + \frac{1}{3q^2+2q-1}}.
\end{equation}
Notice $3q^2+2q-1=(3q-1)(q+1)>0$ precisely when $q>1/3$, and $F(q)>2r$ for $q\in(1/3,1)$.

\begin{theorem}[Wall-runner upper bound]\label{thm:wall-upper}
Let $1/3 < q < 1$, and let $g=(X,Y):\R\to[-b,b]\times\R$ be any unit-speed patrol with bounded longitudinal excursion $\sup_{t\in\R}|Y(t)|\le H < \infty$ (including every periodic patrol). If $g$ blocks all upward straight-line wall runners $p_{\pm,\zeta}(t) = (\pm b,\, qt - \zeta)$ ($\zeta\in\R$), then $W = 2b \le F(q)$.
\end{theorem}

Rather than summing potentially overlapping phase spans from arbitrary visits, we map the patrol onto a Euclidean phase cylinder, uncross self-intersections without increasing Euclidean length, and extract nondecreasing lifted contact order by Jordan separation.

\begin{lemma}[Phase-cylinder speed budget and inverse parameterization]\label{lem:phase-cylinder}
Fix $q\in(0,1)$, $\alpha:=\sqrt{1-q^2}$, half-width $b>0$, and period $T>0$ with $P:=qT$.
\begin{enumerate}
    \item For any $T$-periodic unit-speed patrol $g=(X,Y)$ in $[-b,b]\times\R$, the curve $\Gamma(t) := \bigl(\alpha X(t),\, qt - Y(t)\bigr)$ projects to a winding-number-$+1$ loop on the cylinder $[-\alpha b,\alpha b]\times(\R/P\Z)$ of Euclidean length $L(\Gamma) \le T = P/q$, and $g$ detects the straight wall runner $(\pm b, qt-\zeta)$ iff $\Gamma$ enters the elliptical disk
    \[
        \frac{(u \mp \alpha b)^2}{\alpha^2} + \min_{k\in\Z}(z-\zeta-kP)^2\le r^2.
    \]
    \item Conversely, every rectifiable winding-$+1$ loop $\Gamma_0(s) = (u(s), z(s))$ ($s\in[0,L]$ parameterized by Euclidean arclength, $z(L)=z(0)+P$, $|u(s)|\le \alpha b$) with length $L \le P/q$ lifts to a $T$-periodic unit-speed patrol $g_0$ ($T=P/q$) tracing the exact same cylinder image.
\end{enumerate}
\end{lemma}
\begin{proof}
For~(1), since $X'^2+Y'^2\le 1$ a.e.\ and $1-qY' \ge 1-q > 0$, the pointwise algebraic identity
\[
    (1 - qY')^2 - \bigl[\alpha^2 X'^2 + (q - Y')^2\bigr] = (1-q^2)\bigl(1 - X'^2 - Y'^2\bigr) \ge 0
\]
yields $\|\Gamma'(t)\|_2 = \sqrt{\alpha^2 X'^2 + (q-Y')^2} \le 1 - qY'(t)$ a.e. Integrating over $[0,T]$ gives $L(\Gamma) \le T - q(Y(T)-Y(0)) = T = P/q$, and $(X(t)\mp b)^2 + (qt-\zeta-Y(t))^2 \le r^2$ translates directly into the stated ellipse condition.

For~(2), since $u'(s)^2 + z'(s)^2 = 1$ a.e.\ and $q\in(0,1)$, the minimum time-rate density $\tau_{\min}(s) := \frac{1 - q z'(s)}{1 - q^2} \ge \frac{1-q}{1-q^2} = \frac{1}{1+q} > 0$ satisfies $\int_0^L \tau_{\min}(s)\,ds = \frac{L - qP}{1-q^2} \le \frac{P/q - qP}{1-q^2} = \frac{P}{q} = T$. Setting $t'(s) := \tau_{\min}(s) + c_0 > 0$ with constant $c_0 := \frac{T - \int_0^L \tau_{\min}}{L} \ge 0$ defines a bi-Lipschitz homeomorphism $t:[0,L]\to[0,T]$. Setting $X(t(s)) := u(s)/\alpha$ and $Y(t(s)) := qt(s) - z(s)$ gives $Y(T)=Y(0)$, $X'(t) = \frac{u'(s)}{\alpha t'(s)}$, and $Y'(t) = q - \frac{z'(s)}{t'(s)}$. Using $u'(s)^2 = 1 - z'(s)^2$ and $t'(s) = \tau_{\min}(s) + c_0$, direct expansion yields
\[
    \alpha^2 t'(s)^2\bigl(1 - X'(t)^2 - Y'(t)^2\bigr) = \alpha^2 c_0(2+\alpha^2 c_0) \ge 0
\]
almost everywhere, since $c_0\ge 0$ and $\alpha^2>0$. Thus $X'(t)^2 + Y'(t)^2 \le 1$ a.e.
\end{proof}

\begin{lemma}[Cylindrical uncrossing with arbitrarily small coverage loss]\label{lem:uncrossing}
Let $\Gamma$ be a rectifiable loop of winding number $+1$ and Euclidean length $L$ in $[-\alpha b, \alpha b]\times(\R/P\Z)$ with $b > r$, whose elliptical neighborhoods of radius $r$ cover both wall circles $\{\pm \alpha b\}\times(\R/P\Z)$. Then for every $\varepsilon > 0$, there exists a \emph{simple} (embedded) winding-$+1$ loop $\widetilde{\Gamma}$ in $[-\alpha b,\alpha b]\times(\R/P\Z)$ of Euclidean length at most $L$ whose elliptical neighborhoods of radius $r+\varepsilon$ cover both wall circles. Moreover $\widetilde{\Gamma}$ can be chosen so that the images of $\Gamma$ and $\widetilde{\Gamma}$ are within Hausdorff distance $\varepsilon$ of each other in the detection metric $d\bigl((u,z),(u',z')\bigr)^2=(u-u')^2/\alpha^2+\min_{k\in\Z}(z-z'-kP)^2$ of Lemma~\ref{lem:phase-cylinder}(1).
\end{lemma}
\begin{proof}
First compress horizontally via $u \mapsto (1-\eta)u$ for small $\eta\in(0,1)$. This moves points by at most $\eta b$ in the physical detection metric of Lemma~\ref{lem:phase-cylinder}(1), leaves a positive boundary margin $\eta\alpha b$ inside $[-\alpha b,\alpha b]$, and \emph{strictly} reduces Euclidean length: since $b > r$, covering both walls forces positive horizontal variation $\int |u'|\,ds > 0$, so the line element strictly shrinks on a set of positive measure. Next, approximate the compressed curve by an inscribed polygon and then perturb its vertices within the strict length and boundary margins. This gives a generic polygonal immersion with the same winding number $+1$, transverse double crossings only (no triple points or overlapping edges), length strictly less than $L$, image inside $[-\alpha b,\alpha b]$, and Hausdorff error small enough that $\eta b$ plus polygonal error is less than $\varepsilon/2$.

Place pairwise disjoint disks of radius $\delta$ with $2\delta/\alpha < \varepsilon/2$ inside $(-\alpha b,\alpha b)\times(\R/P\Z)$ around the finitely many double points, avoiding all vertices and non-incident edges. Inside each disk, replacing the four radial half-edges by either of the two noncrossing adjacent chord pairs does not increase Euclidean length (by the triangle inequality) and keeps every replaced point within $2\delta$ of a new chord. Any choice of noncrossing pairings decomposes the graph into finitely many pairwise disjoint simple closed curves. If more than one component results, connectedness of the original immersion implies that two distinct components meet at some common crossing disk; switching the noncrossing pairing inside that single disk merges those two components into one without introducing intersections (compare the A-trail formulation in~\cite{AndersenBouchetJackson1996}). Iterating yields a single simple closed curve. Because local smoothing preserves mod-$2$ homology on the cylinder, the resulting simple loop has odd winding number; since every simple closed curve on a cylinder has winding number in $\{0,\pm 1\}$, orienting it upward gives winding number $+1$, length $\le L$, and wall coverage at radius $r+\varepsilon$.

For the last assertion, each step changes the image by a small Hausdorff distance in the detection metric. The compression moves every point by at most $\eta b$. The inscribed polygon and the perturbation of its vertices contribute the polygonal error. The old and the new curve agree outside the disks, and both meet every disk, so the changes inside the disks contribute at most $2\delta/\alpha$. By the choices above the sum is less than $\varepsilon$.
\end{proof}

\begin{proof}[Proof of Theorem~\ref{thm:wall-upper}]
Assume $b > r$ (otherwise $2b \le 2r < F(q)$).

\smallskip\noindent\textbf{Step 1: Periodic patrols.} Suppose $g$ is $T$-periodic. Fix $0 < \varepsilon < b - r$. Applying Lemma~\ref{lem:phase-cylinder}(1), Lemma~\ref{lem:uncrossing}, and Lemma~\ref{lem:phase-cylinder}(2) yields a $T$-periodic unit-speed patrol $\tilde{g} = (\widetilde{X},\widetilde{Y})$ in $[-b,b]\times\R$ whose phase-cylinder curve $\widetilde{\Gamma}$ is \emph{simple} and covers both wall circles with capture radius $r_\varepsilon := r + \varepsilon < b$. Lift $\widetilde{\Gamma}$ to an infinite simple curve $\widehat{\Gamma}$ in $[-\alpha b,\alpha b]\times\R$. For each lifted phase $z\in\R$, let $q_+(z)$ be the rightmost point of $\widehat{\Gamma}\cap([-\alpha b,\alpha b]\times\{z\})$ and $q_-(z)$ the leftmost; retaining $q_\pm(z)$ preserves all wall coverage because horizontal distance to $\pm\alpha b$ only decreases at fixed phase $z$. Along the oriented simple curve $\widehat{\Gamma}$, the selected points $\{q_+(z)\}$ (and likewise $\{q_-(z)\}$) appear in \emph{nondecreasing lifted phase order}: horizontal segments from $q_+(z_1), q_+(z_2)$ ($z_1 < z_2$) to a vertical line $\{u = \alpha b + 1\}$ are disjoint crosscuts of the strip to the right of $\widehat{\Gamma}$, so cutting along one crosscut forms a topological rectangle whose opposite sides must share endpoint order by Jordan separation. More explicitly, truncate the lift between two horizontal lines far below and above the chosen phase levels, taking the last visit to the bottom line before the two contacts and the first visit to the top line after them. Properness and the limits $z(t)\to\pm\infty$ at the two ends allow the bottom line to lie below all phases after the earlier contact and the top line above all phases before the later contact. Thus the truncated arc stays between the two lines. The resulting simple crosscut joins the bottom to the top of a rectangle and the rightward horizontal connectors miss the crosscut except at their endpoints. Disjoint connectors preserve the boundary order along its right-hand component. The left wall is identical. Passing to closures preserves this weak order.

Take the closure of the rightmost selections satisfying $\widetilde X\ge b-r_\varepsilon$, and the closure of the leftmost selections satisfying $\widetilde X\le -b+r_\varepsilon$; because $b > r_\varepsilon$, the left and right zones are separated by horizontal gap $2(b - r_\varepsilon) > 0$. By rectifiability, one period contains finitely many alternating left- and right-wall contact groups ($2N$ groups total for some integer $N\ge 1$). In group $j\in\{1,\dots,2N\}$, let $u_j, v_j$ be points of the group attaining the maximum of $p + R_\varepsilon$ and minimum of $p - R_\varepsilon$ respectively, where $p = qt - \widetilde{Y}$ is the center phase and $R_\varepsilon(\widetilde{X}) = \sqrt{r_\varepsilon^2 - (b - |\widetilde{X}|)^2}$. We claim $v_j$ can always be chosen no later along $\widetilde{\Gamma}$ than $u_j$: if $u_j$ precedes $v_j$, nondecreasing lifted phase order gives $p(u_j)\le p(v_j)$, while extremality gives $p(v_j)+R_\varepsilon(v_j)\le p(u_j)+R_\varepsilon(u_j)$ and $p(u_j)-R_\varepsilon(u_j)\ge p(v_j)-R_\varepsilon(v_j)$, i.e., $R_\varepsilon(u_j)-R_\varepsilon(v_j)\ge p(v_j)-p(u_j)\ge 0$ and $R_\varepsilon(v_j)-R_\varepsilon(u_j)\ge p(v_j)-p(u_j)\ge 0$. Hence $p(u_j)=p(v_j)$ and $R_\varepsilon(u_j)=R_\varepsilon(v_j)$, so swapping $(u_j,v_j)$ puts the minimizer before the maximizer.

Thus each group $j\in\{1,\dots,2N\}$ has chronologically ordered marker times $t_j^- \le t_j^+$ covering phase length $q(t_j^+-t_j^-) + (\widetilde{Y}(t_j^-)-\widetilde{Y}(t_j^+)) + R_\varepsilon(t_j^-) + R_\varepsilon(t_j^+)$. Let $A := \sum_{j=1}^{2N}(t_j^+-t_j^-) \ge 0$, $B := T - A = \sum_{j=1}^{2N}(t_{j+1}^- - t_j^+) \ge 0$ (with $t_{2N+1}^- := t_1^-+T$), $Z := \sum_{j=1}^{2N}(\widetilde{Y}(t_j^-)-\widetilde{Y}(t_j^+))$, and $a := \frac{1}{4N}\sum_{j=1}^{2N}(|\widetilde{X}(t_j^-)| + |\widetilde{X}(t_j^+)|) \in [b-r_\varepsilon, b]$. Concavity of $d\mapsto\sqrt{r_\varepsilon^2 - d^2}$ bounds the sum of the $4N$ marker radii by $4NR$, where $R := \sqrt{r_\varepsilon^2 - (b-a)^2}$. On each inter-group transit $[t_j^+, t_{j+1}^-]$, $\widetilde{X}(t_j^+)$ and $\widetilde{X}(t_{j+1}^-)$ have opposite signs, so horizontal travel is $|\widetilde{X}(t_j^+)| + |\widetilde{X}(t_{j+1}^-)|$, summing over $j=1,\dots,2N$ to $4Na$, while vertical shifts telescope to $\sum_{j=1}^{2N}(\widetilde{Y}(t_{j+1}^-)-\widetilde{Y}(t_j^+)) = Z$. By the Euclidean triangle inequality over $[t_j^+,t_{j+1}^-]$ and unit speed on $[t_j^-,t_j^+]$, full coverage of both length-$qT$ wall circles implies
\begin{equation}\label{eq:wall-group-bounds}
    2qT \le qA + Z + 4NR, \qquad A \ge |Z|, \qquad T - A = B \ge \sqrt{(4Na)^2 + Z^2}.
\end{equation}
Since $2qT - qA - Z = qA + 2qB - Z \ge q|Z| + 2q\sqrt{(4Na)^2+Z^2} - Z$, setting $h := Z/(4N) \in \R$ and dividing by $4N$ yields
\[
    2q\sqrt{a^2+h^2} + q|h| - h \le R.
\]
If $h < 0$, replacing $h$ by $0$ strictly decreases the left-hand side; for $h\ge 0$, $|h|=h$ gives $\Phi_q(a,h) := 2q\sqrt{a^2+h^2} - (1-q)h \le R$. Because $0 < 1-q < 2q$ when $1/3 < q < 1$, $\Phi_q(a,\cdot)$ is strictly convex on $[0,\infty)$ and attains its global minimum at
\begin{equation}\label{eq:opt-h}
    h_*(a) = \frac{(1-q)a}{\sqrt{3q^2+2q-1}}, \qquad \min_{h\ge 0}\Phi_q(a,h) = a\sqrt{3q^2+2q-1}.
\end{equation}
Writing $c := 3q^2+2q-1 > 0$, $a\sqrt{c} \le \sqrt{r_\varepsilon^2 - (b-a)^2}$ gives $c a^2 + (b-a)^2 \le r_\varepsilon^2$. By Cauchy--Schwarz on $(1/\sqrt{c},\,1)\cdot(\sqrt{c}\,a,\,b-a)$, $b = a + (b-a) \le r_\varepsilon\sqrt{1+c^{-1}}$. Letting $\varepsilon\downarrow 0$ proves $2b \le 2r\sqrt{1+c^{-1}} = F(q)$.

\smallskip\noindent\textbf{Step 2: Nonperiodic bounded-excursion patrols.} Suppose $\sup_{t\in\R}|Y(t)|\le H$ without periodicity. Put $K := H + r$ and examine $[0,L]$ for large $L$ with $P := qL > 2K$. Any wall phase $\zeta\in(K, P-K)$ satisfies $|qt-\zeta| > H+r$ for all $t\notin[0,L]$ (since $|Y|\le H$), so $\Gamma|_{[0,L]}$ (of Euclidean length $\le L + q|Y(L)-Y(0)| \le L + 2qH$) covers both walls across $(K, P-K)$. Close $\Gamma|_{[0,L]}$ on the circumference-$P$ cylinder by joining $(\alpha X(L), P-Y(L))$ to $(\alpha X(0), P-Y(0))$ (length $\le 2\sqrt{\alpha^2 b^2 + H^2}$), and attach at its base point two contractible polygonal detours inside $[-\alpha b,\alpha b]\times(\R/P\Z)$ visiting the vertical seam segments $\{\pm\alpha b\}\times[-K,K]$ and returning (total length $\le 8\alpha b + 8K + 4H$). The resulting closed loop has winding number $+1$, full wall coverage at radius $r$, and Euclidean length $L_{\mathrm{geom}} \le L + C_0$, where $C_0$ depends only on $(q,b,H,r)$.

Set $q_L := P / \max(L, L_{\mathrm{geom}}) \le q$ and $\alpha_L := \sqrt{1-q_L^2} \ge \alpha$; then $q_L\to q$ as $L\to\infty$. Read the exact same cylinder loop as data for speed ratio $q_L$, half-width $b_L := \alpha b/\alpha_L$, and capture radius $r$: since $\alpha_L\ge\alpha$, the old $\alpha$-ellipses fit inside the new $\alpha_L$-ellipses, so full wall coverage persists and $L_{\mathrm{geom}}\le P/q_L$. Lemma~\ref{lem:phase-cylinder}(2) supplies a periodic unit-speed patrol at $(q_L, b_L)$; for all large $L$, $q_L > 1/3$, so Step~1 yields $2(\alpha b/\alpha_L) \le F(q_L)$. Letting $L\to\infty$ proves $2b\le F(q)$.
\end{proof}

\subsection{Wall-optimal bowties and the obstacle-wavefront criterion}
\label{sec:wavefront}

For $q\in(1/3,1)$ and $c := 3q^2+2q-1 > 0$, define the wall-optimal symmetric four-leg bowtie parameters (normalized to $r=1, V=1$) by
\begin{equation}\label{eq:bowtie-params}
    b(q) = \sqrt{1 + \frac{1}{c}}, \quad a(q) = \frac{1}{\sqrt{c(1+c)}}, \quad h(q) = \frac{1-q}{c\sqrt{1+c}}, \quad R(q) = \frac{1}{\sqrt{1+c}}, \quad \ell(q) = \frac{2q}{c\sqrt{1+c}}.
\end{equation}
Notice that as $q\uparrow 1$, $c\to 4$, $F(q)=2b(q)r\to\sqrt{5}\,r=W_{\mathrm{bf}}(1)$, and $h(q)\to 0$, merging~\eqref{eq:bowtie-params} continuously into Washburn's horizontal sweep~\eqref{eq:washburn-bf}. At every $q\in(1/3,1)$,~\eqref{eq:bowtie-params} attains the upper bound $F(q)$ against straight wall runners alone: on each downward wall leg the vertical disk chord has half-length $R(q)$ while its center phase advances by $2(1+q)h(q)$, sweeping total phase $2R(q)+2(1+q)h(q)=4q(h(q)+\ell(q))=qT$, an entire period.

For $q$ above $q_0\approx 0.4719588902$---the unique root of $27q^4-10q^2+4q-1=0$ in $(1/3,1)$---the two-level separator of Theorem~\ref{thm:two-level} already certifies~\eqref{eq:bowtie-params} against all maneuvering intruders: we have $G = 2q\ell(q)-(1-q)h(q) = R(q)$, $\delta = b(q)-a(q) = \sqrt{c/(1+c)}$, $\delta^2+G^2 = 1$, and $s(q) = (1+q)h(q) = \frac{1-q^2}{c\sqrt{1+c}} \le 1$ if and only if
\[
    c^2(1+c) - (1-q^2)^2 = (1+q)^2\bigl[q(3q-1)^2(3q+2) - (1-q)^2\bigr] = (1+q)^2\bigl(27q^4 - 10q^2 + 4q - 1\bigr) \ge 0.
\]
In particular, since $27(59/125)^4 - 10(59/125)^2 + 4(59/125) - 1 = \frac{59372}{244140625} > 0$, Theorem~\ref{thm:two-level} proves $\Wopt(q)=F(q)$ throughout $[59/125,1)=[0.472,1)$ (and throughout $[q_0,1)$).

Below $q_0$, $h(q)(1+q) > 1$ on~\eqref{eq:bowtie-params}, so a flat horizontal front resets outside $\overline{B}_1$. Instead, we attach a curved eikonal diffraction fan directly to the descending disk boundary and hand it off into the convex space-time capture tube of the succeeding diagonal leg. For any general bowtie $(a,h)$ ($r=1, V=1$, $\ell=\sqrt{a^2+h^2}$, diagonal leg length $L_{\mathrm{diag}} := 2\ell$), set
\begin{equation}\label{eq:diffraction-defs}
    \alpha = \sqrt{1-q^2}, \qquad c_{\mathrm{leg}} = \frac{a}{\ell}, \qquad d_{\mathrm{leg}} = \frac{h}{\ell}, \qquad m = q d_{\mathrm{leg}} + \alpha c_{\mathrm{leg}}, \qquad \kappa = \alpha d_{\mathrm{leg}} - q c_{\mathrm{leg}}.
\end{equation}
Notice that $0 < q < d_{\mathrm{leg}} < 1 \iff h^2 > q^2(a^2+h^2) \iff \alpha d_{\mathrm{leg}} > q c_{\mathrm{leg}} \iff \kappa > 0$. Define the terminal diffraction angle $\theta_*$ and elliptic-type fan integral $p$ by
\begin{equation}\label{eq:theta-p}
    \theta_* := \arctan\!\left(\frac{q\kappa}{1 + q m}\right), \qquad p := \int_0^{\theta_*} \frac{d\theta}{\cos\theta + \sqrt{q^2 - \sin^2\theta}},
\end{equation}
and the \emph{rounded-escape ratio}
\begin{equation}\label{eq:R-round}
    \mathcal{R}_{\mathrm{round}}(a,h;q) := \frac{h(m + q)}{1 + \kappa\cos\theta_* - m\sin\theta_* + (m+q)p}.
\end{equation}

\begin{lemma}[Continuous attached diffraction fan]\label{lem:attached-fan}
Suppose $0 < q < d_{\mathrm{leg}} < 1$ ($\kappa > 0$), and let the guard move downward along $(a, -t)$ for $t\in[0,p]$. Define $\theta(s)\in[0,\theta_*]$ implicitly for $s\in[0,p]$ by $\int_0^{\theta(s)}\frac{d\theta}{\cos\theta+\sqrt{q^2-\sin^2\theta}} = s$, and set $X_{\mathrm{fan}}(s) := (a - \cos\theta(s),\, \sin\theta(s) - s)$ with velocity $X_{\mathrm{fan}}'(s) = (u(s), w(s))$. Then $u^2+w^2=q^2$, $u'>0$, $w\ge qm > 0$, $(u(0),w(0))=(0,q)$, $(u(p),w(p))=q(\kappa,m)$, and for each $t\in[0,p]$ the tangent family
\[
    P_s(t) := X_{\mathrm{fan}}(s) + (t - s) X_{\mathrm{fan}}'(s), \qquad 0 \le s \le t,
\]
forms a spatial graph $y=f(t,x)$, smooth on its unlit interior with continuous first derivatives up to the contact endpoint, over $x\in[a-1,\, a-\cos\theta(t)]$ satisfying the eikonal identity $f_t = q\sqrt{1+f_x^2}$ outside $\overline{B}_1(a,-t)$, joining continuously to the flat front $y=qt$ at $x=a-1$ and ending at $X_{\mathrm{fan}}(t)\in\partial B_1(a,-t)$.
\end{lemma}
\begin{proof}
Differentiating $\theta'(s) = \cos\theta + \sqrt{q^2-\sin^2\theta} > 0$ gives $u = \sin\theta\,\theta'$, $w = \cos\theta\,\theta' - 1$, and $u^2+w^2 = \theta'^2 - 2\cos\theta\,\theta' + 1 = q^2$. Writing $(u,w)=(q\sin\psi, q\cos\psi)$ yields $\tan\theta = \frac{q\sin\psi}{1+q\cos\psi}$ with $\frac{d\theta}{d\psi} = \frac{q(q+\cos\psi)}{1+2q\cos\psi+q^2} > 0$; as $\psi$ increases from $0$ to $\arcsin\kappa = \arccos m \in (0,\pi/2)$ (note $\kappa^2+m^2=1$), $\theta$ increases strictly from $0$ to $\theta_*$, with $L := \sqrt{1+2q\cos\psi+q^2}$ matching $\cos\theta + \sqrt{q^2-\sin^2\theta}$ since $L - \cos\theta = \frac{q^2+q\cos\psi}{L} > 0$. Thus $w \ge qm > 0$ and $u' > 0$. For $0\le s < t$, $\partial_s P_s(t) = (t-s)u'(s)(1, -u(s)/w(s))$ has strictly positive $x$-component and slope $f_x = -u/w$, so $f_t = w - f_x u = (w^2+u^2)/w = q^2/w = q\sqrt{1+f_x^2}$. Relative to the current guard $(a,-t)$, $P_s(t) - (a,-t) = (-\cos\theta(s),\sin\theta(s)) + (t-s)(u(s), w(s)+1)$; these two vectors are orthogonal ($\theta' = \cos\theta + \sqrt{q^2-\sin^2\theta} \implies -\cos\theta\,u + \sin\theta(w+1) = 0$), so $\|P_s(t)-(a,-t)\|_2 > 1$ for $s<t$. By Cauchy--Schwarz and $f_t = q\sqrt{1+f_x^2}$, no speed-$q$ trajectory can cross $y=f(t,x)$ upward outside the captured endpoint $X_{\mathrm{fan}}(t)$.
\end{proof}

\begin{theorem}[Obstacle-wavefront sufficient criterion]\label{thm:obstacle-criterion}
Normalize $r=1, V=1$. Suppose a four-leg bowtie $(a,h)$ ($b = a+\delta$, $\ell=\sqrt{a^2+h^2}$, $L_{\mathrm{diag}}=2\ell$, $G = 2q\ell - (1-q)h$) satisfies $0 < q < d_{\mathrm{leg}} < 1$, $\delta\ge 0$, $0 < G < 1$, and $\delta^2+G^2 \le 1$. Put $C := \cos\theta_*$, $S := \sin\theta_*$, $K := \kappa C - mS + (m+q)p$, $v := h - p$, and for $j\in\{0,1\}$ define $(x_0,y_0) := (-1, h+qp)$, $(x_1,y_1) := (-C, h-p+S)$, and
\[
    N_j := d_{\mathrm{leg}} x_j + c_{\mathrm{leg}} y_j, \qquad U_j := \frac{-m x_j + \kappa y_j}{\alpha}.
\]
If all four condition groups hold:
\begin{align}
    h(m + q) &\le 1 + K \qquad\bigl(\text{i.e., }\mathcal{R}_{\mathrm{round}}(a,h;q) \le 1\bigr), \label{eq:obs-3}\\
    0 < p &< h, \qquad a + b \ge 1, \label{eq:obs-4}\\
    q v &\le 1, \qquad \alpha^2 v \le 2(G - q), \label{eq:obs-5}\\
    N_j &\ge 0, \qquad 0 \le U_j \le L_{\mathrm{diag}} \qquad (j = 0, 1), \label{eq:obs-6}
\end{align}
then $(a,h)$ captures every speed-$q$ upward transit trajectory in $[-b,b]\times\R$.
\end{theorem}
\begin{proof}
Parameterize the right-hand downward leg $(a,-\tau)$ by centered time $\tau\in[-h,h]$. During $\tau\in[0,p]$, use Lemma~\ref{lem:attached-fan}'s fan $P_s(\tau)$ on $[a-1, a-\cos\theta(\tau)]$, extended leftward by the flat front $y = q\tau$ on $[-b, a-1]$ (inside $[-b,b]$ since $a-1\ge -b$ by~\eqref{eq:obs-4}) and rightward by the chord from $X_{\mathrm{fan}}(\tau)$ to $W_-(\tau):=(b, -\tau-G)$ (inside $\overline{B}_1(a,-\tau)$ by convexity since $\delta^2+G^2\le 1$). During $\tau\in[-p,0]$, use the vertically reflected incoming fan at time $-\tau$, joined to the same wall point $W_-(\tau)$. The graph coincides with $y=0$ on $[-b,a-1]$ and the chord $(a-1,0)\to(b,-G)$ at $\tau=0$, giving a continuous eikonal barrier across $[-p,p]$.

Next, shift to departure time $t := \tau - h$, so the succeeding leftward diagonal leg runs over $u\in[0, L_{\mathrm{diag}}]$ as $g_{\mathrm{diag}}(u) = (a - c_{\mathrm{leg}} u,\, -h + d_{\mathrm{leg}} u)$, and handoff $\tau = p$ is $t = -v$. The single-leg spacetime capture tube $\mathcal{C}_{\mathrm{diag}} := \{(t,P) : \exists u\in[0,L_{\mathrm{diag}}],\, \|P - g_{\mathrm{diag}}(u)\|_2 + q|u-t| \le 1\}$ is convex in $(t,P)$ (by convexity of $\|\cdot\|_2$ and $|\cdot|$ under joint interpolation of $(t,P,u)$), and any intruder meeting $\mathcal{C}_{\mathrm{diag}}$ cannot avoid $\overline{B}_1(g_{\mathrm{diag}}(u))$ at time $u$ by $\|p(u)-p(t)\|_2 \le q|u-t|$. For any point $P = (a+x, -h+y)$ at time $t = -v$ ($v>0$), the witness $U := (-m x + \kappa y)/\alpha$ satisfies $U - (-v) \ge 0$ when $U\ge 0$, and direct orthogonal expansion along $(c_{\mathrm{leg}}, d_{\mathrm{leg}})$ gives $\|P - g_{\mathrm{diag}}(U)\|_2 = N/\alpha$ whenever $N := d_{\mathrm{leg}} x + c_{\mathrm{leg}} y \ge 0$, yielding the exact affine identity
\begin{equation}\label{eq:obs-witness}
    \|(a+x, -h+y) - g_{\mathrm{diag}}(U)\|_2 + q(U + v) = qv + \kappa x + m y.
\end{equation}
At handoff $\tau=p$ ($t=-v$), the fan $s\mapsto P_s(p) = (a+x_s, -h+y_s)$ runs from $(x_0,y_0)=(-1, h+qp)$ ($s=0$) to $(x_1,y_1)=(-C, h-p+S)$ ($s=p$). Differentiating along $s\in[0,p]$ using $\partial_s P_s(p) = (p-s)u'(s)(1, -u/w)$ with $0 \le u/w \le \kappa/m$ and $d_{\mathrm{leg}} m - c_{\mathrm{leg}}\kappa = q > 0$ gives:
\begin{align*}
    \frac{d}{ds}(qv + \kappa x_s + m y_s) &= (p-s)u'\Bigl(\kappa - m\frac{u}{w}\Bigr) \ge 0, \qquad
    \frac{dN_s}{ds} = (p-s)u'\Bigl(d_{\mathrm{leg}} - c_{\mathrm{leg}}\frac{u}{w}\Bigr) \ge 0, \\
    \frac{dU_s}{ds} &= -\frac{(p-s)u'}{\alpha}\Bigl(m + \kappa\frac{u}{w}\Bigr) \le 0.
\end{align*}
Hence~\eqref{eq:obs-6} implies $N_s \ge N_0 \ge 0$ and $0 \le U_1 \le U_s \le U_0 \le L_{\mathrm{diag}}$ for all $s\in[0,p]$, and~\eqref{eq:obs-3} bounds the maximum cost at $s=p$ by $qv - \kappa C + m(h-p+S) = h(m+q) - K \le 1$. Thus the entire handoff fan lies inside $\mathcal{C}_{\mathrm{diag}}$. Likewise, for $t = -w \in [-v,0]$ ($w\in[0,v]$), the wall point $W_-(h-w) = (b, -h+w-G)$ satisfies $\delta^2 + (w-G)^2 \le 1 - 2Gw + w^2 \le (1-qw)^2$ and $1-qw\ge 0$ by~\eqref{eq:obs-5} (since $(1-qw)^2 - (1-2Gw+w^2) = w[2(G-q) - \alpha^2 w] \ge w[2(G-q)-\alpha^2 v] \ge 0$), so $(-w, W_-(h-w))\in\mathcal{C}_{\mathrm{diag}}$ via witness $u=0$. By convexity of the $t=-v$ slice of $\mathcal{C}_{\mathrm{diag}}$, replacing the fan-and-chord at $t=-v$ by the straight chord from $(a-1, qp)$ to $W_-(p)$ sweeps a region wholly inside $\mathcal{C}_{\mathrm{diag}}$.

Finally, interpolate in spacetime from $\bigl(-v, (a-1, qp)\bigr)\in\mathcal{C}_{\mathrm{diag}}$ to $\bigl(L_{\mathrm{diag}}, (-a, h+G)\bigr)\in\mathcal{C}_{\mathrm{diag}}$ (the latter has $\|(-a,h+G)-g_{\mathrm{diag}}(L_{\mathrm{diag}})\|_2 = G \le 1$ at $u=L_{\mathrm{diag}}$). Because $q(h+L_{\mathrm{diag}}) = qh + 2q\ell = h + G$, the interpolated point $U(t)$ lies on the rising flat level $F_+(t) = q(h+t)$ throughout $t\in[-v, L_{\mathrm{diag}}]$, with $x_U(0) \le a$. The spacetime reflection $L(t):=-U(L_{\mathrm{diag}}-t)$ lies in the same diagonal tube on the lower level $F_-(t)=-h-G+qt$. For $0\le t\le L_{\mathrm{diag}}$, the affine horizontal coordinates obey $x_U(t)-x_L(t)=x_U(0)-a\le0$. Use $F_+$ on $[-b,x_U(t)]$, the captured chord $U(t)L(t)$, and $F_-$ on $[x_L(t),b]$.

During $-v\le t\le0$, use $F_+$ to the left of $U(t)$ and join $U(t)$ to $W_-(h+t)$ in the same convex tube. At departure the swept triangle has vertices $U(0)$, $L(0)=(a,-h-G)$ and $(b,-h-G)$, all in the tube. After arrival use the spacetime reflection of this construction, joining the left wall at offset $+G$ to $L(t)$; the arrival triangle is the reflection of the departure triangle. At the late handoff $\tau=-p$ on the next downward leg, reverse the fan-to-chord replacement in the tube. Switch the wall offset from $+G$ to $-G$ across the triangle formed by the disk contact point and the two wall points at offsets $\pm G$. This triangle lies in the current closed disk. The next incoming fan now has the required wall endpoint. Repeating with left and right interchanged gives a bounded periodic separator throughout the corridor.
\end{proof}

\subsection{The critical speed and exact attainment on \texorpdfstring{$[q_a, 1)$}{[q\_a, 1)}}

Specialize~\eqref{eq:diffraction-defs}--\eqref{eq:R-round} to the wall-optimal bowtie~\eqref{eq:bowtie-params}, where $c_{\mathrm{leg}} = \frac{\sqrt{c}}{2q}$ and $d_{\mathrm{leg}} = \frac{1-q}{2q}$, and write $\mathcal{R}_{\mathrm{round}}(q) := \mathcal{R}_{\mathrm{round}}(a(q),h(q);q)$. Observe the exact factorization
\[
    4q^2\bigl(\alpha^2 d_{\mathrm{leg}}^2 - q^2 c_{\mathrm{leg}}^2\bigr) = (1-q^2)(1-q)^2 - q^2(3q^2+2q-1) = (1+q)(1-2q)(1-q+2q^2).
\]
Because $1+q > 0$ and $1-q+2q^2 = 2(q-\frac{1}{4})^2 + \frac{7}{8} > 0$, we have $\kappa(q) > 0 \iff q \in (1/3, 1/2)$, with $\kappa(1/2) = 0$ (above $q=1/2$, $d_{\mathrm{leg}} < q$ so the obstacle-fan regime does not arise and Theorem~\ref{thm:two-level} already applies via $[59/125,1)\supset[1/2,1)$).

\begin{lemma}[Critical speed $q_a$ and wavefront certification on $\lbrack q_a, 59/125\rbrack$]\label{lem:qa-root}
On $(1/3, 1/2]$, where $\kappa(q)\ge 0$, the function $\mathcal{R}_{\mathrm{round}}(q)$ is continuous, satisfies $\lim_{q\downarrow 1/3}\mathcal{R}_{\mathrm{round}}(q) = +\infty$ and $\mathcal{R}_{\mathrm{round}}(1/2) = 2/\sqrt{7} < 1$, and crosses $1$ at a unique root $q_a$ enclosed in
\begin{equation}\label{eq:qa-enclosure}
    0.457892036651 < q_a < 0.457892036652, \qquad F(q_a) \in [3.36785112514,\, 3.36785112516]\,r.
\end{equation}
Moreover, for every $q \in [q_a, 59/125] = [q_a, 0.472]$, the wall-optimal bowtie~\eqref{eq:bowtie-params} satisfies all hypotheses~\eqref{eq:obs-3}--\eqref{eq:obs-6} of Theorem~\ref{thm:obstacle-criterion}.
\end{lemma}
\begin{proof}
Continuity on $(1/3,1/2]$ is immediate from $\sin\theta_*(q) \le \frac{q\kappa}{1+qm} < q$, and as $q\downarrow 1/3$, $c\downarrow 0$ sends $h(q)\to+\infty$ while $p(q)$ stays bounded, giving $\mathcal{R}_{\mathrm{round}}(q)\to+\infty$; at $q=1/2$, $\kappa=0\implies \theta_*=0, p=0$, $c(1/2)=3/4$, $h(1/2) = \frac{4}{3\sqrt{7}}$, and $m(1/2)+1/2 = \frac{3}{2}$, giving $\mathcal{R}_{\mathrm{round}}(1/2) = \frac{2}{\sqrt{7}} < 1$. Rigorous degree-6 rational Taylor-polynomial enclosures of $\arctan$ and $(1-t)^{-1/2}$ (machine-checked over $\R$ in \texttt{AttainmentRoot.lean} and verified independently in \path{anc/verify_certificates.py}) prove $\mathcal{R}_{\mathrm{round}}(0.457892036651) > 1 > \mathcal{R}_{\mathrm{round}}(0.457892036652)$. Furthermore, outward-rounded rational interval arithmetic in \path{anc/verify_certificates.py} certifies $1 + K - h(m+q) < 0$ on $(1/3, 0.456]$, $\frac{d}{dq}[1+K - h(m+q)] > 0$ on $[0.456, 0.460]$, and $1 + K - h(m+q) > 0$ on $[0.460, 1/2]$ (global uniqueness of $q_a$ on $(1/3, 1/2]$), together with strict auxiliary inequalities~\eqref{eq:obs-4}--\eqref{eq:obs-6} across $[q_a,59/125]$. Inequality~\eqref{eq:obs-3} is an equality at $q_a$ and holds strictly to its right.
\end{proof}

\begin{theorem}[Exact global optimum for $q_a \le q < 1$]\label{thm:qa-exact}
For every speed ratio $q \in [q_a, 1)$, the wall-optimal four-leg bowtie~\eqref{eq:bowtie-params} captures every speed-$q$ upward transit trajectory across $[-b(q),b(q)]\times\R$. Combining this with Theorem~\ref{thm:wall-upper}, the exact minimax corridor width over all bounded-excursion patrols is
\begin{equation}\label{eq:thm1-formula}
    \Wopt(q) = F(q) = 2r\sqrt{1 + \frac{1}{3q^2+2q-1}}, \qquad \forall\, q \in [q_a, 1).
\end{equation}
\end{theorem}
\begin{proof}
Theorem~\ref{thm:wall-upper} gives $\Wopt(q)\le F(q)$ on $(1/3,1)$. On $[q_a, 59/125]$, Lemma~\ref{lem:qa-root} and Theorem~\ref{thm:obstacle-criterion} prove that~\eqref{eq:bowtie-params} guards $[-b(q),b(q)]\times\R$. On $[59/125, 1)$, $s(q)=(1+q)h(q)\le 1$ and $\delta^2+R(q)^2=1$ hold by the computation following~\eqref{eq:bowtie-params}, so Theorem~\ref{thm:two-level} proves that the same bowtie guards $[-b(q),b(q)]\times\R$. Since $[q_a, 59/125]\cup[59/125, 1) = [q_a, 1)$,~\eqref{eq:thm1-formula} holds throughout $[q_a, 1)$.
\end{proof}

\subsection{Exact widths at fixed excursion}

A vertical displacement weight refines the midline bound and identifies the optimal patrol for a range of height caps. Use units $r=V=1$ throughout the proof and the parametrization below.

\begin{theorem}[Excursion-constrained width]\label{thm:fixed-height}
Fix $0<q<1$. For every $0\le u\le\min(q,1-q)$ put $k=\sqrt{4q^2-u^2}$. Every patrol that guards half-width $b$ with $|Y|\le H$ satisfies
\begin{equation}\label{eq:height-dual}
 kb\le\sqrt{1+k^2}+(1-q-u)H.
\end{equation}
To parametrize the heights where equality is attained, take
\[
 0\le d\le\min\{1/2,(1-q)/(2q)\},\qquad c=\sqrt{1-d^2},\qquad E=2q-(1-q)d>0 .
\]
Set $u=2qd$, $k=2qc$, $G=(1+k^2)^{-1/2}$, $\ell=G/E$, $H=d\ell$ and $a=c\ell$. Whenever $(1+q)H\le1$,
\[
 \Wopt(q,H)=2(a+kG).
\]
For $1/2\le q<1$ this curve gives the optimum up to the height of the wall-optimal bowtie~\eqref{eq:bowtie-params}. Above that cap the optimum is $F(q)$. For each $0<q<1$ the parametrization gives an interval of exact widths starting at $H=0$.
\end{theorem}
\begin{proof}
For $b\le1$ the bound follows from $k<\sqrt{1+k^2}$. For $b>1$ consider a complete positive excursion of the horizontal coordinate $X$ from zero and back, of duration $D$, that detects a runner on the right wall. Choose the wall contacts giving the smallest and largest detected phases. In chronological order write their times as $t_1,t_2$, their guard coordinates as $(X_i,Y_i)$ and their chord radii as $R_i=\sqrt{1-(b-X_i)^2}$. Write $\sigma=+1$ if the minimum phase is first and $-1$ otherwise. The span of the detected phases is
\[
 I=R_1+R_2+\sigma\bigl(q(t_2-t_1)-(Y_2-Y_1)\bigr).
\]
Since $k^2+u^2=4q^2$, the speed bound gives $2q-uY'\ge k|X'|$. Its integral before $t_1$ and after $t_2$ is at least $k(X_1+X_2)$. On the middle interval, of duration $t$ and vertical displacement $z$, the signed phase difference exceeds the integral of $2q-uY'$ by
\[
 -(2-\sigma)qt+(u-\sigma)z .
\]
Here $t\ge|z|$ and $|z|\le2H$. For $\sigma=+1$ this expression is at most $(1-q-u)|z|$. For $\sigma=-1$ it is at most $(1+u-3q)|z|$, which is bounded by $(1-q-u)|z|$ since $u\le q$. Finally $R_i+k(b-X_i)\le\sqrt{1+k^2}$. Combining the three parts gives
\[
 I\le\int_0^D(2q-uY')\,dt+2\bigl(\sqrt{1+k^2}-kb+(1-q-u)H\bigr).
\]
Reflection gives the same estimate on a left excursion. On a long window $[0,L]$, both walls require a total phase span at least $2qL-O(1)$ as in the proof of Theorem~\ref{thm:midline}. Each of the at most two truncated excursions contributes at most $qt+2H+2$ over duration $t$, while its integral is at least $qt$ since $u\le q$. Unused time has nonnegative integrand. The full integral is $2qL-u(Y(L)-Y(0))=2qL+O(1)$. Detecting arbitrarily late phases on both walls requires infinitely many switches. Each complete active excursion takes at least $2(b-1)$ time, so their number on $[0,L]$ is finite and tends to infinity with $L$. The parenthesis in the last display must be nonnegative, giving~\eqref{eq:height-dual}.

For attainment, take the four-leg patrol with half-amplitude $a$, height $H$ and wall inset $\delta=kG$. Its wall chord has half-height $2q\ell-(1-q)H=G$, and $\delta^2+G^2=1$. Theorem~\ref{thm:two-level} gives capture under $(1+q)H\le1$. The identities $ka-(1-q-u)H=G$ and $G+k\delta=\sqrt{1+k^2}$ give equality in~\eqref{eq:height-dual}. The height $H(d)$ increases strictly while $E>0$, since $d/E$ is strictly increasing and $G$ is nondecreasing. For $q\ge1/2$ the upper endpoint $d=(1-q)/(2q)$ gives the parameters of~\eqref{eq:bowtie-params} and satisfies $(1+q)H\le 1$ by the computation following~\eqref{eq:bowtie-params}, since $[1/2,1)\subset[59/125,1)$. At all larger caps the same patrol attains $F(q)$.
\end{proof}

In physical units, expanding~\eqref{eq:height-dual} along the equality curve at $H=0$, using $u=4q^2\sqrt{1+4q^2}\,H/r+O_q((H/r)^2)$, gives
\[
 \Wopt(q,H)=W_{\mathrm{bf}}(q)+\frac{1-q}{q}H
     -2q\sqrt{1+4q^2}\,\frac{H^2}{r}+O_q(H^3/r^2).
\]
For $q=1/4$ and $H=4r/5$, the increasing height curve selects the unique $d\in(0,1/2)$ satisfying $(2-3d)^2(5-d^2)=100d^2$. The two-level condition holds with equality and $\Wopt(1/4,4r/5)/r\approx6.6545215541$.

\section{Below the diffraction threshold (\texorpdfstring{$0 < q < q_a$}{0 < q < q\_a})}
\label{sec:below-qa}

Below $q_a$, $\mathcal{R}_{\mathrm{round}}(q)>1$, so the attached diffraction fan can no longer be handed off to the succeeding diagonal leg at width $F(q)$. We first construct explicit escaping transits against the canonical four-leg bowtie across $(1/3,q_a)$, then promote that single-patrol failure to a period-uniform quadratic gap $F(q)-\Wopt(q)=\Theta\bigl(r(q_a-q)^2\bigr)$, establish the sharp leg-count hierarchy ($4,5,8,12$ legs), and finally obtain computer-certified midrange bounds alongside the $\Theta(\varepsilon^{-1/4})$ triangulation law at $q_a$.

\subsection{Escape from the wall-optimal bowtie below \texorpdfstring{$q_a$}{q\_a}}
\label{subsec:canonical-escape}
Use \(r=V=1\). For the canonical bowtie~\eqref{eq:bowtie-params}
write \(D=3q^2+2q-1\), \(c=a/\ell\), \(d=h/\ell\), and
\(G=2q\ell-(1-q)h\). Set time zero at the middle of the downward
return at \(x=a\). For \(t\ge0\) the remaining right return ends at
\(h\), the next diagonal ends at \(h+2\ell\), and the left return ends
at \(3h+2\ell\). The guard satisfies
\(g(-t)=(g_x(t),-g_y(t))\).

\begin{theorem}[Escape with positive clearance]\label{thm:canonical-escape}
For every \(1/3<q<q_a\) there are a speed-\(q\) upward transit in
\([-b,b]\times\mathbb R\) and a number \(\rho(q)>1\) such that
\(|p(t)-g(t)|\ge\rho(q)\) at every real \(t\). The path can be taken
vertical outside a finite time interval. Here \(g\) is the canonical
wall-optimal bowtie at this \(q\).

Moreover, for \(0.4578\le q\le q_a\) the rounded path in the proof is a
speed-\(q\) upward path in the same strip whose distance from \(g(t)\)
is at least \(\mathcal R_{\mathrm{round}}(q)\ge1\) at every time. It has
velocity \((0,q)\) for \(|t|\ge S\), where \(S<3h+2\ell\), and its height at
\(t=\pm S\) is \(\pm21/10\).
\end{theorem}
\begin{proof}
We cover $(1/3,q_a)$ with three explicit trajectory families: a straight
interior lane on $(1/3,0.4563]$, a V-shaped lateral-and-vertical path on
$[0.4563,0.4578]$, and a circular-arc-and-tangent diffraction path on
$[0.4578,q_a]$.

In the first range \(1/3<q\le0.4563\)
let \(s=h(1+q)\), \(k=(d-q)/c>0\),
\(n=\sqrt{1+k^2}\), and \(L=s/(n+k)\).
In the \((X,qt-Y)\) plane, the right return centered at phase zero
is the segment from \((a,-s)\) to \((a,s)\).
The adjacent diagonal lines have equations
\(z=\pm(s-k(a-x))\); their distance from \((a-L,0)\) equals
\((s-kL)/n=L\). The right return is distance \(L\), and the left
column is at least \(2a-L\) away. Consecutive copies of each
diagonal are translated vertically by \(P=4q(h+\ell)=2(s+G)\).
Their ordinates on \([-a,a]\) have magnitude at most \(s\)
before translation, since \(s-2ak=G>0\). Hence all other copies
are at least \(P-s>s\) away. It follows that the straight lane
\(p(t)=(a-L,qt)\) has clearance at least \(L\) whenever \(L\le a\).
To check \(1<L<a\) on this whole range, set \(w=\sqrt{q/(1+q)}\).
Substitution gives
\[
  \frac La=\frac{1-q}{1-2q+2qw}<1 .
\]
The denominator is positive and \(w>1/2\). The derivative of this
ratio has the sign of \(1-2w-2q(1-q)w'<0\). Since \(a\) decreases
and \(L(0.4563)>1\) by rational radical enclosure, the assertion
holds throughout.

For \(0.4563\le q\le0.4578\), fix \(e=39/500\) and
\((u,v)=q(2e,1-e^2)/(1+e^2)\), so \(u^2+v^2=q^2\). Define
\[
 \begin{gathered}
 w=1+v,\quad z=d-v,\quad j=c+u,\qquad
 M=\sqrt{u^2+w^2},\quad N=\sqrt{j^2+z^2},\\
 R=\frac{h((1+d)u+cw)}{N+zM/w},\qquad x_0=a-RM/w .
 \end{gathered}
\]
Use \(p(t)=(x_0+u\min\{|t|,S\},vt)\), where \(S=21/(10v)\).
On the right return the complete relative line is
\((x_0-a+ut,wt)\) for \(t\ge0\), at perpendicular distance
\((a-x_0)w/M=R\) from the origin. At departure its displacement is
\((x_0-a+uh,wh)\); the succeeding relative line has velocity
\((j,-z)\). Its signed normal distance is
\([h((1+d)u+cw)-(a-x_0)z]/N=R\).
The exact inequalities for the remaining pieces are
\[
 \begin{gathered}
 R>1+1/3000,\quad -b<x_0,\quad x_0+uS<b,\quad
 a+x_0>R,\\
 h+2\ell<S<3h+2\ell,\qquad 21/10-h>R .
 \end{gathered}
\]
These hold by rational interval substitution on 64 equal subintervals
of \([0.4563,0.4578]\). Until the tail starts, the left return is
horizontally at least \(a+x_0\) away. Thereafter the vertical
separation is at least \(21/10-h\). Reflection supplies negative times.

Finally take \(0.4578\le q\le q_a\). Put
\(\alpha=\sqrt{1-q^2}\), \(m=qd+\alpha c\), \(\kappa=\alpha d-qc\),
\((u,v)=q(\kappa,m)\), and
\[
 \theta_*=\arctan\frac{u}{1+v},\qquad
 J(\theta)=\int_0^\theta
     \frac{d\phi}{\cos\phi+\sqrt{q^2-\sin^2\phi}},\qquad
 K=\kappa\cos\theta_*-m\sin\theta_*+(m+q)J(\theta_*).
\]
Take \(R=h(m+q)/(1+K)=\mathcal R_{\mathrm{round}}(q)\) and put
\(\tau=RJ(\theta_*)\). By Lemma~\ref{lem:qa-root}, \(R>1\) for
\(q<q_a\) and \(R=1\) at \(q=q_a\).
For \(0\le t\le\tau\) let
\[
 t=RJ(\theta),\qquad p(t)=(a-R\cos\theta,R\sin\theta-t).
\]
This arc has constant distance \(R\) from the guard. Setting
\(w=R\dot\theta=\cos\theta+\sqrt{q^2-\sin^2\theta}\) gives
\(|p'|^2=w^2-2w\cos\theta+1=q^2\). At release the velocity is
\((u,v)\): the vector \((u,1+v)\) has direction
\((\sin\theta_*,\cos\theta_*)\) and length equal to the positive
branch \(w\). Continue on its absolute tangent
\[
 p(t)=(a-R\cos\theta_*+u(t-\tau),vt+B),\qquad
 B=R\sin\theta_*-(1+v)\tau .
\]
This is tangent to the return circle in relative coordinates.
On the next diagonal the normal \((\kappa,m)\) is unit and
\(-c\kappa+dm=q\); both guard and intruder have normal velocity \(q\).
Their constant signed separation is \(h(m+q)-RK=R\).
At \(S=(21/10-B)/v\) freeze the horizontal coordinate and continue
vertically at speed \(q\). The following checks hold on the closed interval
\([0.4578,0.4579]\), which contains the required range:
\[
 \begin{gathered}
  0<\tau<h,\quad a>R,\quad h+2\ell<S<3h+2\ell,\quad 21/10-h>R,\\
  -b<a-R,\qquad a-R\cos\theta_*+u(S-\tau)<b .
 \end{gathered}
\]
The horizontal coordinate increases up to the tail. On the left
return the horizontal separation before the tail is at least
\(2a-R>R\); after the tail the vertical separation suffices.
Reflect the whole path in time and height to complete the transit.
Every piece keeps distance at least \(R\) from the guard. This proves the
assertion on the rounded path, including the case \(q=q_a\).

For a reproducible enclosure of the rounded-arc time let
\(\sigma=\sin\theta_*\). Substitution \(\sin\phi=\sigma t\) gives
\[
 J(\theta_*)=\frac{\sigma}{1-q^2}
  \left(1-q\int_0^1
   \frac{\sqrt{1-(\sigma^2/q^2)t^2}}{\sqrt{1-\sigma^2 t^2}}\,dt\right).
\]
For both binomial factors the parameter \(x\) lies in \([0,1/10)\).
Truncate each after degree nine in \(t^2\), multiply, and integrate
the polynomial term by term. The absolute coefficient bounds by one
give each factor a uniform remainder at most \(x^{10}/(1-x)\).
Both factors and truncated factors have magnitude below two, so twice
the sum of the two remainders bounds the product error. All radicals
use outward integer-square-root enclosures. The read-only ancillary
checker \path{anc/escape/verify_escape.py} replays every range,
including the seven V-shaped-path inequalities and finite-piece rounded margins.
\end{proof}

For each fixed \(q\) the positive clearance gives a robust finite
crossing as well. Truncate the vertical tails at heights beyond
\(\pm(h+2)\). On this compact interval small uniform changes of the
guard preserve avoidance. If the wall positions move slightly,
projection of the intruder's horizontal coordinate onto the new
corridor preserves its speed bound and costs at most the wall displacement
in clearance. Safe tail extensions are available for uniformly bounded
guards when the truncation heights are outside their full sensing slab.

\subsection{Period-uniform gap and quadratic loss near \texorpdfstring{$q_a$}{q\_a}}
\label{sec:global-gap}

At the straight-wall ceiling the optimal bowtie admits an escaping transit for
every $1/3<q<q_a$ by Theorem~\ref{thm:canonical-escape}. To pass from this
single patrol to a statement about the supremum, we must also rule out
patrols with increasing periods or increasing excursion.

\begin{theorem}[Gap uniform over all periods]\label{thm:global-gap}
For every $1/3<q<q_a$ there is an $\epsilon(q)>0$ such that
\[
   \Wopt(q)\le F(q)-r\epsilon(q).
\]
Together with Theorem~\ref{thm:qa-exact}, this identifies $[q_a,1)$ as
the full set of speeds in $(1/3,1)$ where the straight-wall ceiling is attained
as a supremum.
\end{theorem}

\begin{theorem}[Quadratic loss at the threshold]\label{thm:global-quadratic}
There are constants $k_1,k_2,\eta>0$ such that
\[
 k_1r(q_a-q)^2\le F(q)-\Wopt(q)\le k_2r(q_a-q)^2
 \qquad(q_a-\eta<q<q_a).
\]
The upper inequality holds with $k_2=175$ throughout $[0.4578,q_a]$.
Proposition~\ref{prop:tilted} lowers this constant to $51.605$ with
an ordinary bowtie. Theorem~\ref{thm:bent-deficit} gives $48.7$ with a
twelve-leg patrol.
\end{theorem}

\subsubsection{Control of ordered wall contacts}

Use units $r=V=1$ and restrict $q$ to a compact interval
$J\subset(1/3,1/2)$. Write
\[
 D=3q^2+2q-1,\quad b_*=\sqrt{1+D^{-1}},\quad
 a_*=\frac{1}{\sqrt{D(1+D)}},\quad R_*=\frac{1}{\sqrt{1+D}},
\]
\[
 c=\frac{\sqrt D}{2q},\quad d=\frac{1-q}{2q},\quad
 \ell_*=a_*/c,\quad h_*=d\ell_*,\quad T_*=4(h_*+\ell_*).
\]
These are the parameters of~\eqref{eq:bowtie-params}, where the present $D$
is written $c$. In this section $c$ and $d$ are the direction cosines of a
transfer, so $c^2+d^2=1$.
In each half-cycle the bowtie descends at one pole for time $2h_*$ and
transfers to the other pole at velocity $(\pm c,d)$ for time $2\ell_*$.

Let $g=(X,Y)$ be a periodic guard detecting both straight-wall families,
with period $T$, half-width $b$, $P=qT$, $\alpha=\sqrt{1-q^2}$ and
$\Gamma=(\alpha X,z)$, where $z=qt-Y$. Theorem~\ref{thm:wall-upper}
gives $\Delta=b_*-b\ge0$. We work with $\Delta$ small enough that $b>1$.

\begin{lemma}[Ordered comparison at the original radius]\label{lem:gap-order}
There is a speed-bounded periodic comparison patrol with the same entire
lifted phase image as $g$, up to a phase translation. It has an even number
$N$ of alternating wall groups, with chronological markers
$t_i^-\le t_i^+\le t_{i+1}^-$. Write
\[
 p_i^\pm=z(t_i^\pm),\qquad x_i^\pm=|X(t_i^\pm)|,\qquad
 r_i^\pm=\sqrt{1-(b-x_i^\pm)^2}
\]
for the comparison patrol. Then $p_i^-\le p_i^+\le p_{i+2}^-$,
$p_{i+N}^\pm=p_i^\pm+P$, and the periodically extended intervals
$[p_i^--r_i^-,p_i^++r_i^+]$ cover each wall's phase line.
\end{lemma}
\begin{proof}
Fix the original period. Apply Lemma~\ref{lem:uncrossing} with excess
radius $e$ tending to zero. By its last assertion the embedded comparison
images are within Hausdorff distance $e$ of the original cylinder image, so
they converge to that entire image. Their
lengths are at most the original length, so Lemma~\ref{lem:phase-cylinder}
lifts each to a speed-bounded patrol of period $T$. Choose the outermost
contacts, groups and chronological markers as in Step~1 of the proof of
Theorem~\ref{thm:wall-upper}, with radius $1+e$ in place of
$r_\varepsilon$. That step provides three things. One period contains
finitely many alternating wall groups. Each group has a marker $t_i^-$
minimizing $p-R$ over the group and a marker $t_i^+\ge t_i^-$ maximizing
$p+R$, so that the intervals $[p_i^--r_i^-,p_i^++r_i^+]$ of the groups at
a wall cover its phase line. The selected contacts at each wall appear in
nondecreasing lifted phase order, which gives
$p_i^-\le p_i^+\le p_{i+2}^-$. Start the period at a left marker and subtract
the initial physical height. The resulting translations of phase modulo $P$
have a convergent subsequence, and the normalized patrols converge uniformly
along a further subsequence by equicontinuity. At excess radius $e$ the
opposite active zones are separated by $2(b-1-e)$, so the number of groups
is bounded for this fixed period. Extract a constant count and convergent
marker times. Closed interval coverage and the lifted order pass to the
limit. Its cylinder image agrees with the translated original image. Each
patrol advances by $P$ per period, so its entire lifted image consists of
all the translates of its image over one period. This yields equality of
the lifted images as well.
\end{proof}

On the comparison patrol put $u_i=t_i^+-t_i^-$, $v_i=t_{i+1}^--t_i^+$ and
\[
 z_i=Y(t_i^-)-Y(t_i^+),\quad w_i=Y(t_{i+1}^-)-Y(t_i^+),\quad
 l_i=x_i^++x_{i+1}^-,\quad A=\sum_i u_i,\quad Z=\sum_i z_i=\sum_i w_i .
\]
Use $j$ to index all $2N$ marker endpoints. Define the coverage excess $O$,
the marker defects $L_j$, and the two speed defects $f_i,e_i$ by
\[
 \begin{aligned}
 O&=qA+Z+\sum_jr_j-2qT, & L_j&=\sqrt D(x_j+\Delta)-r_j,\\
 f_i&=u_i-z_i, & e_i&=v_i-cl_i-dw_i .
 \end{aligned}
\]
Coverage gives $O\ge0$ and the speed bound gives $f_i,e_i\ge0$.
At $a_*$ the graph of $x\mapsto\sqrt{1-(b_*-x)^2}$ has tangent $\sqrt D x$.
Its second derivative is at most $-1$ on the interior of the active interval,
so $L_j\ge(x_j+\Delta-a_*)^2/2$, including the endpoints by continuity.
Using $\sum_i l_i=\sum_j x_j$, $2qc=\sqrt D$ and $q+2qd=1$ gives the exact identity
\begin{equation}\label{eq:gap-defects}
 O+\sum_j L_j+q\sum_i f_i+2q\sum_i e_i=2N\sqrt D\,\Delta .
\end{equation}

\begin{lemma}[Mean error on each pair of legs]\label{lem:gap-mean}
With
\[
 H_i=|x_i^--a_*|+|x_i^+-a_*|+
 |u_i-2h_*|+|v_i-2\ell_*|+\sqrt{u_if_i}+\sqrt{v_ie_i},
\]
one has $\sum_i H_i\le C_J N\sqrt\Delta$ for $0\le\Delta\le\Delta_J$. The
constants $C_J$ and $\Delta_J>0$ depend only on $J$, not on the patrol or
its period.
In particular $|T-NT_*/2|\le C_JN\sqrt\Delta$.
\end{lemma}
\begin{proof}
Constants in the proof are uniform over $J$. By~\eqref{eq:gap-defects} and
Cauchy--Schwarz, $\sum_j|x_j-a_*|\le CN\sqrt\Delta$.
Since $r_j=\sqrt D(x_j+\Delta)-L_j$ and $R_*=\sqrt D a_*$, the same bound
holds for $\sum_j|r_j-R_*|$. Let $\sigma_i$ be the sign of the horizontal
transfer displacement. The physical speed bound gives
\[
 \int_{\mathrm{transfer}\ i}|(\sigma_iX',Y')-(c,d)|^2\,dt\le2e_i,
 \qquad
 \int_{\mathrm{return}\ i}|(X',Y')-(0,-1)|^2\,dt\le2f_i .
\]
Also $v_i\le(2cb+e_i)/(1-d)$, so $\sum_i v_i\le CN$. Applying
Cauchy--Schwarz to the first integral and summing bounds
$\sum_i(|cv_i-l_i|+|dv_i-w_i|)$ by $CN\sqrt\Delta$. Thus for
$K_i=qv_i-w_i$, $K_*=2q\ell_*-2h_*$,
\[
 \sum_i(|v_i-2\ell_*|+|K_i-K_*|)\le CN\sqrt\Delta .
\]

Put $G_i=p_{i+2}^--p_i^+\ge0$. Replacing every endpoint radius by $R_*$
uncovers total measure at most $\sum_j|r_j-R_*|$ over the two phase
circles. After the replacement both endpoints are ordered on each wall,
so the uncovered measure equals $\sum_i(G_i-2R_*)_+$.
Telescoping on both walls gives
\[
 \sum_i(G_i-2R_*)=2P-qA-Z-2NR_*=\sum_j(r_j-R_*)-O .
\]
These two estimates imply $\sum_i|G_i-2R_*|\le CN\sqrt\Delta$.
Finally,
\[
 G_{i-1}=K_{i-1}+(1+q)u_i-f_i+K_i,\qquad R_*-K_*=(1+q)h_*,
\]
so $\sum_i|u_i-2h_*|\le CN\sqrt\Delta$. In particular $\sum_i u_i\le CN$.
Cauchy--Schwarz and~\eqref{eq:gap-defects} now bound the two square-root
sums in $\sum_i H_i$ by $CN\sqrt\Delta$.
\end{proof}

\subsubsection{Recovery of the original patrol in a phase interval}

The comparison may have placed points of the original image in a different
temporal order. We recover the original timing over a long interval and
control all visits from times outside it.

\begin{lemma}[A synchronized crossing]\label{lem:gap-crossing}
Fix a positive integer $m$ and a compact $J\subset(1/3,1/2)$. There are
constants $\Delta_0>0$ and $C$, depending only on $J$ and $m$, such that
the following holds for every $q\in J$ and every periodic guard $g$ as
above, of any period, that detects both straight-wall families at
half-width $b=b_*-\Delta$ with $0\le\Delta\le\Delta_0$. The original
patrol $g$ has a lifted phase interval $I=[L,U]$ and times $t_0<t_1$ with
$z(t_0)=L$, $z(t_1)=U$ and $L<z(t)<U$ for $t_0<t<t_1$. We call $[t_0,t_1]$
a crossing. It has the following properties.
The width $U-L$ differs from $mqT_*$ by at most $C\sqrt\Delta$.
The duration $t_1-t_0$ differs from $mT_*$ by at most $C\sqrt\Delta$.
After the changes of origin $t'=t-t_0$ and $Y'=Y-(qt_0-L)$ the original
patrol is within $C\sqrt\Delta$ in physical position of the ordinary bowtie
started at its left-return midpoint $( -a_*,0)$ throughout the crossing.
Its total time at phases in $(L,U)$ outside the crossing is at most
$C\sqrt\Delta$.
\end{lemma}

\begin{proof}
Number the left groups $2j$. Their time and phase cuts are
$s_j=(t_{2j}^-+t_{2j}^+)/2$ and $C_j=(p_{2j}^-+p_{2j}^+)/2$.
Both lists are nondecreasing. Repeat the original patrol and its ordered
comparison over $6m$ copies of the period and use the enlarged phase
circumference. Their lifted images still agree, all sums and counts scale
by the same factor, and the resulting $N/2$ left cycles divide exactly
into blocks of $m$ cycles, with at least six blocks. Use $T,N,P$ henceforth
for the enlarged quantities. For each block between successive chosen
left midpoints, assign the time interval between its $s$ cuts and the phase
slab between its $C$ cuts. Let $H[j]$ be the sum of $H_i$ over the groups
touching the $j$th block, including both partial endpoint groups. Then
$\sum_jH[j]\le CN\sqrt\Delta$.

Here are the relevant properties of the ideal $m$-cycle arc in phase space.
Set $S=2(1+q)h_*$ and $K=K_*<0$. From the start of a left return the
left and right column phases are $[0,S]$ and $[S+K,2S+K]$. At horizontal
fraction $x\in[0,1]$ measured from the left pole, the two transfer phases
are $S+Kx$ and $P_0-Kx$, where $P_0=qT_*=2(S+K)$. The difference and
its complement modulo $P_0$ lie between $S+2K=2R_*>0$ and $S<P_0$.
The infinite lifted ideal curve is therefore simple. The phase line at the
left midpoint meets it only there. The closest transfer is $S/2+K=R_*>0$
away in phase. Between two such midpoint cuts the intervening arc stays
inside the slab, with the return at fixed $X=-a_*$ as its only part within $R_*/2$ of either boundary.

\begin{figure}[tbp]
 \centering
 \includegraphics[width=0.86\linewidth]{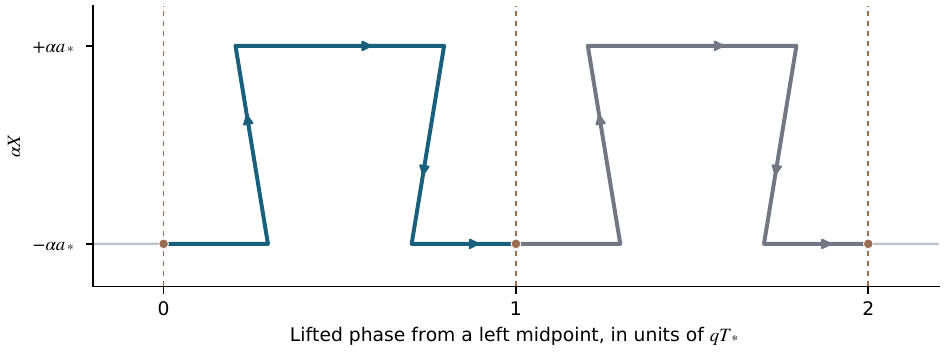}
 \caption{Two successive phase intervals of the ordinary bowtie at
 $q=0.45785$. The horizontal axis is lifted phase $z=qt-Y$, and the
 vertical axis is the scaled transverse position $\alpha X$. A transfer
 moves backward in phase while a downward return moves forward. Each cut
 through a left-return midpoint meets the curve once, and consecutive
 midpoint cuts differ by the period advance $qT_*$.}
 \label{fig:phase-interval}
\end{figure}

Choose a small constant $\kappa_0>0$. A block is regular when
$H[j]\le\kappa_0$. For each regular block, translate its initial $C$ cut
to zero and align the time on each of its legs linearly with the ideal leg.
The comparison's phase path is within $CH[j]$ of the corresponding ideal
arc in this parametrization. Indeed the integrated direction errors bound
every partial displacement error by $\sqrt{2u_i f_i}$ or $\sqrt{2v_i e_i}$,
and aligning a leg changes the ideal displacement by at most its duration
error. The initial horizontal error comes from the marker errors and the
partial return displacement. Linear interpolation of $z$ on a return also
shows $|z(s_j)-C_j|\le C H_{2j}$. Summing over the at most $4m+1$
legs of a block proves the assertion, including a slab-width error at most $CH[j]$.

All exceptional blocks together number at most $CN\sqrt\Delta$. They
occupy at most $CN\sqrt\Delta$ comparison time by Lemma~\ref{lem:gap-mean},
including exceptionally long legs. Project the images of their time
intervals onto the phase circle and enlarge the projected intervals by
one unit at each end. Their total length is at most $CN\sqrt\Delta$ since
$|z'|\le1+q$. Discard regular slabs meeting these enlarged intervals and
regular slabs adjacent in the block list to an exceptional block. For small
$\kappa_0$, regular slab widths have a uniform positive lower bound
$l=\tfrac12 m\min_{q\in J}qT_*$. Their interiors are disjoint by the
ordered $C$ cuts. A projected interval of length $L'$ can therefore meet
at most $L'/l+4$ regular slabs, splitting at the seam if needed. Altogether
$O(N\sqrt\Delta)$ blocks are discarded. For small $\Delta$ the retained
count is at least a positive constant times $N$.

For a retained slab $j$ let $E_j=H[j-1]+H[j]+H[j+1]$. Its entire original
phase image in the closed slab lies in a $CE_j$ neighborhood of its
translated ideal arc. To see this, use equality of the full lifted images
from Lemma~\ref{lem:gap-order}, absorbing the constant translation in the
choice of cuts. The discarded exceptional time blocks cannot visit the slab.
Any regular time block stays within $CH[\cdot]$ of its own slab. Choose
$\kappa_0$ so small that this error cannot cross an intervening regular
neighbor and that $CE_j\le3C\kappa_0<R_*/2$ on every retained slab,
uniformly on $J$. Here $C$ is the constant of the regular-block estimate
above, which does not depend on $\kappa_0\le1$. Both immediate neighbors of a retained slab are regular. They
can enter it only within $CE_j$ of the shared boundary in phase, where
the ideal neighbor has only its endpoint return at fixed transverse position. On this leg the
distance to the endpoint equals the phase difference. Including the
cut-width error gives the claimed $CE_j$ confinement and places every
visit to a boundary within $CE_j$ of its specified ideal endpoint.
The argument applies in every lifted copy.

In one lift $I=[L,U]$ let $t_1$ be the original patrol's first hit of $U$
and $t_0$ the last hit of $L$ before $t_1$. These exist because $Y$ is
bounded and $z(t)\to\pm\infty$ at the corresponding time ends. The
crossing stays strictly inside $I$ on $(t_0,t_1)$. Write
\[
 D_I=t_1-t_0,\qquad
 M_I=\int_{\R}{\bf1}_{\{L<z(t)<U\}}\,dt.
\]
Translation of $I$ by $P$ translates all its visits by $T$, so $M_I$
equals the time in all the open copies over one period. Hence
$\sum_{\rm retained}M_I\le T$, even when the guard waits on a cut.

On any subinterval of an original crossing,
\begin{equation}\label{eq:gap-clock}
 (1-q^2)\,\Delta t\ge\operatorname{len}(\Gamma)-q\,\Delta z .
\end{equation}
Indeed $|\Gamma'|\le1-qY'=1-q^2+qz'$ by the physical speed bound. Equality
holds along the unit-speed ideal legs, so the ideal arc length is $L_*=mT_*$
and its forward arclength $s$ corresponds to physical time
$t_*(s)=[s-q(\gamma_z(s)-\gamma_z(0))]/(1-q^2)$, where $\gamma$ denotes
the ideal phase polygon parametrized by arclength. This time is strictly
increasing since $|d\gamma_z/ds|\le1$.

The inverse arclength coordinate of the ideal polygon is Lipschitz on its
image with a constant $K_0$ uniform over $J$. On adjacent edges this follows
from the positive angle between the two outgoing rays at their common
vertex. Nonadjacent edges remain a positive distance apart. There are
finitely many edges of positive lengths, with fixed incidences and continuous
vertices on $J$, and the simplicity margins above are uniform. Extend the
coordinate to all phase points by truncating
$\inf_{0\le v\le L_*}(v+K_0|x-\gamma(v)|)$ to $[0,L_*]$, and denote the
result by $s(x)$. It is continuous, agrees with arclength on the arc, and
$|x-\gamma(s(x))|\le(1+K_0)\delta$ whenever $x$ lies in its $\delta$
neighborhood.

Take $\delta=CE_j$ on the principal crossing. Its initial and final
coordinates are within $C\delta$ of $0,L_*$. For any crossing point
$x$, continuity of $s$ supplies visits in order to all intervening ideal
vertex coordinates between the start and $s(x)$, choosing successive first
hits. Each visit is within $C\delta$ of the specified vertex. Chord lengths
between consecutive hits, including the two endpoints, sum to at least
$s(x)-C'\delta$, because the ideal pieces are straight and their number
is bounded. Omitting up to $C\delta$ of initial arclength costs at most this
length. The same argument from $x$ to the end gives a suffix length at least
$L_*-s(x)-C'\delta$. Both estimates allow backtracking.
Equation~\eqref{eq:gap-clock} on these prefix and suffix paths now gives
\[
 D_I\ge mT_*-CE_j,\qquad
 |t-t_0-t_*(s(x))|\le (D_I-mT_*)_++CE_j.
\]
The ideal patrol is uniformly Lipschitz in physical time, also past the
ends of its ideal block by periodic extension. The formulas
$X=\Gamma_x/\alpha$ and $Y=qt-z$ convert these estimates into a physical
position error at most $C(E_j+(D_I-mT_*)_+)$ after the stated translations.

If $n$ is the retained count, then $|T-nmT_*|\le CN\sqrt\Delta$ by the
duration estimate and the discarded count. Sum the preceding crossing-time
bounds and use $\sum M_I\le T$, $M_I\ge D_I$, and
$\sum_{\rm retained}E_j\le CN\sqrt\Delta$. The result is
\[
 \sum_{\rm retained}\big(E_j+(D_I-mT_*)_++M_I-D_I\big)\le CN\sqrt\Delta .
\]
As $n$ is at least a positive multiple of $N$, one slab has all three
errors at most $C\sqrt\Delta$. This gives every conclusion of the lemma.
\end{proof}

\subsubsection{Escaping the original patrol}

\begin{proof}[Proof of Theorem~\ref{thm:global-gap}]
Fix $q\in(1/3,q_a)$, and take an escaping transit of the ordinary
bowtie from Theorem~\ref{thm:canonical-escape}. Cut it at two times at
which its height is below $-h_*-2$ and above $h_*+2$, respectively.
Extend the endpoints vertically at upward speed $q$ on the two time tails.
The tails have clearance greater than two against the ordinary bowtie.
Continuity and closed-disk capture give clearance at least $1+\beta$
on the remaining compact segment for some $\beta>0$, which we may take
at most one. The evader phase $\zeta(t)=qt-y(t)$ is nondecreasing and
constant on each tail.

The transit of Theorem~\ref{thm:canonical-escape} is centered at the
midpoint of a right return, and an ideal block starts at the midpoint of a
left return. Place the center of the transit at time $T_*/2+jT_*$ of the
block for an integer $j\ge0$. Choose $j$ and $m$ so that its central
segment lies inside an ideal $m$-cycle block and its
entire phase range is more than two units from the phase boundaries. This can
be done simultaneously. Each increase of $j$ adds $T_*$ to the time
locations and $qT_*$ to the phase locations, so first move both lower
ends beyond the initial cuts, then choose $m$ with room at both upper ends.

Suppose a periodic patrol fully guards half-width $b=b_*-\Delta$ with
$\Delta$ sufficiently small. Use the selected crossing from
Lemma~\ref{lem:gap-crossing}, translating the evader by the same time and
height origins. Clamp its horizontal coordinate to $[-b,b]$. This changes
its position by at most $\Delta$, preserves its speed bound, and preserves
its phase. On the crossing the original guard and the ordinary bowtie are
within $C\sqrt\Delta$ in physical position. Thus the evader escapes there
when $C\sqrt\Delta+\Delta<\beta$.

To exclude capture at all other times, a visit to
$[L+\rho,U-\rho]$ outside $[t_0,t_1]$ costs at least
$\rho/(1+q)$ additional time in the open slab. Indeed a visit before
$t_0$ entails a passage from the lower boundary to the core by
$z(-\infty)=-\infty$, and a visit after $t_1$ entails a passage from
the core to the upper boundary by $z(+\infty)=+\infty$. Take a fixed
$\rho>0$ small enough, for example $1/2$, that the one-unit enlargement
of the evader's entire phase range lies in this core for small $\Delta$.
Lemma~\ref{lem:gap-crossing} makes the extra occupation smaller than
$\rho/(1+q)$. At any other time the phase difference then exceeds one,
hence so does the physical vertical separation. The resulting full transit
escapes, a contradiction. This excludes a fixed interval of widths below
the ceiling, uniformly over all periods. Theorem~\ref{thm:periodic-approximation}
transfers the bound to all bounded-excursion patrols. Restoring units
multiplies the width loss by $r$.
\end{proof}

\begin{proof}[Proof of Theorem~\ref{thm:global-quadratic}]
Use $J=[0.4578,q_a]$ at $r=1$ and write
$\mathcal R(q)=\mathcal R_{\mathrm{round}}(a_*,h_*;q)$ as
in~\eqref{eq:R-round}. By the last assertion of
Theorem~\ref{thm:canonical-escape}, the rounded transit has distance at
least $\mathcal R(q)$ from the ordinary bowtie at every time. It moves
vertically at speed $q$ outside its two outer joins, where its absolute
height is $21/10$, so its phase is constant there. Its join times and phase range are uniformly bounded on
$J$, and $qT_*$ has a uniform positive lower bound. One choice of $m$
therefore accommodates the central segments and phase buffers throughout
$J$. Lemma~\ref{lem:gap-crossing} gives uniform constants with this $m$.
For any successful periodic half-width with small $\Delta$, the preceding
escape argument requires
\[
                  \mathcal R(q)-1\le C\sqrt\Delta+\Delta .
\]
The rational enclosure in \path{anc/gap/check_wall_gap_upper.py} gives
$-7.514<\mathcal R'(q)<-7.435$ on $[0.4578,0.4579]$.
For the derivative enclosure, a binomial tail of degree $n=9$ and parameter $0\le x\le x_0<1/10$ has derivative bounded by
\[
 |x'|\sum_{k>n} k x_0^{k-1}
 =|x'|\,x_0^n\frac{(n+1)-nx_0}{(1-x_0)^2}.
\]
This follows by differentiating the uniformly convergent series on the compact parameter interval. Combining it with the value remainders for the two factors and the product rule encloses the derivative of the integral in the escape proof.
Since $\mathcal R(q_a)=1$ by Lemma~\ref{lem:qa-root}, it follows that
$\mathcal R(q)-1\ge7.435(q_a-q)$ on $J$. This proves the quadratic lower
bound for small $\Delta$. Widths with larger $\Delta$ satisfy the same
bound upon decreasing $k_1$. Taking suprema and applying
Theorem~\ref{thm:periodic-approximation} proves the global lower bound.

For the upper bound set $\lambda=1/\mathcal R(q)$ and scale $a_*,h_*$
by $\lambda$ while holding the speed and leg direction fixed. The fan angle,
integral $p$, and angular coefficients in~\eqref{eq:R-round} depend only
on the leg direction and $q$. Therefore the handoff ratio of the scaled
bowtie equals one exactly. Its wall displacement $G$ is $\lambda R_*$. Choose
the inset $\sqrt{1-\lambda^2R_*^2}$ to satisfy the wall condition at
equality. The resulting half-width is
\[
 B_q(\lambda)=\lambda a_*+\sqrt{1-\lambda^2R_*^2}.
\]
The accompanying checker recomputes 19 margins on the full scaled
family over $[0.4578,0.4579]$. They comprise every inequality among the
hypotheses of Theorem~\ref{thm:obstacle-criterion} other than those that
hold by construction, namely $\delta\ge0$, the wall condition
$\delta^2+G^2\le1$ and the handoff inequality~\eqref{eq:obs-3}, the last
two with equality. Three further margins are $1-q$ and the quantities
$m$ and $\kappa$ of~\eqref{eq:diffraction-defs}.
The checker encloses every margin above $0.0332$
using rational intervals and the binomial integral with its remainder.
This certifies an all-transit guard at half-width $B_q(\lambda)$.

Since $B_q(1)=b_*$ and $B_q'(1)=a_*-R_*^2/\sqrt{1-R_*^2}=0$, Taylor's
theorem and $B_q''(\lambda)=-R_*^2(1-\lambda^2R_*^2)^{-3/2}$ give, for
$0.4578\le q\le q_a$ where $\lambda\le1$,
\[
 2b_*-2B_q(\lambda)\le
 \frac{\sqrt{1+D}}{D^{3/2}}(1-\lambda)^2
 \le \frac{\sqrt{1+D}}{D^{3/2}}(7.514)^2(q_a-q)^2 .
\]
Here $1-\lambda=(\mathcal R-1)/\mathcal R\le\mathcal R-1$. The same
checker bounds the displayed coefficient by $175$. Rescaling gives the
announced upper estimate. At $q=0.45785$ the checker gives the two outer
inequalities in the following bracket. Theorem~\ref{thm:global-gap} gives
the middle inequality.
\[
        3.3682502470r<\Wopt(0.45785)<F(0.45785)<3.3682505519r .
\]
\end{proof}

\subsubsection{A tilted bowtie}\label{subsec:tilted}

The scaled bowtie in the preceding proof keeps the direction of the legs.
Changing the direction as well lowers the constant by a factor of more than
three. Use units $r=V=1$ as before. In this subsection $t$ denotes the
speed deficit $q_a-q$, not a time.

Consider a four-leg bowtie with half-amplitude $a$ and leg direction
$d_{\mathrm{leg}}=h/\ell$, and write
$c_{\mathrm{leg}}=\sqrt{1-d_{\mathrm{leg}}^2}$. The quantities
$m,\kappa,\theta_*,p$ of~\eqref{eq:diffraction-defs}
and~\eqref{eq:theta-p} depend only on $d_{\mathrm{leg}}$ and $q$. Put
\[
 A(d_{\mathrm{leg}},q)=
 \frac{c_{\mathrm{leg}}\bigl(1+\kappa\cos\theta_*-m\sin\theta_*+(m+q)p\bigr)}
      {d_{\mathrm{leg}}(m+q)},\qquad
 k(d_{\mathrm{leg}},q)=\frac{2q-(1-q)d_{\mathrm{leg}}}{c_{\mathrm{leg}}}.
\]
Since $h=ad_{\mathrm{leg}}/c_{\mathrm{leg}}$, the rounded-escape
ratio~\eqref{eq:R-round} of the bowtie equals $a/A$ and its wall chord is
$G=ak$. With the inset $\delta=\sqrt{1-a^2k^2}$ its half-width is
\[
 B(a,d_{\mathrm{leg}};q)=a+\sqrt{1-a^2k^2}.
\]

\begin{proposition}[Tilted bowtie]\label{prop:tilted}
Let $0.4578\le q\le q_a$ and $t=q_a-q$. The four-leg bowtie with
\[
 d_{\mathrm{leg}}=\frac{1-q}{2q}-3.735\,t,\qquad
 a=A(d_{\mathrm{leg}},q)
\]
and the inset above captures every upward transit in half-width
$B(a,d_{\mathrm{leg}};q)$. Consequently, in physical units,
\[
 0\le F(q)-\Wopt(q)\le 51.605\,r(q_a-q)^2 .
\]
Moreover
\[
 \limsup_{q\uparrow q_a}\frac{F(q)-\Wopt(q)}{r(q_a-q)^2}\le C_{\mathrm{tilt}},
 \qquad 51.28352020<C_{\mathrm{tilt}}<51.28352025 .
\]
\end{proposition}

\begin{proof}
For $a=A$ the handoff inequality~\eqref{eq:obs-3} is an equality, and the
inset makes the wall condition an equality. The checker
\path{anc/gap/check_tilted_gap.py} covers $[0.4578,q_a]$ by eight boxes,
with $q_a$ ranging over its enclosure~\eqref{eq:qa-enclosure}. On each box
it encloses the 19 margins described in the proof of
Theorem~\ref{thm:global-quadratic} above $0.0332$.
Theorem~\ref{thm:obstacle-criterion} therefore gives the guard.

Fix $q$ and write $x=(a,d_{\mathrm{leg}})$ and $x_*=(a_*,d)$ for the
wall-optimal bowtie. Let primes denote derivatives in $d_{\mathrm{leg}}$.
Then $k'=(2qd_{\mathrm{leg}}-(1-q))/c_{\mathrm{leg}}^3$ vanishes at $x_*$,
where also $k=\sqrt D$ and $a_*k^2=\delta$. Hence both partial derivatives
$\partial_aB=1-ak^2/\delta$ and
$\partial_{d_{\mathrm{leg}}}B=-a^2kk'/\delta$ vanish at $x_*$. The entries
of $-\nabla^2B$ are
\[
 \frac{k^2}{\delta^3},\qquad
 \frac{akk'(1+\delta^2)}{\delta^3},\qquad
 a^2\Bigl(\frac{(k')^2}{\delta^3}+\frac{kk''}{\delta}\Bigr).
\]
At the root, $A(d,q_a)=a_*$ because
$\mathcal R_{\mathrm{round}}(q_a)=1$. The checker encloses the derivative
in $q$ of $a-a_*$ along the family. Its absolute value is at most some
$\sigma<2.431$ on all boxes, so $|a-a_*|\le\sigma t$ by the mean value
theorem, and $|d_{\mathrm{leg}}-d|=3.735\,t$. In each box the checker
bounds the three entries above on the coordinatewise hull of the points
$x$ and $x_*$, which contains the segment joining them. With upper bounds
$M_{11},|M_{12}|,M_{22}$ of these entries, Taylor's theorem at $x_*$ gives
\[
 2b_*-2B(x;q)\le
 \bigl(M_{11}\sigma^2+2|M_{12}|\,\sigma\cdot3.735+M_{22}\,3.735^2\bigr)t^2
 <51.605\,t^2 .
\]

For the limit let $M=-\nabla^2B(x_*;q_a)$, which is diagonal with positive
entries. Write $\mathcal R(x;q)=a/A(d_{\mathrm{leg}},q)$ for the
rounded-escape ratio of the bowtie $x$, so that
$\mathcal R(x_*(q);q)=\mathcal R_{\mathrm{round}}(q)$. Let
$\nabla\mathcal R$ denote its gradient in $x$ at $(x_*;q_a)$ and
$E_0=-\mathcal R_{\mathrm{round}}'(q_a)>0$.
For a real number $v_2$ take $d_{\mathrm{leg}}=d+v_2t$ and
$a=A(d_{\mathrm{leg}},q)$. This family is smooth in $t$, has
$\mathcal R=1$, and equals $x_*(q_a)$ at $t=0$. Differentiating
$\mathcal R=1$ at $t=0$ shows $x-x_*(q)=vt+O(t^2)$ with
$\nabla\mathcal R\cdot v=-E_0$. The 19 margins are positive at
$t=0$, so the bowtie is a guard for small $t$, and Taylor's theorem gives
$(2b_*-2B)/t^2\to v^{\mathsf T}Mv$. On the line
$\nabla\mathcal R\cdot v=-E_0$ the least value of $v^{\mathsf T}Mv$ is
\[
 C_{\mathrm{tilt}}
 =\frac{E_0^2}{(\nabla\mathcal R)^{\mathsf T}M^{-1}\nabla\mathcal R},
 \qquad\text{attained at}\qquad
 v=-\frac{E_0\,M^{-1}\nabla\mathcal R}
         {(\nabla\mathcal R)^{\mathsf T}M^{-1}\nabla\mathcal R},
\]
whose second coordinate is $v_2\approx-3.7349$. The checker encloses
$C_{\mathrm{tilt}}$ in the stated interval.
\end{proof}

At $q=0.45785$ a tilted bowtie with $d_{\mathrm{leg}}=0.591904$ and
$h=0.80063570625788\,r$ has a strict handoff, and the same checker gives
$\Wopt(0.45785)>3.3682504612\,r$. This raises the lower endpoint of the
bracket at the end of the preceding proof. Theorem~\ref{thm:bent-deficit}
identifies $C_{\mathrm{tilt}}$ as the exact limiting coefficient for the
ordinary class and gives a smaller upper bound for the unrestricted game.

\subsection{When ordinary bowties are optimal}\label{sec:ordinary-transition}

Let $W_4(q)$ be the supremal one-way width for an \emph{ordinary
symmetric bowtie}. Its two vertical strokes and two straight diagonals
join $(\pm a,\pm h)$ up to translation. We allow either orientation,
transverse offsets, and $h=0$ (a transverse sweep). Each leg is run at
speed $V$. A smaller common speed gives the same class supremum, since
an escape against the faster patrol can be slowed by the same factor.

\begin{theorem}[Optimality transition for ordinary bowties]\label{thm:ordinary-transition}
An eight-leg periodic patrol guards a width strictly greater than
$W_4(q)$ at each $0<q<q_a$. The one-way widths satisfy
\[
 \begin{aligned}
 W_4(q)&<\Wopt(q) &&(0<q<q_a),\\
 W_4(q)&=\Wopt(q)=F(q) &&(q_a\le q<1),\\
 W_4(q)&=\Wopt(q)=r\sqrt{4+q^{-2}} &&(q\ge1).
 \end{aligned}
\]
\end{theorem}

There are two constructions behind the strict inequality.
For $0<q\le2/5$ we keep the vertical strokes and use three straight
pieces on each crossing. The resulting patrol beats the upper bounds
from straight runners for both orientations of the ordinary bowtie.
For $2/5\le q<q_a$, a runner that rounds the end of a stroke puts a
necessary bound on the height of every competitive ordinary bowtie.
At a class optimum that bound is tight. Bending each stroke inward
allows a first-order increase in height and amplitude, while the extra
time spent on the stroke is quadratic in the bend. Appendix~\ref{app:ordinary}
proves both comparisons. The wall-optimal bowtie attains the supremum
on $[q_a,1)$, and the transverse sweep does so on $[1,\infty)$, by
Theorems~\ref{thm:qa-exact} and~\ref{thm:midline}. The theorem compares a
specified four-leg class with all bounded-excursion patrols. The
minimum number of legs among globally optimal patrols below $q_a$
remains unknown.

\begin{figure}[tb]
\centering
\includegraphics[width=\linewidth]{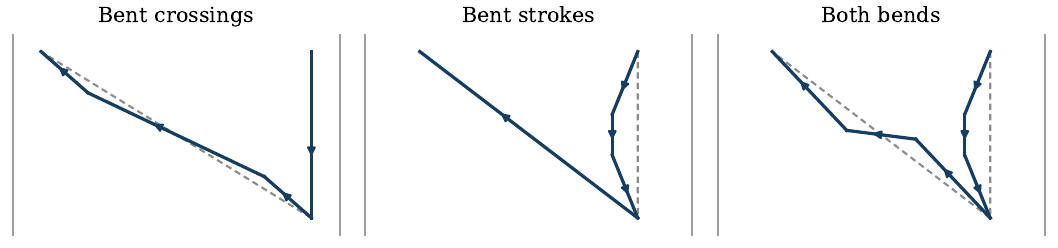}
\caption{One half-period of each construction, with arrows in
traversal order. Horizontal reflection gives the other half.
The dashed path is the ordinary stroke and crossing through the
same endpoints. The first panel shows the patrol of
Appendix~\ref{app:slow-comparison} at $q=1/4$. The other two panels
enlarge the bends for visibility.}
\label{fig:bent-patrols}
\end{figure}

\begin{theorem}[Size of the gain near the transition]\label{thm:bent-deficit}
Put $\varepsilon=q_a-q$. For $0.4578\le q\le q_a$,
\[
 W_4(q)\le F(q)-50r\varepsilon^2,\qquad
 \Wopt(q)\ge F(q)-48.7r\varepsilon^2 .
\]
Moreover, writing $C_4=C_{\mathrm{tilt}}$ from Proposition~\ref{prop:tilted},
\[
 \lim_{q\uparrow q_a}\frac{F(q)-W_4(q)}{r\varepsilon^2}=C_4,\qquad
 \limsup_{q\uparrow q_a}\frac{F(q)-\Wopt(q)}{r\varepsilon^2}
 \le C_{12},
\]
where $51.28352020<C_4<51.28352025$ and
$48.42508649<C_{12}<48.42508654$.
\end{theorem}

For this estimate we use three pieces on each stroke \emph{and} on
each crossing. Appendix~\ref{app:variation} computes their separate
contributions to the quadratic coefficient. At $q=0.45785$ the
twelve-leg patrol guards width greater than $3.368250466264r$,
while every ordinary bowtie has width below $3.368250463511r$.
A second twelve-leg patrol at the same speed has a rational
150-triangle certificate and guards width greater than
$3.3682504657890r$. The coefficient for the unrestricted optimum
is still unknown.

\subsection{Four arbitrary straight legs just below the critical speed}
Consider next all periodic polygonal patrols whose closed itinerary has
at most four straight legs. Each leg is traversed monotonically at speed
at most $V$, with speed changes and waits allowed and with vertices placed
anywhere in the corridor, repeating on a common period. Let $W_{\le 4}(q)$
denote the one-way width supremum over this general four-leg class, so
$W_4(q)\le W_{\le 4}(q)\le\Wopt(q)$. While Theorem~\ref{thm:ordinary-transition}
gives $W_4(q)=\Wopt(q)$ for $q\ge q_a$ and beats $W_4(q)$ with eight legs
for $0<q<q_a$, one may ask whether breaking symmetry or varying speed
allows four straight legs to improve on $W_4(q)$ just below $q_a$, and how
many legs are actually needed to surpass it.

\begin{theorem}[Arbitrary periodic four-leg patrols]\label{thm:four}
There exists $\varepsilon_0>0$ such that for $q_a-\varepsilon_0<q<q_a$,
\[
 W_{\le4}(q)=W_4(q)<\Wopt(q).
\]
Throughout this interval, every four-leg patrol attaining $W_4(q)$ is an
ordinary full-speed symmetric bowtie up to corridor symmetries. In particular,
any periodic polygonal patrol attaining $\Wopt(q)$ there requires at least
five straight legs, without assuming that an unrestricted optimizer exists.
\end{theorem}

Appendix~\ref{app:four-legs} proves Theorem~\ref{thm:four} by showing that,
near $q_a$, the wall-phase penalty for asymmetry or waiting dominates the
amplitude gained from relaxing the rounded-escape height constraint. A single
extra vertex on one wall stroke then suffices to beat the entire four-leg class.

\begin{theorem}[A fifth leg just below the transition]\label{thm:five}
For every $q<q_a$ sufficiently close to $q_a$, five is the least number of
straight legs among periodic patrols guarding a strictly wider corridor than
every periodic patrol with at most four straight legs: $W_{\le 4}(q)$ is
strictly exceeded by a full-speed five-leg patrol of width
\[
       F(q)-r C_5\varepsilon^2+O(r\varepsilon^3),
       \qquad \varepsilon=q_a-q,\qquad 50.84211365<C_5<50.84211368 .
\]
\end{theorem}

\begin{figure}[htbp]
\centering
\includegraphics[width=.96\linewidth]{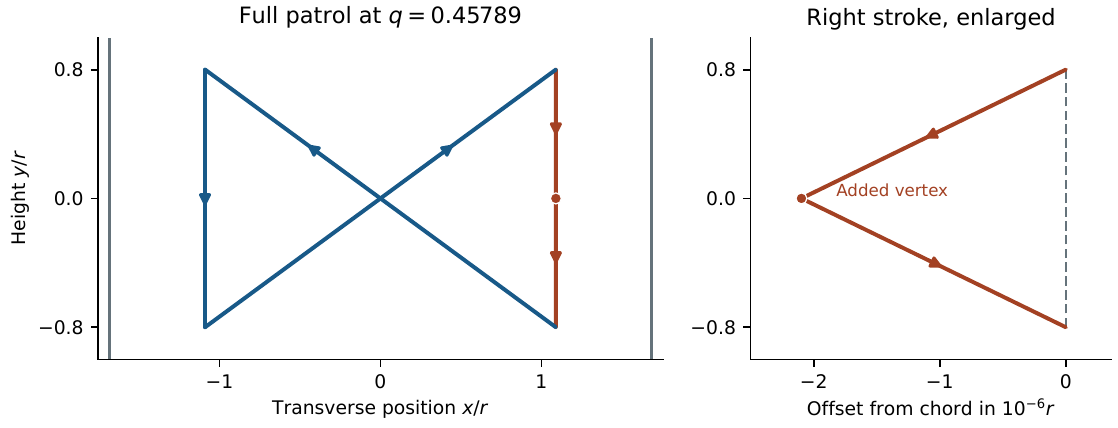}
\caption{A five-leg patrol at $q=0.45789$, with outer vertical lines marking
the corridor walls. Left: the full itinerary in true corridor proportions.
Right: the right wall stroke with only the horizontal coordinate magnified around
the stroke, resolving the inward midpoint displacement of about
$2.097\times 10^{-6}r$ from the dashed straight chord. At $q=0.45789$ the
interval check in Appendix~\ref{app:five-leg} certifies a strict gain over
every ordinary symmetric bowtie, while Theorem~\ref{thm:five} establishes strict
gain over all periodic four-leg patrols on a left neighborhood of $q_a$.}
\label{fig:fifth-leg}
\end{figure}

Figure~\ref{fig:fifth-leg} illustrates the construction: the right downward
stroke bends inward at its vertical midpoint while the left stroke and both
crossings remain single straight segments, with slightly unequal wall insets
matching the two stroke durations. Appendix~\ref{app:five-leg} certifies both
wall handoffs and computes $C_5$, which lies strictly below the shared four-leg
coefficient $C_4=51.2835202\ldots$ of Theorems~\ref{thm:bent-deficit}
and~\ref{thm:four}. Because Theorem~\ref{thm:four} relies on a compactness
localization, the neighborhood on which five legs beat all asymmetric four-leg
patrols is existential; at the explicit speed $q=0.45789$ shown in
Figure~\ref{fig:fifth-leg}, the interval check in Appendix~\ref{app:five-leg}
verifies that this five-leg patrol ($W>3.36787047176625\,r$) strictly exceeds
the ordinary-bowtie supremum ($W_4(0.45789)<3.36787047176461\,r$).

\subsection{Spacetime surface certificates and midrange width bounds}\label{sec:upper}
\subsubsection{Lower bounds from oriented spacetime surfaces}

We certify rigorous lower bounds $\Wopt(q) \ge W_{\mathrm{cert}}$ for multi-leg periodic patrols in exact rational arithmetic via triangulated spacetime surface cuts (translating $x\in[-b,b]$ by $+b$ to $[0,W]$ to match certificate coordinates):

\begin{proposition}[Triangulated spacetime surface-cut certificate]\label{prop:surface-cut}
Let $g:\R\to[0,W]\times\R$ be a $T$-periodic piecewise-affine patrol ($T>0, W>0$) with speed at most $1$ on every segment. Let $\Sigma$ be a finite oriented triangulated surface in $\R^3_{(t,x,y)}$ whose vertices satisfy $0\le t\le T$ and $0\le x\le W$ such that, when edge incidences are read modulo $T$ in the time coordinate:
\begin{enumerate}
    \item \emph{Closedness off the corridor walls:} Every directed edge of $\Sigma$ is paired with an oppositely oriented partner modulo $T$, except edges lying entirely in the left wall $\{x=0\}$ or right wall $\{x=W\}$.
    \item \emph{Unit negative projection degree:} Summing the signed areas of all triangles projected onto $(t,x)$ equals $-TW$ (equivalently, $\sum_{\triangle\subset\Sigma} n_y(\triangle) = -2TW$ where $\mathbf{n}(\triangle) = (p_1-p_0)\times(p_2-p_0)$).
    \item \emph{Facewise causal or single-segment capture certificate:} Every triangle $\triangle = \operatorname{conv}\{v_0,v_1,v_2\}\subset\Sigma$ ($v_i = (t_i, P_i)$ with $P_i\in[0,W]\times\R$) satisfies at least one of:
    \begin{itemize}
        \item \emph{(Type C --- Causal):} Its normal $\mathbf{n}=(n_t,n_x,n_y)$ satisfies $n_t \ge 0$ and $n_t^2 \ge q^2(n_x^2 + n_y^2)$; or
        \item \emph{(Type S --- Capture):} There exist witness times $w_0,w_1,w_2$ lying on \emph{one common affine segment} of $g$ such that $\|P_i - g(w_i)\|_2 + q|w_i - t_i| \le r$ for each $i\in\{0,1,2\}$.
    \end{itemize}
\end{enumerate}
Then $g$ captures every upward transit trajectory of speed at most $q$ in $[0,W]\times\R$, proving $\Wopt(q)\ge W$.
\end{proposition}
\begin{proof}
Extend $\Sigma$ $T$-periodically along $t\in\R$ to an infinite polyhedral $2$-chain $\widetilde{\Sigma}\subset\R\times[0,W]\times\R$. By Condition~1, $\partial\widetilde{\Sigma}$ is supported entirely in $\R\times\{0,W\}\times\R$, and by compactness of $\Sigma$, $|y|$ is bounded by some $Y_{\max}$ on $\widetilde{\Sigma}$. In the cylinder of $(t\bmod T,x)$ the projected chain has boundary supported on its two walls. Its multiplicity is therefore constant at generic interior points. Integrating that multiplicity and using Condition~2 gives $-1$. Thus in the lift the oriented intersection number of a vertical fiber from below to above the chain is $-1$. Every interior curve with endpoints on opposite sides of the height range has the same number: it is homotopic to such a fiber with its endpoints moving only in the two disjoint outer half-spaces.

To justify transverse intersections even for a wall-following or grazing intruder, suppose an undetected path exists and choose a finite crossing from below $-Y_{\max}-2$ to above $Y_{\max}+2$. Include a time buffer $r/q+1$ at both ends. On this compact enlarged interval the path has strictly positive excess distance over $r$. Contract its spatial coordinates slightly toward $(W/2,0)$ to gain strict speed and wall margins, then approximate by a piecewise linear time graph. Sufficiently small generic perturbations of the vertices preserve both margins, the excess distance and the endpoint height inequalities. Only finitely many faces meet the time window. The perturbed crossing avoids their edges and meets the faces transversely. A Type~S witness for any such intersection lies within the enlarged interval since $q|w-t|\le r$.

On any Type~C face, $\|p'(\tau)\|_2\le q$ and Cauchy--Schwarz give $\mathbf{n}\cdot c'(\tau) = n_t + n_x x'(\tau) + n_y y'(\tau) \ge n_t - q\sqrt{n_x^2+n_y^2} \ge 0$ a.e., so $c$ can never cross a Type~C face with negative orientation. Hence $c$ must meet some Type~S face $\triangle$ at some time $\tau_*$, say $(\tau_*, p(\tau_*)) = \sum_{i=0}^2 \lambda_i (t_i, P_i)$ ($\lambda_i\ge 0$, $\sum\lambda_i=1$). Because $g$ is affine on the common segment containing $w_0,w_1,w_2$, the map $(t,P,w)\mapsto \|P - g(w)\|_2 + q|w - t|$ is jointly convex, so $w_* := \sum_{i=0}^2 \lambda_i w_i$ satisfies $\|p(\tau_*) - g(w_*)\|_2 + q|w_* - \tau_*| \le r$. By $\|p(w_*) - p(\tau_*)\|_2 \le q|w_* - \tau_*|$ and the triangle inequality, $\|p(w_*) - g(w_*)\|_2 \le r$.
\end{proof}

A 64-leg periodic patrol at $q=1/4$ has an oriented certificate with 837 vertices and 1,592 faces. Its exact width exceeds the following decimal endpoint:
\begin{equation}\label{eq:q25-lower}
 \Wopt(1/4)>6.944691333740\,r .
\end{equation}
Six 18-leg patrol certificates give the following strict lower endpoints for $\Wopt(q)/r$:
\[
\begin{array}{c|rrrrrr}
q&1/4&7/25&3/10&1/3&7/20&2/5\\\hline
\text{lower}&6.9446905&6.100184&5.631526&4.97565174&4.694534&3.990757
\end{array}
\]
The rational patrol coordinates, triangle incidences and face witnesses accompany the checker in \texttt{anc/}.
\subsubsection{Flux across straight lanes and parity at the walls}

Normalize $r=1, V=1$, $b>1$, and let $J := [-b+1, b-1]$ be the interior gap separating the left wall footprint $[-b,-b+1]$ from the right wall footprint $[b-1,b]$.
\begin{enumerate}
    \item \textbf{Weighted Galilean lane flux:} For $a,c\ge 0$, define $\Phi_{a,c}(z) := \max\!\bigl(c,\, \frac{a|z|}{\sqrt{1-z^2}}\bigr)$ for $z\in(-1,1)$ and its continuous odd primitive $P_{a,c}:[-1,1]\to\R$ ($P_{a,c}(-z)=-P_{a,c}(z)$) by setting, for $z\in[0,1]$ and $D := \sqrt{a^2+c^2}$,
    \[
        P_{a,c}(z) := \begin{cases} c z, & z D \le c, \\ D - a\sqrt{1-z^2}, & z D > c. \end{cases}
    \]
    For any bounded nonnegative lane density $\lambda:[-b,b]\to[0,\infty)$ (extended by $0$ outside $[-b,b]$), when the guard is at $X\in[-b,b]$ with velocity $v=(v_x,v_y)$ ($v_x^2+v_y^2\le 1$), pairing the two boundary arcs of $\partial B_1(X,Y)$ at horizontal offset $z = x - X \in (-1,1)$ gives the instantaneous $\lambda$-weighted Galilean sweep rate
    \begin{equation}\label{eq:lane-flux}
        K_\lambda(X, v) := \int_{-1}^{1} \lambda(X+z)\,\Phi_{|v_x|, |q-v_y|}(z)\,dz = \int_{[-b,b]\cap[X-1,X+1]} \lambda(x)\,d_x P_{|v_x|,|q-v_y|}(x-X).
    \end{equation}
    \item \textbf{Two-wall parity credits and transverse weight $\beta(X)$:} For finitely many cutoffs $c_k\in(0,1]$ ($s_k := \sqrt{1-c_k^2}$) with weights $\eta_k\ge 0$ ($E := \sum_k \eta_k$), define the wall-footprint functions
    \[
        f_{c_k}(X) := \begin{cases} \dfrac{u}{2\sqrt{1-u^2}}, & b-1 < |X| \le b,\ u := b-|X| \in [0, s_k), \\ 0, & \text{otherwise}, \end{cases}
    \]
    choose any bounded measurable gap density $\rho:J\to[0,\infty)$ (extended by $0$ outside $J$) satisfying $\int_J \rho(x)\,dx \ge 4\sum_k c_k \eta_k$, and set $\beta(X) := 4\sum_k \eta_k f_{c_k}(X) + \rho(X)$.
\end{enumerate}

\begin{theorem}[Galilean flux and two-wall free-gap parity inequality]\label{thm:free-gap}
Let $b>1$, $q\in(0,1)$, $\lambda\ge 0$, $\eta_k\ge 0$, $\rho\ge 0$ with $\int_J \rho(x)\,dx \ge 4\sum_k c_k\eta_k$, and $\mu\in\R$. For every unit-speed patrol $g=(X,Y)$ on $[-b,b]\times\R$ with bounded longitudinal excursion $\sup_{t\in\R}|Y(t)|\le H < \infty$ that blocks all upward straight vertical runners $p_{x,\zeta}(t)=(x, qt-\zeta)$,
\begin{equation}\label{eq:unrestricted-cap}
    \sup_{X\in[-b,b],\, \|v\|_2\le 1} \Biggl\{ K_\lambda(X,v) + \Bigl(\sum_k \eta_k\Bigr)|q-v_y| + \beta(X)|v_x| - \mu v_y \Biggr\} \ge q\int_{-b}^b \lambda(x)\,dx + 3q\sum_k \eta_k.
\end{equation}
\end{theorem}
\begin{proof}
Fix a long window $[0,L]$ ($qL > 2H+4$). Every phase $\zeta\in[H+1, qL-H-1]$ on every lane $x\in[-b,b]$ satisfies $|qt-\zeta| > H+1 \ge |Y(t)|+1$ outside $[0,L]$, so $g|_{[0,L]}$ must sweep $\lambda$-weighted $(x,\zeta)$-measure at least $(qL - 2H - 2)\int_{-b}^b\lambda(x)\,dx$, giving
\[
 \int_0^L K_\lambda(X(t),g'(t))\,dt\ge(qL-2H-4)\int_{-b}^b\lambda(x)\,dx.
\]
Indeed, a point in the swept disk union that lies outside both endpoint disks encounters the moving boundary at least twice, apart from a set of area zero. The area formula applied to $(t,\theta)\mapsto (X(t),qt-Y(t))+(\cos\theta,\sin\theta)$ integrates the absolute normal velocity with multiplicity. Half of that integral is $\int_0^L K_\lambda\,dt$ by pairing the two arcs at each lane. The two endpoint disks counted with weight one half have weighted area at most $2\int\lambda$ in total (points in exactly one endpoint disk require one crossing). Also $|\int_0^L \mu Y'(t)\,dt| = |\mu(Y(L)-Y(0))| \le 2|\mu|H$.

It remains to prove the temporal credit inequality $E\int_0^L |q-Y'|\,dt + \int_0^L \beta(X)|X'|\,dt \ge 3qEL - E(6H+14)$ where $E = \sum_k \eta_k$. Partition $[0,L]$ into \emph{alternating wall-labelled blocks} $B$, where the active label ($\mathrm{L}$ or $\mathrm{R}$) switches only when $X(t)$ enters the opposite wall footprint $\{X < -b+1\}$ or $\{X > b-1\}$ (assigning any initial unlabelled prefix to the first visited wall). Discard the first and last blocks from the travel estimate; every remaining interior block $B$ enters its own footprint across $|X|=b-1$, reaches some maximum vertical chord radius $r_B := \max_{t\in B\cap\{|X|\ge b-1\}}\sqrt{1-(b-|X(t)|)^2} \in [0,1]$, exits $|X|=b-1$, and subsequently makes a complete crossing of $J=[-b+1,b-1]$ to the opposite footprint. On $B$'s own footprint, $b-|X|$ descends from $1$ to $\sqrt{1-r_B^2}$ and returns to $1$, so total variation gives
\[
    \int_{B\cap\{|X|>b-1\}} f_{c_k}(X)|X'|\,dt \ge 2\int_{\sqrt{1-r_B^2}}^{1} \frac{u\,\mathbf{1}_{\{\sqrt{1-u^2}>c_k\}}}{2\sqrt{1-u^2}}\,du = (r_B - c_k)_+,
\]
while the disjoint subsequent crossing of $J$ contributes $\int_{\text{gap}_B} \rho(X)|X'|\,dt \ge \int_J \rho(x)\,dx \ge 4\sum_k c_k\eta_k$. Summing over all interior blocks (and paying at most $8E$ for the two discarded endpoint radii $r_{B_{\mathrm{first}}}, r_{B_{\mathrm{last}}}\le 1$) gives
\[
    \int_0^L \beta(X)|X'|\,dt \ge 4\sum_{B\text{ int}}\sum_k \eta_k\bigl((r_B - c_k)_+ + c_k\bigr) \ge 4E\sum_{B\text{ all}} r_B - 8E.
\]
Independently, consider the phase map $\phi(t) := qt - Y(t)$ on $[0,L]$. Because $\phi(0)\le H$ and $\phi(L)\ge qL-H$, every regular value $\zeta \in [H+1, qL-H-1]$ has odd preimage count $n(\zeta) = \#\phi^{-1}(\zeta) \ge 1$. Let $\kappa(\zeta)\in\{0,1,2\}$ count how many of the two wall labels $\{\mathrm{L},\mathrm{R}\}$ have some block $B$ whose $\phi$-image $\phi(B)$ contains $\zeta$. Because $\mathrm{L}$-blocks and $\mathrm{R}$-blocks are disjoint in $[0,L]$, $\kappa(\zeta)=2$ forces at least two preimages of $\zeta$, hence $n(\zeta)\ge 3$ by odd parity; thus $n(\zeta) \ge 2\kappa(\zeta) - 1$ for a.e.\ $\zeta\in[H+1,qL-H-1]$. Meanwhile, each block $B$ can only detect wall phases within distance $r_B$ of its interval $\phi(B)$, and full coverage of $[H+1,qL-H-1]$ on both walls requires the two-wall uncovered measure $\int_{H+1}^{qL-H-1}(2-\kappa(\zeta))\,d\zeta$ to be covered by the $2r_B$ endpoint extensions across all blocks: $\int_{H+1}^{qL-H-1}(2-\kappa(\zeta))\,d\zeta \le 2\sum_B r_B$. By the Banach indicatrix theorem,
\[
    V_z := \int_0^L |q - Y'(t)|\,dt \ge \int_{H+1}^{qL-H-1} n(\zeta)\,d\zeta \ge \int_{H+1}^{qL-H-1} \bigl(3 - 2(2-\kappa(\zeta))\bigr)\,d\zeta \ge 3(qL - 2H - 2) - 4\sum_B r_B.
\]
Multiplying by $E$, adding the $\beta(X)$ bound ($-4E\sum_B r_B + 4E\sum_B r_B = 0$), combining with $K_\lambda$ and $-\mu Y'$, dividing by $L$, and letting $L\to\infty$ proves~\eqref{eq:unrestricted-cap}.
\end{proof}

\begin{proposition}[Vertical-drift floor of stateless flux and two-wall credits for $q\le1/3$]\label{prop:q13-floor}
At $q=1/3$, for every symmetric $\lambda\ge 0$, $\eta_k\ge 0$, $\rho\ge 0$, and $\mu\in\R$, if $W = 2b \le 6r$, then the left-hand supremum in~\eqref{eq:unrestricted-cap} is at least the right-hand demand. Conversely, at $W = 6.12r > 6r$, there exist rational multipliers for which the supremum drops strictly below demand, certifying $\Wopt(1/3) \le 6.12r < W_{\mathrm{area}}(1/3) = 2\sqrt{10}\,r \approx 6.3246r$.

More generally, let $0<q\le1/3$ and $n=\lfloor1/q\rfloor$. If $W=2b\le2nr$, then for every $\lambda\ge0$, $\eta_k\ge0$, $\rho\ge0$ and $\mu\in\R$ the left-hand supremum in~\eqref{eq:unrestricted-cap} is at least the right-hand demand. In particular, at $q=1/4$ no choice of multipliers in Theorem~\ref{thm:free-gap} excludes a width $W\le8r$.
\end{proposition}
\begin{proof}
Evaluate the bracket $\Psi(X,v)$ in~\eqref{eq:unrestricted-cap} at $v_\pm = (0,\pm 1)$. Since $v_x=0$, $P_{0,c}(z)=cz$, so $K_\lambda(X,v_\pm) = |q\mp 1| I(X)$ where $I(X):=\int_{[-b,b]\cap[X-1,X+1]}\lambda(x)\,dx$. Averaging $v_+$ and $v_-$ cancels $\mu v_y$: $\frac{1}{2}(\Psi(X,v_+)+\Psi(X,v_-)) = I(X) + \sum_k \eta_k$. When $2b\le 6$, $[-b,b]$ is covered by the three length-$2$ intervals centered at $-b+1, 0, b-1$, so $\max_X I(X) \ge \frac{1}{3}\int_{-b}^b\lambda$, matching $\frac{1}{3}\int_{-b}^b\lambda + \sum_k \eta_k$ at $q=1/3$. At $W = \frac{153}{25}\,r = 6.12r$, our exact rational certificate verifies $\sup\Psi < \text{Demand}$.

For the general statement let $2b\le2n$ and put $X_j=-b+(2j-1)b/n$ for $1\le j\le n$. These points lie in $[-b,b]$, and because $b/n\le1$ the $n$ intervals $[X_j-1,X_j+1]$ cover $[-b,b]$. Hence $\max_jI(X_j)\ge\frac1n\int_{-b}^b\lambda\ge q\int_{-b}^b\lambda$. Also $\sum_k\eta_k\ge3q\sum_k\eta_k$, since $q\le1/3$. For the maximizing $j$ the average $I(X_j)+\sum_k\eta_k$ of $\Psi(X_j,v_+)$ and $\Psi(X_j,v_-)$ is therefore at least the demand, and so is the larger of the two values. Symmetry of $\lambda$ is not used.
\end{proof}

Using exact rational Pythagorean cutoffs $(c_k,s_k) = \bigl(\frac{2mn}{m^2+n^2},\,\frac{m^2-n^2}{m^2+n^2}\bigr)$ and outward-rounded enclosures over $[-b,b]\times\{\|v\|_2\le 1\}$ (verified independently by \path{anc/verify_certificates.py} and replayed in Lean~4 in \texttt{ContinuumUpperCert.lean}), Theorem~\ref{thm:free-gap} establishes three unrestricted ceilings:
\begin{itemize}
    \item At $q = 37/100$: $\Wopt(37/100) \le \frac{107}{20}\,r = 5.35r < F(37/100) \approx 5.5266r$. By monotonicity of $\Wopt(q)$, $\Wopt(q) \le 5.35r < F(q)$ holds across the entire interval $[0.37,\, q_*)$ where $F(q_*)=5.35r$ ($q_* = (\sqrt{4+4800/9849}-1)/3\approx 0.37278$). This gives a numerical separation on that interval.
    \item At $q = 1/3$: $\Wopt(1/3) \le \frac{153}{25}\,r = 6.12r$ (finite as $F(q)\to\infty$ for $q\downarrow1/3$, and below $W_{\mathrm{area}}(1/3)\approx 6.3246r$).
    \item At $q = 1/4$: $\Wopt(1/4) \le \frac{81}{10}\,r = 8.10r$, strictly below $W_{\mathrm{area}}(1/4)=2\sqrt{17}\,r\approx 8.2462r$.
\end{itemize}
By Proposition~\ref{prop:q13-floor}, Theorem~\ref{thm:free-gap} cannot exclude any width up to $6r$ at $q=1/3$ or up to $8r$ at $q=1/4$. The last two ceilings are within $0.12r$ and $0.10r$ of these limits. The limits belong to this family of inequalities. They do not provide patrols of width $6r$ or $8r$.

\subsubsection{Upper bounds at bounded excursion}
\label{sec:bh2}

The height cap also bounds the vertical displacement available on each crossing. A supporting direction couples this displacement to transverse travel.

\begin{theorem}[Pythagorean inter-wall chord inequality]\label{thm:bh2}
Let $g=(X,Y):\R\to[-b,b]\times\R$ ($b>1$, $r=1$) be any unit-speed patrol with $\sup_{t\in\R}|Y(t)|\le H$ blocking straight wall runs at speed $q\in(0,1)$. Fix any unit pair $(\alpha_0,\beta_0)$ ($\alpha_0,\beta_0\ge 0$, $\alpha_0^2+\beta_0^2=1$) and any cutoff $(c_k,s_k)$ ($c_k>0, s_k\ge 0, c_k^2+s_k^2=1$) satisfying the side condition
\[
    s_k \ge q \alpha_0 c_k.
\]
Define
\begin{equation}\label{eq:bh2-defs}
    H_{\mathrm{eff}} := (1 - q\beta_0)H - q\alpha_0(b-1), \qquad \psi_k := c_k - q\alpha_0(1 - s_k), \qquad \mathrm{req}_k := \max(0,\, H_{\mathrm{eff}} + \psi_k).
\end{equation}
Then every weight $\eta_{\mathrm{bh},k}\ge 0$ adds $+2q\,\eta_{\mathrm{bh},k}$ to the right-hand demand of~\eqref{eq:unrestricted-cap} with \emph{zero} $|q-v_y|$ penalty, requiring interior gap mass $4\,\mathrm{req}_k\,\eta_{\mathrm{bh},k}$ on $J=[-b+1,b-1]$ and wall-footprint weight $4\,\eta_{\mathrm{bh},k} f_{\mathrm{bh},k}(b-|X|)|v_x|$, where $f_{\mathrm{bh},k}(u) := \max\!\bigl(0,\frac{u}{2\sqrt{1-u^2}}-\frac{q\alpha_0}{2}\bigr)\mathbf{1}_{(0,s_k)}(u)$. Moreover, on every wall-optimal four-leg bowtie ($q>1/3$), choosing $(\alpha_0,\beta_0)=(a/\ell, h/\ell)$ achieves exact algebraic equality in the inter-wall chord step.
\end{theorem}
\begin{proof}
Decompose $[0,L]$ ($qL > 2H+4$, so both walls must be visited) into alternating left/right wall-labelled blocks $B_1,\dots,B_{2N}$ ($N\ge 1$) inside disjoint intervals $[\tau_i^-,\tau_i^+]$, separated by full crossings $[\tau_i^+,\tau_{i+1}^-]$ of $J=[-b+1,b-1]$ (retain an even run of complete blocks with one further block supplying the last exit-to-entry interval). Every interval of length $>2(H+1)/q$ visits both closed footprints: use on each wall the phase $q$ times that interval's midpoint. Hence the initial and final pieces, including at most one discarded complete block to ensure even parity, have uniformly bounded length for fixed $(q,H,b)$. Their phase spans and all endpoint terms below are $O_{q,H,b}(1)$. In particular $t_{2N+1}^- - t_1^-= L-O_{q,H,b}(1)$. On each block $B_i$, let $t_i^-, t_i^+ \in B_i\cap\{|X|\ge b-1\}$ minimize $qt - Y(t) - \sqrt{1-(b-|X(t)|)^2}$ and maximize $qt - Y(t) + \sqrt{1-(b-|X(t)|)^2}$, and write $u_i^\pm := b - |X(t_i^\pm)| \in [0,1]$. Every phase covered on $B_i$'s wall during $B_i$ lies in an interval of length $I_i = \sqrt{1-(u_i^-)^2} + \sqrt{1-(u_i^+)^2} + q(t_i^+-t_i^-) + (Y(t_i^-)-Y(t_i^+))$. Using the telescoping identity $\sum_{i=1}^{2N}(t_i^+-t_i^-) = (t_{2N+1}^- - t_1^-) - \sum_{i=1}^{2N}(t_{i+1}^--t_i^+)$, full two-wall coverage $\sum_{i=1}^{2N} I_i \ge 2qL - O(H)$ rearranges (after subtracting $qL$) to
\[
    qL - O(H+b) \le \sum_{i=1}^{2N} \Bigl[\sqrt{1-(u_i^-)^2} + \sqrt{1-(u_i^+)^2} + \bigl(Y(t_i^-)-Y(t_i^+)\bigr) - q(t_{i+1}^- - t_i^+)\Bigr].
\]
Since $t_i^+ \le \tau_i^+ < \tau_{i+1}^- \le t_{i+1}^-$ and $X(t_i^+), X(t_{i+1}^-)$ lie on opposite walls, $|X(t_{i+1}^-)-X(t_i^+)| = 2(b-1) + (1-u_i^+) + (1-u_{i+1}^-)$. By unit guard speed and Cauchy--Schwarz with $(\alpha_0,\beta_0)$,
\[
    t_{i+1}^- - t_i^+ \ge \alpha_0\bigl(2(b-1) + (1-u_i^+) + (1-u_{i+1}^-)\bigr) + \beta_0\bigl(Y(t_{i+1}^-) - Y(t_i^+)\bigr).
\]
Summing over $i=1,\dots,2N$, $\sum_i (1-u_{i+1}^-) = \sum_i (1-u_i^-) + O(1)$ and $\sum_i (Y(t_{i+1}^-)-Y(t_i^+)) = \sum_i (Y(t_i^-)-Y(t_i^+)) + O(H)$; since $1-q\beta_0 > 0$ and $Y(t_i^-)-Y(t_i^+)\le 2H$, multiplying by $2$ gives
\[
    2qL - O(H+b) \le \sum_{i=1}^{2N} \Bigl[ 2F(u_i^-) + 2F(u_i^+) + 4H_{\mathrm{eff}} \Bigr], \qquad F(u) := \sqrt{1-u^2} - q\alpha_0(1-u).
\]
Now use the side condition $s_k \ge q\alpha_0 c_k$: for every $u\in[s_k, 1]$, $\sqrt{1-u^2} \le c_k$, so
\[
    F(u) - F(s_k) = (u - s_k)\!\left(q\alpha_0 - \frac{u + s_k}{\sqrt{1-u^2} + c_k}\right) \le (u - s_k)\!\left(q\alpha_0 - \frac{2s_k}{2c_k}\right) \le 0,
\]
whence $F(u) \le F(s_k) = \psi_k$ on $[s_k,1]$ without any footprint term. For $u \in [0, s_k)$, $-F'(w) = \frac{w}{\sqrt{1-w^2}} - q\alpha_0 \le 2 f_{\mathrm{bh},k}(w)$, so $F(u) \le \psi_k + 2\int_u^{s_k} f_{\mathrm{bh},k}(w)\,dw$. Letting $u_i^* := \min_{t\in B_i\cap\{|X|\ge b-1\}}(b-|X(t)|) \le \min(u_i^-, u_i^+)$, we obtain $2F(u_i^-) + 2F(u_i^+) + 4H_{\mathrm{eff}} \le 4\,\mathrm{req}_k + 8\int_{\min(u_i^*,s_k)}^{s_k} f_{\mathrm{bh},k}(w)\,dw$. On the single excursion inside $B_i$ attaining $u_i^*$, $b-|X|$ traverses $[\min(u_i^*,s_k), s_k]$ twice, contributing $\int_{\tau_i^-}^{\tau_i^+} 4 f_{\mathrm{bh},k}(b-|X|)|X'|\,dt \ge 8\int_{\min(u_i^*,s_k)}^{s_k} f_{\mathrm{bh},k}(w)\,dw$, while each gap crossing $[\tau_i^+,\tau_{i+1}^-]$ contributes $\int_{\tau_i^+}^{\tau_{i+1}^-}\rho_{\mathrm{bh},k}(X)|X'|\,dt \ge 4\,\mathrm{req}_k$. Dividing by $L\to\infty$ proves the $+2q\,\eta_{\mathrm{bh},k}$ demand, and on any bowtie with $(\alpha_0,\beta_0)=(a/\ell,h/\ell)$, Cauchy--Schwarz along $(2a,2h)$ holds with equality.
\end{proof}

At $q=1/4$ the five bounded-excursion certificates and the two analytic inequalities give the following ceilings on $\Wopt(1/4,H)/r$. The decimal entries are rounded upward.
\[
\begin{array}{l|ccccc}
H/r&4/5&49/50&1&5/4&3/2\\\hline
\text{flux certificate}&36/5&153/20&77/10&79/10&399/50\\
\text{Theorem~\ref{thm:midline}}&6.873&7.413&7.473&8.223&8.973\\
\text{Theorem~\ref{thm:fixed-height}}&6.654522&7.111252&7.161548&7.783853&8.396359
\end{array}
\]
The last row uses $u=2/15,0.14961,0.15127,0.16962,0.18456$, respectively, in~\eqref{eq:height-dual}. It gives the strongest of these three ceilings at the first four heights. At the first height the equality curve also gives the exact optimum. At $H=3r/2$ the flux certificate is strongest and excludes width $8r$ for all $H\le3r/2$ by monotonicity. For arbitrary bounded excursion the unrestricted certificate above gives $8.10r$.

\subsection{The cost of triangulation at the critical speed}\label{sec:face-rate}

The optimal guard at $q_a$ has four straight legs. Its separator in Theorem~\ref{thm:obstacle-criterion} follows a circular arc during each downward return. How many planar faces are needed to certify a width close to this optimum?

Normalize $r=V=1$. Let $N(\varepsilon)$ be the smallest number of triangles in a certificate of Proposition~\ref{prop:surface-cut} for a full width at least $F(q_a)-\varepsilon$ at speed $q_a$. The minimum ranges over all periodic polygonal patrols and all allowed real vertices, periods and affine-leg witness times. Uniform contraction of both space and time certifies any narrower width with the same number of faces.

\begin{theorem}[Sharp triangle count at the critical speed]\label{thm:face-rate}
As $\varepsilon\downarrow0$,
\[
 N(\varepsilon)=\Theta(\varepsilon^{-1/4}).
\]
For every integer $n\ge1$ there is a four-leg patrol with a certificate of at most $20n+46$ faces and a full width $W_n$ satisfying
\[
 0<F(q_a)-W_n<6.107\times10^{-10}n^{-4}.
\]
Conversely, there are absolute constants $c_0,\varepsilon_0>0$ such that every $N$-face certificate with width deficit $\varepsilon\le\varepsilon_0$ satisfies $\varepsilon\ge c_0/(N+1)^4$. In physical units the width deficits are multiplied by $r$.
\end{theorem}

In this statement $\Theta(\cdot)$ denotes equality of orders of magnitude. In the rest of the subsection the letter $\Theta$ is the terminal fan angle defined below.

For example, 66 faces suffice to get within $10^{-9}r$, 146 within $10^{-12}r$, and 606 within $10^{-15}r$. The ancillary 66-face example consists of 32 causal faces and 34 capture faces. Its rational coordinates certify width greater than $3.36785112455336r$ at $q=0.457892036651454584>q_a$. The two critical-root sign tests and the exact surface check give a deficit below $5.952\times10^{-10}r$ at $q_a$. Monotonicity in speed allows this example to use rational data even though the theorem is stated at the critical root.

\subsubsection{Travel time along chords}
Take $q=q_a$ and use the critical bowtie of~\eqref{eq:bowtie-params} and~\eqref{eq:diffraction-defs}. Its half-amplitude and height are $a,h$, its diagonal half-length is $\ell$, and $c=a/\ell$, $d=h/\ell$, $\alpha=\sqrt{1-q^2}$, $m=qd+\alpha c$, $\kappa=\alpha d-qc$. Its wall chord has half-height $G=2q\ell-(1-q)h$ and inset $\delta=\sqrt{1-G^2}$. The full width is $2(a+\delta)$ and $a=G^2/\delta$. Put
\[
 \Theta=\arctan\frac{q\kappa}{1+qm},\quad C=\cos\Theta,\quad S=\sin\Theta,\quad
 f(\theta)=\frac{1}{\cos\theta+\sqrt{q^2-\sin^2\theta}},\qquad
 p=\int_0^\Theta f(\theta)\,d\theta .
\]
The critical handoff equality is
\begin{equation}\label{eq:face-critical}
 h=p+\frac{1+\kappa C-mS}{m+q}.
\end{equation}
During a descending return centered at time zero the guard is $(a,-t)$. In the guard frame write $z(\theta)=(-\cos\theta,\sin\theta)$. The velocity of a speed-$q$ intruder relative to the guard along the tangent at angle $\theta$ is $L(\sin\theta,\cos\theta)$, where $L=1/f(\theta)$ is the larger root of $|L(\sin\theta,\cos\theta)-(0,1)|=q$. Thus $p$ is the contact-arc travel time. A chord joining angles $u,v$ has direction $(\sin\mu,\cos\mu)$ and length $2\sin(l/2)$, where $l=v-u$ and $\mu=(u+v)/2$. Its fastest relative travel time is $2\sin(l/2)f(\mu)$.

\begin{lemma}[Polygonal chord-time deficit]\label{lem:face-chord}
For any partition of $[0,\Theta]$ with $k$ gaps $l_i$ and midpoints $\mu_i$ let $P=\sum_i 2\sin(l_i/2)f(\mu_i)$ and $E=p-P$. With $g=f''+f$,
\[
 g(\theta)=\frac{q^2}{(q^2-\sin^2\theta)^{3/2}},\qquad
 \frac{1-\Theta^2/80}{24q}\sum_i l_i^3
 \ \le E\le\frac{\max g}{24}\sum_i l_i^3 .
\]
For $n$ equal angular gaps,
\[
 n^2 E \longrightarrow L_0=\frac{\Theta^2}{24}\int_0^\Theta g(\theta)\,d\theta .
\]
\end{lemma}
\begin{proof}
Rationalize $f=(\cos\theta-\sqrt{q^2-\sin^2\theta})/(1-q^2)$ and differentiate. On either half of a gap the difference between $f$ and the solution of $y''+y=0$ matching its value and derivative at the midpoint is the integral of $g$ with nonnegative sine kernel. Integration over both halves bounds the gap error between $\min g$ and $\max g$ on the gap times $l-2\sin(l/2)$. Here $\min g\ge 1/q$ and $l^3(1-l^2/80)/24\le l-2\sin(l/2)\le l^3/24$. A Riemann sum gives the stated limit.
\end{proof}

\subsubsection{A finite fan and its certified width}
Use equal angular gaps and the $P,E$ of Lemma~\ref{lem:face-chord}. Set
\[
 h_n=h-E,\quad \lambda=h_n/h,\quad a_n=\lambda a,\quad
 \ell_n=\lambda\ell,\quad G_n=\lambda G,\quad
 \delta_n=\sqrt{1-G_n^2},\quad b_n=a_n+\delta_n .
\]
In particular $h_n-P=h-p$. The wall chord is tight, and the terminal polygonal handoff is tight by \eqref{eq:face-critical} with $(h,p)$ replaced by $(h_n,P)$.

\paragraph{The outgoing return.}
Traverse the successive chords on the unit disk in the return frame at their fastest speeds, starting at $z(0)$ at time zero. In the laboratory frame chord $i$ has velocity $q e_i$ with $e_i=(\sin\psi_i,\cos\psi_i)$. These directions turn monotonically from the vertical toward $(\kappa,m)$, staying strictly between the two endpoint directions. Indeed $\tan\theta=q\sin\psi/(1+q\cos\psi)$ is strictly increasing on the angular range, and the relevant $\theta$ is the chord midpoint. Let $e_0=(0,1)$ for the initial flat. Each contact segment belongs to a causal plane
\[
 e_i\cdot(X-X_i)=q(t-t_i),
\]
where $(t_i,X_i)$ is its initial event. At a change of direction extend the two planes along their line of intersection into future time. Its spatial velocity is $q(e_{i-1}+e_i)/(1+e_{i-1}\cdot e_i)$. Its transverse component is strictly between those of $qe_{i-1}$ and $qe_i$ (equivalently it is the secant slope of the strictly convex function $q\sqrt{1+s^2}$ at the two slopes $s=\tan\psi$). Thus each intersection starts at the current contact, stays to its left, and stays to the right of the preceding intersection. These lines and the contact polygon bound the $n$ planar pieces of the fan. Its slice at $P$ is a continuous polyline with slopes between $0$ and $-\kappa/m$.

Extend the flat $y=qt$ to the left wall. On the right join contact to $(b_n,-t-G_n)$ in the return tube. Triangulate each strip between consecutive contact times. At each corner the actual time is a witness on this return leg. Convexity of its spacetime tube gives both triangles even when their four corners are not coplanar. For the incoming half reflect $t,y$ and use the wall offset $+G_n$. At the center a single captured triangle switches that offset to $-G_n$. Reverse the vertex order under this $(t,y)$ reflection so that each face of nonzero projected area has $n_y<0$. The reflected causal planes then have $n_t>0$.

\paragraph{The diagonal handoff.}
At time $P$, use the local coordinates $(x,y)$ relative to the next departure position $(a_n,-h_n)$, and $v=h_n-P$. The fan begins at $(-1,h_n+qP)$ and ends at $(-C,v+S)$. The cost to the diagonal tube for the witness $U=(-mx+\kappa y)/\alpha$ equals
\[
 qv+\kappa x+my
 \quad\text{if }0\le U\le2\ell_n\text{ and }dx+cy\ge0 .
\]
Along the terminal polyline this cost and $dx+cy$ are nondecreasing, and $U$ is nonincreasing, by the slope bounds above and $dm-c\kappa=q$. The terminal cost equals one by the discrete handoff equality. The two endpoint tests for $U$ and $dx+cy$ suffice for the whole polyline. The remaining strict conditions from Theorem~\ref{thm:obstacle-criterion} are $0<P<h_n$, $a_n+b_n\ge1$, $qv\le1$, $\alpha^2v\le2(G_n-q)$ and $0<G_n<1$. The companion interval checker verifies these for every $n\ge1$ using $0\le E\le (\max g)\Theta^3/24$. Consequently the wall endpoint and whole terminal polyline admit witnesses on the same diagonal. The terminal edges include the flat segment from $x=a_n-1$ to the first intersection ray. Join each edge to the wall point to reset to the endpoint chord. Convexity of the diagonal tube supplies the capture witnesses.

The bridge to the next return is the affine bridge in the proof of Theorem~\ref{thm:obstacle-criterion} with $P$ in place of $p$. Explicitly, in departure time $t$, join $(a_n-1,qP)$ at $t=-v$ linearly to $(-a_n,h_n+G_n)$ at $t=2\ell_n$, giving a point $U_+(t)$ on the upper flat $y=q(h_n+t)$. Use its time-space reflection $-U_+(2\ell_n-t)$ for the lower flat. Both endpoints lie in the same convex diagonal tube whenever both are used. From $-v$ to departure the upper point joins the wall point. From departure to arrival the two points join each other. The arrival end is the reflection of the departure end. Let $x_U,x_L$ be the transverse coordinates of the upper and lower points. During the diagonal both are affine with the same slope and $x_U(t)-x_L(t)=x_U(0)-a_n<0$, so the two flats and the strip between them are ordered. The two resets at departure and arrival are captured triangles. Reflect horizontally for the second bridge. The flat, fan and captured strip projections partition the $(t,x)$ rectangle with downward orientation. The instantaneous resets have zero projected area. The boundary matches modulo the period away from the two walls, so this is a certificate in Proposition~\ref{prop:surface-cut}.

\paragraph{Face count and deficit.}
Per half-fan there are $2n-1$ active-plane triangles, three triangles on the rising flat (split at $x=a_n-1$), $2n$ on the contact-to-wall strip and $n+1$ on the handoff reset. There are four half-fans, two central wall resets and two bridges of 16 faces each. A bridge uses four triangles before departure, one at the departure reset, six during the diagonal, one at the arrival reset and four after arrival. This gives $20n+46$. Write $B(\lambda)=\lambda a+\sqrt{1-\lambda^2G^2}$. By $a=G^2/\delta$,
\[
 B'(1)=0,\qquad -B''(\lambda)=G^2(1-\lambda^2G^2)^{-3/2}.
\]
In particular $W_n=2B(\lambda)$ and
\[
 F(q_a)-W_n \le \frac{G^2}{\delta^3h^2}E^2
 \le \frac{G^2}{\delta^3h^2}
       \left(\frac{\Theta^3\max g}{24}\right)^2 n^{-4},\qquad
 \lim_{n\to\infty}n^4(F(q_a)-W_n)=\frac{G^2L_0^2}{\delta^3h^2}.
\]
The exact interval checks in \path{anc/faces/check_fan_face_rate.py} bound the coefficient in the inequality by $6.107\times10^{-10}$ and enclose the limiting coefficient in $(5.96482285,5.96482288)\times10^{-10}$. This limit concerns equal angular gaps in the displayed four-leg family.

\subsubsection{A lower bound for every certified patrol}
There is a geometric reason more general patrols, auxiliary vertices, and time-shifted witnesses cannot improve the exponent. A certificate with few capture vertices also works with part of the sensing disk removed.

\begin{lemma}[A smaller sensing set from the witnesses]\label{lem:witness-set}
Let a polygonal patrol admit an $N$-face unit-radius certificate. Define $D_0=\{v:|v|\le1,\ |v_x|\le1/3\}$. There is a convex set $K\subseteq\overline B_1(0)$, the convex hull of $D_0$ and at most $12N$ points, invariant under both coordinate sign changes, with the following property. Every full transit with velocity in $qD_0$ meets $g(t)+K$ at some time.
\end{lemma}
\begin{proof}
For each captured face and each of its three vertex occurrences let $r_i=q|w_i-t_i|$. When $r_i<1$ retain $z_i=(P_i-g(w_i))/(1-r_i)$. When $r_i=1$ take $z_i=0$, since $P_i=g(w_i)$. All $z_i$ are in the unit disk. Adjoin their coordinate reflections to $D_0$ and take the convex hull. On a face with barycentric coefficients $\lambda_i$, put $t=\sum\lambda_i t_i$ and $w=\sum\lambda_i w_i$. The common affine leg implies
\[
 p(w)-g(w)=\sum_i\lambda_i(1-r_i)z_i+[p(w)-p(t)]
\]
for a transit intersecting the face at $(t,p(t))$. Central symmetry and convexity give $p(w)-p(t)\in q|w-t|D_0\subseteq (\sum_i\lambda_i r_i)D_0$. The coefficients of the retained points total $1-\sum_i\lambda_i r_i$ and the remaining term lies in $(\sum_i\lambda_i r_i)D_0$. Hence the right side belongs to $K$.

For completeness, if a transit avoids $g+K$, trim it between levels below and above the bounded surface, with a buffer of $1/q+1$ for all witnesses. On this compact window it has positive distance from $g+K$. A sufficiently small contraction of both spatial coordinates toward $(0,0)$ in the centered strip creates uniform velocity slack in $qD_0$ and a wall margin. Polygonal interpolation and a small generic perturbation preserve both margins and avoidance. The crossing has signed intersection number $-1$ by the degree and boundary conditions of Proposition~\ref{prop:surface-cut}. A causal face has nonnegative pairing with every allowed velocity, so a transversely negative intersection must lie on a captured face. The preceding calculation contradicts avoidance. This is the soundness argument with its velocity and sensor sets restricted together.
\end{proof}

\begin{lemma}[An escape through the missing circular caps]\label{lem:cap-escape}
For the canonical bowtie at $q_a$ and each set $K$ constructed in Lemma~\ref{lem:witness-set} there exists a $qD_0$ transit whose distance from $g(t)+K$ is uniformly at least $\gamma/(3N+1)^2$, where $\gamma=1.37\times10^{-5}$. Its nonvertical portion and the range of its phase $qt-p_y(t)$ can be chosen uniformly bounded, independently of $N$ and $K$.
\end{lemma}
\begin{proof}
Partition the left return-frame circular arc $z(\theta)=(-\cos\theta,\sin\theta)$, $0\le\theta\le\Theta$, at the directions of all nonzero point generators whose radial projections lie on the arc, including generators strictly inside the disk. Add both endpoints. Each orbit under the two coordinate sign changes contributes at most one direction to this open quadrant. The at most $3N$ original witness occurrences therefore leave $k\le3N+1$ gaps. On each gap of length $l$ the whole set $K$ is below the chord in its outward normal direction $z(\mu)$:
\[
 \sup_{z\in K} z(\mu)\cdot z\le\cos(l/2).
\]
Indeed no retained direction is interior to the gap, and $z(\mu)\cdot v\le 1/3+\sin\Theta<\cos\Theta$ on $D_0$. This inequality persists under convex combinations.

Write $P$ for the chord time and $E=p-P$. Put
\[
 R=\frac{h(m+q)}{1+\kappa C-mS+(m+q)P}=\frac h{h-E}>1 .
\]
During the right downward return follow the \emph{dilated} chords $R z(\theta)$ in the guard frame, starting at $(a-R,0)$ in laboratory coordinates. A chord takes $R$ times its unit-scale travel time, so the terminal time is $\tau=RP$. Its distance from $K$ is at least $(R-1)\cos(l/2)$ by projection on the chord normal. The absolute chord velocities are $q e_i$ with $0\le(e_i)_x\le\kappa<1/3$.

Continue with absolute velocity $(u,w)=q(\kappa,m)$, tangent to the terminal circle in the return frame. On the remainder of that leg the Euclidean distance to the guard is at least $R$. On the next diagonal the unit normal $(\kappa,m)$ has guard and evader normal velocity $q$. The signed separation there is $h(m+q)-R[\kappa C-mS+(m+q)P]=R$. The tangent line has height $wt+B$ with $B=R S-(1+w)\tau$. At $T_s=(2.1-B)/w$ freeze the transverse coordinate, then continue vertically at speed $q$. Before this join on the left return the horizontal separation is at least $2a-R>R$; after it the vertical separation is at least $2.1-h>R$. The conditions just used, including $\tau<h$, $h+2\ell<T_s<3h+2\ell$, and both strip endpoints for the path, hold for every partition. They are strict interval checks using $E\le(\max g)\Theta^3/24$. Reflect in time and height for negative times, using the symmetry of $K$. The vertical tails can start at bounded uniform levels beyond the guard. This gives uniform time and phase bounds for the middle piece.

Since $R-1\ge E/h$, Lemma~\ref{lem:face-chord} and Jensen's inequality $\sum l_i^3\ge\Theta^3/k^2$ bound the global clearance from below by
\[
 \frac{\cos\Theta}{h}\frac{(1-\Theta^2/80)\Theta^3}{24q k^2}
 >\frac{\gamma}{k^2}.
\]
\end{proof}

\paragraph{Transferring the escape to the certified patrol.}
Write $W=F(q_a)-2\Delta$. A valid certificate implies full guarding by the original unit disk, so the period-uniform estimate of Lemma~\ref{lem:gap-crossing} applies to the original patrol with $J=\{q_a\}$, independently of $K$. Let $T_*$ be the canonical period. Lemma~\ref{lem:gap-crossing} starts its block at a left-return midpoint. Place the right-return midpoint of Lemma~\ref{lem:cap-escape} at $T_*/2+jT_*$ in that block. Increasing $j$ by one shifts the escape by $T_*$ in time and its phase by $qT_*$. The uniform time and phase bounds allow one choice of $j$ and a fixed block length that put the middle segment inside the block and the entire phase range more than two units from its ideal phase boundaries. After time and height translations there is a crossing on which the physical patrol stays within $C_0\sqrt{\Delta}$ of that canonical block. Its time at phases in the interior of the block's phase slab outside the crossing is at most $C_0\sqrt\Delta$. Here $C_0$ and the smallness threshold on $\Delta$ are independent of the certificate, its period and $N$.

Shift the evader accordingly and clamp its transverse coordinate to the narrower actual strip. The cost is at most $\Delta$ and the velocity remains in $qD_0$. On the crossing its clearance from the translated sensing set drops by at most $C_0\sqrt\Delta+\Delta$. At other times a unit enlargement of the evader's phase range lies in a fixed interior core of the slab. A guard visit to this core would use at least $\rho/(1+q)$ time outside the crossing in the open slab, where $\rho>0$ is the distance of the core from its edges. Indeed $|(qt-g_y)'|\le 1+q$ and the guard phase tends to the appropriate infinity on each tail. For small $\Delta$ this contradicts the outside-time estimate, exactly as in the proof of Theorem~\ref{thm:global-gap}. The vertical separation exceeds one at all these other times. Thus if $C_0\sqrt\Delta+\Delta<\gamma/(3N+1)^2$, there would be a $qD_0$ transit avoiding $g+K$. Lemma~\ref{lem:witness-set} forbids it. For $\Delta\le1$, this gives $(C_0+1)\sqrt{\Delta}\ge\gamma/(3N+1)^2$. Squaring and using $\varepsilon=2\Delta$ proves the lower bound of Theorem~\ref{thm:face-rate}.

\section{Slow intruders (\texorpdfstring{$q \downarrow 0$}{q -> 0})}\label{sec:asymptotics}
We turn next to the slow-intruder limit $q\downarrow 0$ in the known-direction game, establishing the leading $2r/q$ width together with the exact fixed-excursion scaling $\min\{2,1+h\}$ before refining the unrestricted additive deficit $D_1(q)=2/q-\Wopt(q)/r$.

\begin{theorem}[Leading slow-speed width and fixed-excursion scaling]\label{thm:slow-limit}
As $q\downarrow0$, $q\Wopt(q)/r\to2$, and for sufficiently small $q>0$ a four-leg bowtie guards width $(2/q-2+3q+O(q^2))r$. Under any fixed excursion cap $H=hr$ ($h\ge 0$),
\[
    \lim_{q\downarrow 0}\frac{q\Wopt(q,hr)}{r}=\min\{2,\,1+h\}.
\]
\end{theorem}
\begin{proof}
Normalize $r=V=1$. Set $a=q^{-1}-1-q/2$ and $h_*=1-q+q^2-3q^3$. For small $q$ these quantities are positive and $(1+q)h_*=1-2q^3-3q^4<1$. Writing $\ell=\sqrt{a^2+h_*^2}=q^{-1}-1+O(q^2)$ gives
\[
 G=2q\ell-(1-q)h_*=1-2q^2+O(q^3),\qquad
 \delta=\sqrt{1-G^2}=2q+O(q^2).
\]
Theorem~\ref{thm:two-level} certifies width $2(a+\delta)=2/q-2+3q+O(q^2)$. Washburn's area bound~\eqref{eq:area-bound} gives $q\Wopt(q)\le2\sqrt{1+q^2}$, proving $q\Wopt(q)/r\to 2$.

For a fixed cap $h\ge 0$ in normalized units,~\eqref{eq:midline-bound} and~\eqref{eq:area-bound} give the upper limit $\min\{2,1+h\}$. When $h\ge 1$, the four-leg construction above has half-height $h_*<1\le h$ and achieves leading coefficient $2$. When $0<h<1$, take $\ell=[1+(1-q)h]/(2q)$, $a=\sqrt{\ell^2-h^2}$, and $\delta=0$; once $(1+q)h\le 1$, $G=2q\ell-(1-q)h=1$, so Theorem~\ref{thm:two-level} certifies guarded width $2a$ with $2qa\to 1+h$, while $h=0$ is immediate from Theorem~\ref{thm:midline}(1).
\end{proof}

\subsection{An explicit ten-leg patrol at finite speed}

The four-leg family in Theorem~\ref{thm:slow-limit} has limiting additive loss $2r$ relative to $2r/q$. A short inward kink in each return and two extra corners on each transfer give a smaller loss throughout a fixed speed interval.

\begin{theorem}[Explicit ten-leg slow-speed patrol]\label{thm:slow-hooks}
For every $0<q\le 1/40$, a periodic ten-leg patrol of period $4r/U$ guards a width greater than $(2/q-1.235352896057745)r$.
\end{theorem}
\begin{proof}
Normalize $r=V=1$. Use the exact rationals
\[
 h=1.227689517575272716605,\qquad e=0.029002301252121295625.
\]
Set $L=1+e$, $\rho=\sqrt{h^2+e^2}$, $\alpha=\sqrt{1-q^2}$, and define
\begin{align*}
 k&=h-1+q\rho, &
 v&=\frac{k+q\sqrt{k^2+\alpha^2 L^2}}{\alpha^2}, &
 \tau&=\sqrt{L^2+v^2},\\
 Y&=h-v, &
 x&=\alpha Y/q, &
 a&=x+L,\\
 G&=2-h-q\rho, &
 \delta&=\sqrt{1-G^2}, &
 b&=a+\delta.
\end{align*}
The right return is $(a,h)\to(x+1,0)\to(a,-h)$, followed by the transfer $(a,-h)\to(x,-Y)\to(-x,Y)\to(-a,h)$. Reflect these five unit-speed legs horizontally to close the patrol. Write $D=2(\tau+Y/q)$ for the transfer duration. The stated root of the quadratic for $v$ gives $v=k+q\tau$ (the unsquared sign is positive). Consequently
\[
 Y+q(\rho+\tau)=1,\qquad G=q(D+\rho)-h,\qquad 4\rho+2D=4/q.
\]
The middle leg has velocity $(-\alpha,q)$. Each return leg has duration $\rho$ and the two outer transfer legs have duration $\tau$. All vertices lie in $[-b,b]\times\R$. The parameters satisfy $v,Y>0$, $0<G<1$ and, for $R=\sqrt{1-(\delta+e)^2}$,
\begin{equation}\label{eq:hook-wall}
 R>0,\qquad 1-q\rho>0,\qquad \delta^2+(h-R)^2<(1-q\rho)^2 .
\end{equation}
We give the uniform arithmetic check below.

Set the departure time from $(a,-h)$ to zero. If the guard leg on $J$ is affine, the set of spacetime points $(t,p)$ for which some $w\in J$ satisfies $\|p-g(w)\|+q|w-t|\le1$ is convex, by interpolating both points and witnesses. Every intruder occupancy of this set causes detection. We construct the separator from two horizontal levels
\[
 F_+(t)=q(t+\rho),\qquad F_-(t)=qt-h-G.
\]
At the right-return midpoint $t=-\rho$, the anchor $K=(x,0)$ and both wall points $(b,\pm R)$ lie in the current disk. The anchor also belongs to the capture set of the first transfer leg via its corner at time $\tau$, since $Y+q(\rho+\tau)=1$. The lower wall point belongs to the same set via the departure at time zero, by \eqref{eq:hook-wall}.

Let $U(t)$ interpolate affinely in spacetime through the successive points
\[
 (-\rho,K),\quad(\tau,(x,F_+(\tau))),\quad
 (D-\tau,(-x,F_+(D-\tau))),\quad(D,(-a,h+G)).
\]
Its height is $F_+(t)$. At both transfer corners the upper and lower levels have offsets $+1$ and $-1$ from the guard, by the timing identities. At arrival the offset of $U(D)$ is $G$. Thus each interpolated segment lies in the capture set of its own transfer leg; the first such segment starts at $t=-\rho$. For $0\le t\le D$ the central symmetry of the transfer gives $B(t)=-U(D-t)$ at height $F_-(t)$ in the same leg's capture set as $U(t)$. Their horizontal coordinates satisfy $U_x(t)\le B_x(t)$: on the first leg the difference is affine, nonpositive at departure and zero at the corner, on the middle leg it vanishes, and on the last leg the ordering follows by symmetry. Use the upper flat from the left wall to $U(t)$, the captured chord $U(t)B(t)$, and the lower flat from $B(t)$ to the right wall.

Before departure, join the upper flat at $U(t)$ to the wall trace interpolating from $(-\rho,(b,-R))$ to $(0,(b,-h-G))$. Both endpoints and the joining chord belong to the first leg's capture set, using $\delta^2+G^2=1$ at time zero. The change at departure sweeps only the captured triangle with vertices $U(0),(a,-h-G),(b,-h-G)$. The reflection $(t,p)\mapsto(D-t,-p)$ gives the corresponding chord and reset at arrival. At the next return midpoint, the change from the upper to the lower wall endpoint sweeps the triangle between the reflected anchor and the two wall points, all in the current disk. Horizontal reflection now gives the second half-period. This bounded separator spans the corridor. Away from captured chords, endpoints, and reset triangles it is horizontal and rises at speed $q$, so the difference between intruder height and front height cannot increase there. Every upward transit must cross the separator and is detected.

For the width estimate, put $E(q)=b-1/q$. The timing identity gives the regular expression
\[
 E=-q/(1+\alpha)-\alpha(\rho+\tau)+L+\delta .
\]
At zero, $-2E(0)=2(\rho+\sqrt{L^2+(h-1)^2}-L-\sqrt{1-(2-h)^2})$. Differentiation gives $\tau'=v(\rho+\tau)/(\tau-qv)$ and
\[
 E'=-\frac1{\alpha(1+\alpha)}+\frac{q(\rho+\tau)}{\alpha}
       -\alpha\tau'+\frac{\rho G}{\delta}.
\]
The standalone \path{anc/hooks/check_slow_hooks.py} encloses these radicals using rational intervals and integer square roots. Direct evaluation on the single interval $[0,1/40]$ certifies the positivity conditions above, $\tau-qv>0$, $1/5<E'<4/5$, and $-2E(0)<1.235352896057745$. Only the desingularized expressions are evaluated at zero. Since $2b=2/q+2E(q)>2/q+2E(0)$ for $q>0$, rescaling proves the theorem.
\end{proof}

\subsection{The best ordinary bowtie and the unrestricted additive loss}
The ordinary class of Theorem~\ref{thm:ordinary-transition} has a
sharp limiting loss of two capture radii.
\begin{theorem}[Ordinary bowties at slow speed]\label{thm:slow-ordinary}
For all sufficiently small $q>0$, the ordinary supremum is attained by
a centered, descending-stroke bowtie of unique dimensions. Up to
translations in time and longitudinal position and a horizontal reflection,
its shape is specified by the analytic maximizer in Appendix~\ref{app:slow-ordinary}. As $q\downarrow0$,
\[
 \frac{W_4(q)}r=\frac2q-2+3q-\frac94 q^3-\frac{19}4 q^4+O(q^5).
\]
In particular, for the optimum over all bounded-excursion patrols,
\[
 \liminf_{q\downarrow0}\frac{\Wopt(q)-W_4(q)}r
 \ \ge\ 0.764647103942255 .
\]
\end{theorem}

Appendix~\ref{app:slow-ordinary} locates all competitive geometries
using interior lane phases and optimizes the exact finite-leg contact
constraint at the resulting angles. For unrestricted bounded-excursion
patrols, Theorem~\ref{thm:slow-hooks} gives $D_1(q)<1.235352896057745$
throughout $0<q\le1/40$, where $D_1(q)=2/q-\Wopt(q)/r$, while the next
theorem supplies a strictly positive floor.

\begin{theorem}[Positive unrestricted one-way loss]\label{thm:loss}
For all $0<q\le2^{-24}$,
\[
                  D_1(q)\ \ge\ 2^{-21}.
\]
In particular the unrestricted slow one-way loss has strictly positive
lower limit. The constants are conservative. The value and existence
of a limiting loss remain open.
\end{theorem}

The proof in Appendix~\ref{app:oneway-loss} uses a disk cap swept into
one transverse band during a turn and a second band which remains
untouched. Their reachable-height difference controls reversals
on a window long enough to compare with the entire corridor width.

\section{Bidirectional guarding}\label{sec:diagnostics}

We now turn to the bidirectional (two-way) game, in which the intruder's
direction of transit is unknown in advance: a single prescribed patrol $g$ on
$\mathcal{C}_b=[-b,b]\times\R$ must simultaneously capture every upward speed-$U$
transit of Definition~\ref{def:model} and every downward speed-$U$ transit
obtained by reflecting such a path across $y=0$. For $q=U/V>0$ and $H\ge 0$,
define the excursion-constrained two-way minimax width
\[
    W_{\mathrm{2way}}(q,H):=\sup\bigl\{\,2b>0:\exists\text{ two-way guard }g\text{ on }\mathcal{C}_b\text{ with }\sup_{t\in\R}|Y(t)|\le H\,\bigr\}
\]
and the unrestricted bounded-excursion two-way minimax width
$W_{\mathrm{2way}}(q):=\sup_{H\ge 0}W_{\mathrm{2way}}(q,H)$. When comparing the
known-direction and bidirectional games directly, we write $W_1:=\Wopt$,
$W_2:=W_{\mathrm{2way}}$, and $D_i(q):=2/q-W_i(q)/r$ ($i=1,2$). By the
finite-state discretization in Theorem~\ref{thm:periodic-approximation}---tracking
in each discrete state both the set reachable from the bottom grid row $B$ (kept
disjoint from $C$) and the set reachable from $C$ (kept disjoint from $B$)---periodic
polygonal patrols achieve $W_{\mathrm{2way}}(q,H)$ at every $q>0$ and $H\ge 0$,
and hence achieve $W_{\mathrm{2way}}(q)$ as well.

\subsection{An excursion bound and optimality of horizontal sweeping for \texorpdfstring{$q \ge 1/2$}{q >= 1/2}}\label{sec:twoway}

\begin{theorem}[Two-way width bound]\label{thm:twoway}
If a patrol of height at most $H$ captures every straight wall runner of
speed $U$ in both directions, its width satisfies
\begin{equation}\label{eq:twoway-bound}
  2b\le W_{\mathrm{bf}}(q)+\frac{(1-2q)_+}{q}H,
  \qquad W_{\mathrm{bf}}(q)=r\sqrt{4+q^{-2}} .
\end{equation}
Here $s_+=\max(s,0)$. In particular, for every $H\ge0$ and $q\ge1/2$,
\begin{equation}\label{eq:twoway-exact}
 W_{\mathrm{2way}}(q,H)=W_{\mathrm{2way}}(q)=W_{\mathrm{bf}}(q).
\end{equation}
\end{theorem}

\begin{proof}
Use $r=V=1$ and put
\[
 b_0=\tfrac12\sqrt{4+q^{-2}},\qquad
 a_0=\frac{1}{2q\sqrt{1+4q^2}},\qquad R_0=2qa_0,\qquad
 \delta_0=b_0-a_0=2qR_0 .
\]
It suffices to treat $\Delta=b-b_0>0$. At a time when the right wall is in
range, the slice of the disk along that wall has vertical half-length
$\sqrt{1-(b-X)^2}$. Its supporting tangent at $b-X=\delta_0$ bounds this
half-length by $R_0-2q(b-X-\delta_0)=2q(X-\Delta)$.

Consider a complete positive excursion with $X(0)=X(D)=0$ and $X>0$ on
$(0,D)$, translating the time origin for this calculation. Call it active if
the disk reaches the right wall at some time. On an active excursion the
upward and downward runners $y=\pm(qt-\zeta)$ detected on that wall give
two sets of phases $\zeta$. Write $I_\pm$ for their respective diameters.
For
\[
 A_\pm(t)=qt\mp Y(t)+2qX(t),\qquad B_\pm(t)=qt\mp Y(t)-2qX(t),
\]
the tangent bound gives $I_\pm\le A_\pm(u_\pm)-B_\pm(v_\pm)-4q\Delta$,
where $u_\pm$ and $v_\pm$ maximize $A_\pm$ and minimize $B_\pm$ on $[0,D]$.

For arbitrary fixed choices of these four times, write
$J(q)=\sum_\pm(A_\pm(u_\pm)-B_\pm(v_\pm))-4qD$. The four marked terms at
$u_+,v_+,u_-,v_-$ have respective pairs of signs
$(+,-),(-,+),(+,+),(-,-)$ on the time and height coordinates. At time $t$,
let $k$ be the number of markers later than $t$, let $m$ sum their time
signs, and let $j$ sum their height signs. All four horizontal coefficients
are $2q$. Integrating up to the markers, using the zero sum of their
height signs and $\int_0^D X'=0$, gives
\[
 J(q)=\int_0^D \bigl(2q(k-2)X'(t)+jY'(t)-q(4-m)\bigr)\,dt .
\]
For $q\ge1/2$, the integrand is nonpositive. For $k=0$ or $4$, the norm
of the velocity coefficient equals $4q=q(4-m)$. For $k=1$ or $3$,
$|m|=|j|=1$ and the norm is $\sqrt{4q^2+1}\le3q\le q(4-m)$.
For $k=2$, the horizontal coefficient vanishes. If $|j|=2$, those two
markers have opposite time signs, so $m=0$ and $|j|\le4q$. The other
possibility is $j=0$. In each case use $\|(X',Y')\|\le1$.
For $0<q<1/2$, $J(q)$ is affine for the fixed marks and path. The bounds
$J(0)\le4H$ and $J(1/2)\le0$ give $J(q)\le4(1-2q)H$. Hence on either
side of the corridor, by reflection in $x=0$,
\[
 I_++I_-\le4qD+4(1-2q)_+H-8q\Delta .
\]

On $[0,L]$ each wall and direction must cover all phases
$[H+2,qL-H-2]$ for large $L$, since those runners can be detected only
inside this window. As $b>1$, detection requires an excursion of the
corresponding sign. Each complete active excursion lasts at least
$2(b-1)$, so their count $N(L)$ is finite. There are at most two truncated
excursions. For each direction a truncated piece of length $d$ covers a
phase span at most $qd+2H+2$. Summing the four phase requirements and
using that the excursion lengths sum to at most $L$ gives
\[
 4qL-O(H+1)\le4qL+
    \bigl(4(1-2q)_+H-8q\Delta\bigr)N(L)+O(H+1).
\]
Arbitrarily late phases require detections on both walls at arbitrarily
late times. The patrol therefore makes infinitely many complete active
excursions, so $N(L)\to\infty$. This proves~\eqref{eq:twoway-bound} upon
rescaling. The horizontal patrol of Corollary~\ref{cor:two-level-conseq}(1)
guards every upward transit. Its reflection in $y=0$ is the same patrol,
so it guards in both directions and proves~\eqref{eq:twoway-exact}.
\end{proof}

At half speed ($q=1/2$), Theorem~\ref{thm:twoway} gives
$W_2(1/2)=W_{\mathrm{bf}}(1/2)=2\sqrt{2}\,r$, whereas Theorem~\ref{thm:qa-exact}
gives $W_1(1/2)=F(1/2)=2\sqrt{7/3}\,r$; knowing the direction of transit thus
multiplies the guardable width by $W_1(1/2)/W_2(1/2)=\sqrt{7/6}$, an $8.01\%$
gain.

\subsection{Sharp gain at small height}

\begin{theorem}[Two-way first variation]\label{thm:twoway-slope}
For each fixed $0<q<1/2$, as $H/r\downarrow0$,
\[
 W_{\mathrm{2way}}(q,H)
 =W_{\mathrm{bf}}(q)+\frac{1-2q}{q}H+O_q(H^2/r).
\]
In particular $W_{\mathrm{2way}}(q,H)>W_{\mathrm{bf}}(q)$ for every $H>0$.
\end{theorem}

\begin{proof}
Normalize $r=V=1$ and take $a=a_0$ from the preceding proof. At unit
speed follow the first four legs
\[
 (-a,-H)\longrightarrow(a,H)\longrightarrow(a,-H)
 \longrightarrow(-a,-H)\longrightarrow(-a,H).
\]
Concatenate their reflection in $y=0$ and repeat. With
$\ell=\sqrt{a^2+H^2}$, the half-period is $D=2\ell+2a+4H$ and the
full period is $T=2D$. Thus
\begin{equation}\label{eq:antipodal}
 (X(t+D),Y(t+D))=(X(t),-Y(t)).
\end{equation}
At either pole $x_0=\pm a$, successive vertical visits last $2H$ and their
centers are $D$ apart. A downstroke alternates with an upstroke. Put
\[
 R=qD/2-H,\qquad \delta=\sqrt{1-R^2},\qquad b=a+\delta,\qquad
 B_d=R+(1+q)H,\qquad B_u=R+(1-q)H .
\]
As $H\downarrow0$, $R\to R_0\in(0,1)$ and $\delta\to\delta_0>0$.
Take $H>0$ small enough that
\[
 0<R<1,\qquad qH<\min(R,\delta),\qquad B_d+H<1 .
\]
The separator below uses only these three inequalities. It applies to the
eight-leg patrol with any half-amplitude $a>0$ and height $H>0$ that
satisfy them, and Proposition~\ref{prop:twoway-half} uses it in that
generality. We construct a bounded separator for every upward transit in
$[-b,b]\times\R$.

Assign a level to each wall. For successive centers $c_j,c_{j+1}$ at its
pole, let $B_j$ be $B_d$ or $B_u$ according to the direction of visit $j$.
Between centers its level is
\[
 f(t)=-B_j+q(t-c_j),\qquad c_j<t<c_{j+1}.
\]
Since $qD=B_d+B_u$, this reaches $+B_{j+1}$ just before the next center.
Both wall levels differ from the guard height by at most $B_d+H<1$.
During transfers and upstrokes use the left wall level for $x<X(t)$ and
the right wall level for $x>X(t)$. Join them by a vertical segment at
$x=X(t)$. That segment lies in the current disk.

During a downstroke centered at $c$, use the other wall's level throughout
the corridor on the inward side of the visited pole. On the outer strip,
between that pole and its wall, join the point at the pole on this level
to the wall point at height $Y(t)-R(t-c)/H$. The chord is inside the
disk. Its wall endpoint has horizontal offset $\delta$ and vertical offset
of magnitude at most $R$. At entry the old strip level is $B_d-qH=H+R$,
the height of the chord's wall endpoint. At exit the new strip level is
$-B_d+qH=-H-R$, again the height of that endpoint. At either join, the
region between the two profiles is the triangle formed by the common
wall endpoint and the two heights at the pole. The disk contains all
three vertices and hence this entire triangle.

At the center $c$ of an upstroke, the level on its outer strip drops from
$+B_u$ to $-B_u$. A speed-$q$ path in this strip at $c$ with
$|y(c)|\le B_u$ is captured during the visit. Indeed, project $y(c)$ onto
$[-H,H]$, call the result $s$, and take $w=c+s$. The guard is then at
$(x_0,s)$. Corridor confinement bounds the path's distance from $x_0$
towards the wall by $\delta$. Its distance on the inward side of $x_0$
is at most $q|w-c|\le qH$, since it started on the outer strip.
Also
\[
 |y(w)-s|\le qH+(|y(c)|-H)_+\le qH+(B_u-H)=R ,
\]
where $B_u-H=R-qH>0$. Therefore its squared distance to the guard
at $w$ is at most $\delta^2+R^2=1$.

To finish, consider a path avoiding capture at all times. On each open
flat, its height relative to the flat is nonincreasing, since
$y'(t)\le q$. At a moving vertical join the path must lie strictly above
both flat endpoints or strictly below both, as the entire join is
captured. The same one-sided classification holds across the chords and
their entry and exit triangles. At an upstroke reset a change from below
to above would place the path in the rectangle just proved captured
during that visit. These local statements preclude a crossing from below
to above on every finite interval. More explicitly, between resets a
point of an uncaptured path on the separator can lie only in the interior
of a flat relative to the corridor, where the same rising level applies
in a time neighborhood. At points strictly off the separator its side is
locally constant. The finitely many joins and resets on a finite interval
preserve this implication. The separator is bounded, so an upward transit
would have to cross it. By~\eqref{eq:antipodal}, a downward transit is
captured as well.

Finally,
$R=q(\ell+a)+(2q-1)H=R_0+(2q-1)H+O_q(H^2)$.
Since $\delta'(R_0)=-R_0/\delta_0=-1/(2q)$, the guarded width $2(a+\delta)$
equals $W_{\mathrm{bf}}(q)+(1-2q)H/q+O_q(H^2)$ in normalized units.
Combining this construction with~\eqref{eq:twoway-bound} and restoring
$r$ gives the expansion. Its positive linear coefficient proves strict
gain for small $H$, and monotonicity in $H$ proves the assertion for all
$H>0$.
\end{proof}

\subsection{Unrestricted height just below \texorpdfstring{$q=1/2$}{q=1/2}}

Theorem~\ref{thm:twoway-slope} concerns small heights at a fixed speed.
Optimizing the same eight-leg patrol over both of its parameters gives a
lower bound for the unrestricted two-way width when $q$ is close to $1/2$.

\begin{proposition}[Two-way gain below one half]\label{prop:twoway-half}
For $19/40\le q<1/2$ put $\mu=\sqrt{(3q-1)(1-q)}$ and $\nu=q+\mu$. Then
\[
 W_{\mathrm{2way}}(q)\ \ge\ 2r\sqrt{1+\nu^{-2}}\ >\ W_{\mathrm{bf}}(q).
\]
Consequently
\[
 \liminf_{q\uparrow1/2}
 \frac{W_{\mathrm{2way}}(q)-W_{\mathrm{bf}}(q)}{r(1/2-q)^2}\ \ge\ 4\sqrt2 .
\]
\end{proposition}

\begin{proof}
Normalize $r=V=1$. Use the eight-leg patrol of the preceding proof with
\[
 a=\frac{1}{\nu\sqrt{1+\nu^2}},\qquad H=\frac{(1-2q)a}{\mu}.
\]
Since $\mu^2+(1-2q)^2=q^2$, its diagonal has $\ell=\sqrt{a^2+H^2}=qa/\mu$. Hence
\[
 R=q(\ell+a)-(1-2q)H=\Bigl(q+\frac{q^2-(1-2q)^2}{\mu}\Bigr)a=\nu a
 =\frac{1}{\sqrt{1+\nu^2}},\qquad
 \delta=\frac{\nu}{\sqrt{1+\nu^2}},
\]
and $b=a+\delta=\sqrt{1+\nu^{-2}}$. These are the values obtained by first
minimizing $R$ over $H$ at fixed $a$ and then maximizing $a+\sqrt{1-R^2}$
over $a$.

On the stated range $\mu^2\ge(17/40)(21/40)>(0.47)^2$ and $\mu<q<1/2$, so
$0.94<\nu<1$. Therefore $0.70<R<0.73$, $a=R/\nu<0.78$, $\delta=\nu R>0.65$ and,
by $1-2q\le1/20$, $H<0.083$. The three inequalities of the preceding proof
follow:
\[
 0<R<1,\qquad qH<0.05<\min(R,\delta),\qquad
 B_d+H=R+(2+q)H<0.73+2.5\cdot0.083<1 .
\]
Its separator shows that the patrol guards half-width $b$ in both
directions. Since $\nu<2q$ and $x\mapsto2\sqrt{1+x^{-2}}$ is strictly
decreasing, $2b>2\sqrt{1+(2q)^{-2}}=W_{\mathrm{bf}}(q)$.

For the limit put $\varepsilon=1/2-q$. Then
$\mu=\tfrac12\sqrt{1-4\varepsilon-12\varepsilon^2}
 =\tfrac12-\varepsilon-4\varepsilon^2+O(\varepsilon^3)$ and
$\nu=2q-4\varepsilon^2+O(\varepsilon^3)$. The derivative of
$2\sqrt{1+x^{-2}}$ at $x=1$ is $-\sqrt2$. Thus
$2b-W_{\mathrm{bf}}(q)=4\sqrt2\,\varepsilon^2+O(\varepsilon^3)$.
\end{proof}

For example, at $q=0.475$ the patrol has $H\approx0.0811r$ and guards width
about $2.90807r$ in both directions, against $W_{\mathrm{bf}}(0.475)\approx2.90381r$.
Theorem~\ref{thm:twoway} bounds the two-way width only for patrols of a
given height, so a matching upper estimate of this gain would have to
exclude patrols with large excursions.

\subsection{The value of knowing the direction at slow speeds}
Recall the shorthands $W_1=\Wopt$, $W_2=W_{\mathrm{2way}}$, and
$D_i(q)=2/q-W_i(q)/r$ ($i=1,2$). Each competing patrol in $W_1(q)$ or $W_2(q)$
has bounded longitudinal excursion, with no uniform excursion cap imposed across
patrols or as $q\downarrow 0$.

\begin{theorem}[One-way and two-way slow widths]\label{thm:direction-value}
As $q\downarrow0$,
\[
 \frac{qW_2(q)}r\longrightarrow2,\qquad
 \frac{W_1(q)}{W_2(q)}\longrightarrow1,\qquad
 \frac{W_1(q)-W_2(q)}r\longrightarrow+\infty .
\]
The extra width available when the direction is known satisfies
\[
 0<\liminf_{q\downarrow0}\sqrt q\,\frac{W_1(q)-W_2(q)}r
 \ \le\ \limsup_{q\downarrow0}\sqrt q\,\frac{W_1(q)-W_2(q)}r
 \ \le\ 2\sqrt6 .
\]
Consequently $D_2(q)=\Theta(q^{-1/2})$ and the relative gain
$W_1(q)/W_2(q)-1$ has exact order $\sqrt q$.
\end{theorem}

\begin{figure}[tb]
\centering
\includegraphics[width=\linewidth]{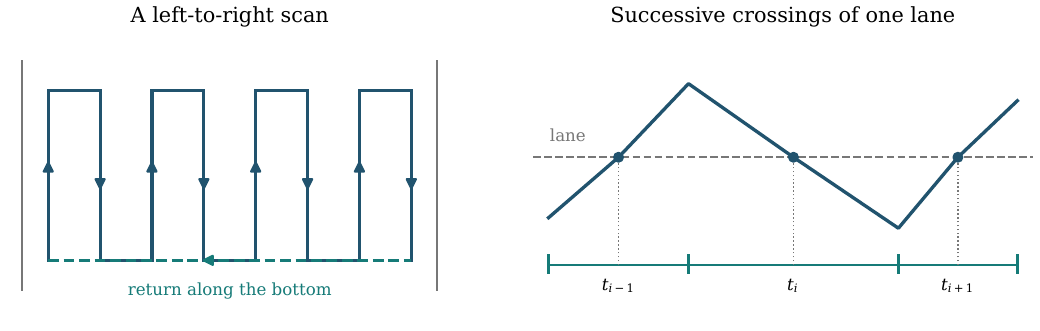}
\caption{Left: one pass of an eight-column two-way scan, sweeping left to right
with a fast bottom return. Right: the timing symmetry forced in a stationary
limit by a hypothetical bounded two-way deficit. Each crossing $t_i$ of a fixed
lane sits at the midpoint of its own erased-height interval of length $a_i$, and
adjacent intervals tile the time axis because $2(t_{i+1}-t_i)=a_i+a_{i+1}$.}\label{fig:scan-crossings}
\end{figure}

\begin{proof}[Construction and reduction to the appendices]
We first construct an elementary vertical-stroke column scan that proves
$qW_2(q)/r\to 2$, gives the explicit upper bound $\limsup_{q\downarrow 0}\sqrt{q}\,D_2(q)\le 8$,
and supplies the finite-height estimate used in Corollary~\ref{cor:two-way-height};
Appendices~\ref{app:direction-loss} and~\ref{app:direction-rate} then prove the
matching $\Omega(q^{-1/2})$ lower bound and sharpen the upper constant from $8$
to $2\sqrt{6}$ via a reclining scan with a diagonal return.

In normalized units $r=V=1$, place an even number $n$ of columns at
$x_j=(j+1/2)s$, $0\le j<n$, inside $[0,ns]$. Scan them from left to right,
alternating vertical strokes between $-h$ and $h$. After each stroke except the
last, move horizontally for time $s$ to the next column; after the final
downstroke, wait at the bottom right and then return horizontally along $y=-h$
at unit speed. Each column block lasts $\tau=2h+s$, and the period is
$P=n(2h+2s)-s$. Figure~\ref{fig:scan-crossings} (left) illustrates one pass.

During pass $k$, the auxiliary transverse sweep $z(t)=v(t-kP)$ with $v=s/\tau$
crosses $[0,ns]$ from left to right, so by continuity it meets the intruder's
transverse coordinate at some time $t_k$ (uniquely when $v>q$). If the intruder's
height at $t_k$ lies in $[-h,h]$, the column block containing $t_k$ reaches that
exact height on its vertical stroke at some time $w$ (assigning the right
endpoint to the final block). At time $w$, the intruder lies within Euclidean
distance $q|t_k-w|+s/2\le q\tau+s/2$ of the guard, ensuring capture whenever
$q\tau+s/2\le1$.

Between consecutive meetings write $d=t_{k+1}-t_k$ and $\rho=q/v$. Because the
horizontal displacement between $t_k$ and $t_{k+1}$ is $v(d-P)$, the speed bound
gives
\[
 |\Delta y|^2\le q^2d^2-v^2(d-P)^2
                  \le q^2P^2/(1-\rho^2),
\]
maximizing the quadratic in $d$. Whenever $qP\le2h\sqrt{1-\rho^2}$, the meeting
heights of an uncaptured trajectory cannot jump from $(-\infty,-h)$ to $(h,\infty)$
or vice versa, blocking complete transits in both directions.

Setting $h=q^{-1/2}$ and
\[
 s=\frac{2(1-2qh)}{1+2q},\qquad \rho=q(2h+s)/s,\qquad
 n=2\left\lfloor\frac{2h(1-\rho^2)+qs}{2q(2h+2s)}\right\rfloor ,
\]
we have $n\ge2$, $\rho<1$, $q\tau+s/2=1$, and $qP\le2h(1-\rho^2)<2h\sqrt{1-\rho^2}$
for all sufficiently small $q$. With $e=\sqrt q$, $s=2-4e+O(e^2)$,
$\rho=e+O(e^2)$, and $qn=1-2e+O(e^2)$ after even-integer rounding, giving
$qW_2(q)/r\ge 2-8\sqrt q+O(q)$. At $q=10^{-4}$, for example, this scan uses
$9{,}806$ columns to guard width $96{,}098{,}800\,r/5{,}001>19{,}215\,r$. More
generally, taking $h=c/\sqrt q$ ($c>0$) yields $s=2-4c\sqrt q+O_c(q)$,
$qn=1-2\sqrt q/c+O_c(q)$, and normalized width
\[
 qns=2-4(c+c^{-1})\sqrt q+O_c(q),
\]
separating the stroke-spacing loss $4c\sqrt q$ from the inter-column and return
loss $4c^{-1}\sqrt q$, minimized at $c=1$.

Combining $W_2\le W_1\le 2r\sqrt{1+q^{-2}}$ from~\eqref{eq:area-bound} with
$qW_2/r\ge 2-8\sqrt q+O(q)$ gives $qW_2(q)/r\to 2$, $W_1(q)/W_2(q)\to 1$, and
$\limsup_{q\downarrow 0}\sqrt q\,D_2(q)\le 8$. Because $D_1(q)=O(1)$ by
Theorem~\ref{thm:slow-hooks} (with $D_1(q)>0$ for small $q$ by
Theorem~\ref{thm:loss}), $(W_1(q)-W_2(q))/r = D_2(q)-D_1(q)$ differs from
$D_2(q)$ by $O(1)$. Appendix~\ref{app:direction-loss} proves $D_2(q)\to+\infty$
by comparing the upward and downward reachable sets erased across a typical
lane, and Appendix~\ref{app:direction-rate} makes that occupation imbalance
quantitative on a single physical time window to prove
$\liminf_{q\downarrow 0}\sqrt q\,D_2(q)>0$, while constructing a reclining scan
with a diagonal return (Lemma~\ref{lem:diagonal-return}) that lowers
$\limsup_{q\downarrow 0}\sqrt q\,(W_1(q)-W_2(q))/r$ to $2\sqrt{6}$. Dividing
$W_1-W_2=\Theta(r/\sqrt q)$ by $W_2\sim 2r/q$ yields $W_1/W_2-1=\Theta(\sqrt q)$.
\end{proof}

Geometrically, the $\Omega(r/\sqrt q)$ lower bound reflects an incompatibility
between symmetric height erasure and asymmetric lane-occupation times. When the
guard crosses a typical lane, bounded or too-small two-way deficit $D_2(q)$ forces
the height intervals erased from intruders arriving from below and from above to
nearly coincide. In a rescaled limit, each crossing of the lane must then sit at
the midpoint of its own erased-height time interval (Figure~\ref{fig:scan-crossings},
right), predicting that the guard spends half its time on each side of the lane.
For a lane one quarter of the way across the corridor, however, the flux budget
requires the guard to spend three quarters of its time on the wider side.

\begin{corollary}[Excursion needed for the full leading width]\label{cor:two-way-height}
For any function $H(q)\ge0$,
\[
 \lim_{q\downarrow0}\frac{qW_{\mathrm{2way}}(q,H(q))}{r}=2
 \quad\Longleftrightarrow\quad
 \lim_{q\downarrow0}H(q)/r=+\infty .
\]
For every fixed finite $h\ge0$, $\limsup_{q\downarrow 0}qW_{\mathrm{2way}}(q,hr)/r<2$,
in contrast to the one-way limit $\lim_{q\downarrow 0}q\Wopt(q,hr)/r=\min\{2,1+h\}$
of Theorem~\ref{thm:slow-limit}.
\end{corollary}
\begin{proof}
When $H(q)/r\to+\infty$, choosing $h=\min\bigl(H(q)/r,q^{-1/2}\bigr)$ in the
vertical-stroke scan above gives $h\to\infty$, $qh\to0$, $s\to2$, $\rho\to0$,
and $qn\to1$ (including even-integer rounding), matching the area
ceiling~\eqref{eq:area-bound}. Conversely, Appendix~\ref{app:fixed-two-direction}
rules out any sequence attaining leading coefficient $2$ under a fixed finite
cap $hr$, and any $H(q)/r$ that fails to diverge admits a bounded subsequence.
\end{proof}

\section{What remains open}
Near the diffraction threshold $q_a$, Theorem~\ref{thm:global-quadratic} establishes quadratic scaling of $F(q)-\Wopt(q)$ without settling whether the normalized deficit
\[
 \lim_{q\uparrow q_a}\frac{F(q)-\Wopt(q)}{r(q_a-q)^2}
\]
converges, or what patrol profile governs its infimum. The known upper bounds descend from $C_4=51.2835202\ldots$ (four straight legs) to $C_5=50.8421136\ldots$ (the single-bend five-leg family of Theorem~\ref{thm:five}) and $C_{12}=48.4250865\ldots$ (Theorem~\ref{thm:bent-deficit}), whereas the positive lower bound $k_1$ from Lemma~\ref{lem:gap-crossing} is existential. Even within five-leg patrols, minimizing the quadratic loss coefficient over vertex placement and stroke asymmetry remains open. Away from $q_a$, it is likewise unknown whether any single patrol attains $\Wopt(q)$ at fixed $0<q<q_a$, or whether an asymmetric four-leg patrol can ever beat $W_4(q)$ outside the left neighborhood of Theorem~\ref{thm:four}.

For constrained longitudinal excursion with $q<1/2$, the exact value of $\Wopt(q,H)$ above the two-level equality cap $(1+q)H=r$ of Theorem~\ref{thm:fixed-height} remains undetermined. In the midrange at $q=1/4$, certified bounds bracket $\Wopt(1/4)$ between $6.94469133\,r$ and $8.10\,r$, leaving open whether any bounded-excursion patrol can guard width $8r$. By Proposition~\ref{prop:q13-floor}, no choice of multipliers in the stateless Galilean flux-and-parity inequality (Theorem~\ref{thm:free-gap}) can push the upper bound below $8r$; crossing that barrier will require coupling consecutive wall-visit times to the interior reachable set. At $q_a$ itself, Theorem~\ref{thm:face-rate} pins down the growth order $N(\varepsilon)=\Theta(\varepsilon^{-1/4})$ for triangulated space-time certificates, leaving open the sharp leading constant $\lim_{\varepsilon\downarrow 0}\varepsilon^{1/4}N(\varepsilon)$: what angular distribution of chords minimizes the deficit per face, and can auxiliary patrol vertices beat the four-leg chord family?

As $q\downarrow 0$ in the one-way game, $D_1(q)=2/q-\Wopt(q)/r$ stays in $[2^{-21},\,1.235352896057745]$ for all sufficiently small $q$ (Theorems~\ref{thm:slow-hooks} and~\ref{thm:loss}). Whether $\lim_{q\downarrow 0}D_1(q)$ exists, and what limiting turn geometry near a wall determines its value, remain open.

In the bidirectional game (Section~\ref{sec:diagnostics}), $W_{\mathrm{2way}}(q)$ is unknown at finite speeds $0<q<1/2$, including the sharp quadratic coefficient of $W_{\mathrm{2way}}(q)-W_{\mathrm{bf}}(q)$ above $4\sqrt{2}\,r(1/2-q)^2$ as $q\uparrow 1/2$ (Proposition~\ref{prop:twoway-half}). As $q\downarrow 0$, Theorem~\ref{thm:direction-value} determines the exact exponent of the directional advantage, leaving open whether $\sqrt{q}\,(W_1(q)-W_2(q))/r$ converges and what leading constant below $2\sqrt{6}$ is optimal. Finally, under a fixed excursion cap $H=hr$ ($h>0$), Corollary~\ref{cor:two-way-height} shows that two-way patrols forfeit a positive fraction of the $2r/q$ leading width; the horizontal sweep and the vertical-stroke scan at fixed stroke half-height $h$ (where $s\to 2$ and $qn\to h/(h+2)$) give
\[
 \max\Bigl\{1,\,\frac{2h}{h+2}\Bigr\}
 \le\liminf_{q\downarrow0}\frac{qW_{\mathrm{2way}}(q,hr)}{r}
 \le\limsup_{q\downarrow0}\frac{qW_{\mathrm{2way}}(q,hr)}{r}<2 ,
\]
and the exact asymptotic coefficient as a function of $h$---even at $h=1$, where $q\Wopt(q,r)/r\to 2$---remains open.

\section*{Verification and AI assistance}
\addcontentsline{toc}{section}{Verification and AI assistance}
The ancillary directory \texttt{anc/} contains exact rational inputs and thirteen standalone Python verifiers, indexed with their hypotheses in \path{anc/README.md}. These scripts verify the critical-root enclosure and obstacle-wavefront margins, the three-regime canonical escape, the scaled and tilted bowties near $q_a$, the seven original triangulated surface-cut lower bounds (64-leg and 18-leg), the eight unrestricted and bounded-height flux-and-parity upper bounds, the fixed-height table at $q=1/4$, the chord-error and 66-face rational certificate at $q_a$ (Section~\ref{sec:face-rate}), the ordinary-bowtie transition covers and 150-face twelve-leg certificate (Theorems~\ref{thm:ordinary-transition} and~\ref{thm:bent-deficit}), and the ten-leg slow-speed guard (Theorem~\ref{thm:slow-hooks}). For the slow-speed and five-leg results they check the critical-root contact margins for Theorem~\ref{thm:four}, the symbolic power series for Theorem~\ref{thm:slow-ordinary}, the rational constants for Theorem~\ref{thm:loss}, both wall handoffs, $C_5$, and the $q=0.45789$ comparison in Appendix~\ref{app:five-leg}, and the finite-lag stationary cancellation identity and $q=10^{-4}$ scan examples in Appendix~\ref{app:direction-rate}.

The accompanying Lean~4 and Mathlib v4.15.0 archive (\path{anc/hallway_guarding_lean4.tar.gz}) contains the wall-budget algebra, selected phase-cylinder and fan identities, critical-root bracket, midline and bounded-height chord inequalities, and kernel replays of the seven original surface-cut and eight continuum upper certificates. The topological uncrossing, periodic polygonal approximation, synchronized-crossing escape transfer, two-sided quadratic estimate, fixed-height attainment, $\Theta(\varepsilon^{-1/4})$ face-count bounds, four- and five-leg transition proofs, and BV reachable-set arguments for the slow-speed one-way and directional bounds are proved analytically in the text, supported by the finite checks above.

OpenAI GPT/Codex, Anthropic Claude, and Google DeepMind Gemini/Antigravity assisted with mathematical exploration, programming, proof auditing, drafting, and Lean tactics. Throughout the paper, analytic proofs and executable finite checks play distinct, explicitly delineated roles.

\appendix
\section{Comparison with every ordinary bowtie}\label{app:ordinary}

This appendix proves the strict inequality $W_4(q)<\Wopt(q)$ for every
$0<q<q_a$ in Theorem~\ref{thm:ordinary-transition}. Subsection~\ref{app:lane-tests}
derives three straight-runner upper bounds on an arbitrary ordinary bowtie;
Subsection~\ref{app:slow-comparison} beats those bounds across $0<q\le 2/5$ with
an explicit eight-leg patrol, and Subsection~\ref{app:middle-comparison}
handles $2/5\le q<q_a$ by localizing any class-maximizing bowtie and bending its
vertical strokes inward.

Normalize $r=V=1$ throughout the appendices. For an ordinary symmetric bowtie
with column half-spacing $a\ge 0$ and stroke half-height $h\ge 0$ (and $\ell=\sqrt{a^2+h^2}>0$),
put $d=h/\ell$, $c=\sqrt{1-d^2}$, and let $\sigma=+1$ for descending vertical
strokes or $\sigma=-1$ for ascending strokes. In a corridor $[-b,b]\times\R$ of
half-width $b$, if the bowtie's axis of symmetry is offset transversely by $\xi$,
its distance to the farther wall is $b_f=b+|\xi|$ and the farther wall inset is
$\delta_f=b_f-a\ge0$. Blocking straight wall runners requires $\delta_f\le1$,
and a degenerate $a=0$ column guards only $b\le1$, whereas every benchmark below
exceeds one.

\subsection{Three tests with straight runners}\label{app:lane-tests}

In the coordinates $(X,qt-Y)$ of a straight runner, one ordinary
period advances by $P=4q(h+\ell)$. Put
\[
 S_\sigma=\sqrt{1+q^2-2\sigma qd},\qquad
 E_+=2q-(1-q)d,\qquad E_-=2q-(1-3q)d .
\]
The tests below apply for $0<q<1/2$ and $b>1$. The other column is
farther than one from the tested wall.

\paragraph{Farther wall.}
Center the nearer stroke at phase zero and test phase $P/2$.
In the forward orientation the two endpoints have phases
$\pm s_0$, $s_0=(1+q)h$. Write $G=P/2-s_0=\ell E_+$.
Whenever $G>0$ the departing and arriving diagonals and the nearer
stroke give all candidates for the least distance. Indeed the departing
diagonal ranges in phase between $s_0$ and $G$, both in $[0,P/2]$,
and reflection gives the arriving one. If $d\ge q$, both the horizontal
and phase distances along the departure are at least the endpoint values.
Thus $\delta_f^2+G^2\le1$. For $d<q$ put $e_0=q-d$.
The squared distance along the departure at elapsed time $u$ is
$(\delta_f+cu)^2+(G-e_0u)^2$. Its minimum on the ray is at the
endpoint if $e_0G\le c\delta_f$. In either endpoint case,
Cauchy--Schwarz gives $b_f\le\sqrt{1+(c/E_+)^2}$. Otherwise
capture on the ray requires $aE_++e_0\delta_f\le S_+$ by distance
to its supporting line. Together with $e_0G>c\delta_f$ this gives
$\delta_f<e_0/S_+$. Since $E_+>e_0$,
\[
 b_f\le \frac{c^2+e_0 E_+}{S_+ E_+}\le\sqrt{1+(c/E_+)^2}.
\]
For $0\le d\le q<1/2$, $c/E_+$ increases with $d$.

For ascending strokes the endpoints have phases of magnitude
$s_0=(1-q)h$, traversed in reverse order. If $G_c=P/2-s_0=\ell E_->0$,
the stroke can capture only if $\delta_f^2+G_c^2\le1$. The departing
diagonal starts at $-s_0$ and ends at $P/2+s_0$. Other copies of
the test phase lie at distance at least $G_c$ throughout the leg.
If the stroke misses, the diagonal or its reflection must reach the
test disk along the ray with $e_0=q+d$ and initial phase distance
$G_a=P/2+s_0=\ell E_a$, $E_a=2q+(1+q)d>e_0$. Its starting endpoint
also misses. The supporting-line calculation above gives
$b_f\le\sqrt{1+(c/E_a)^2}\le\sqrt{1+(2q)^{-2}}$. For $q\ge1/3$,
$E_-\ge2q$, so the stroke branch satisfies this same bound.

\paragraph{Center lane.}
At the bowtie axis the two diagonals have phases spaced $P/2$ apart.
Halfway between them the distances to all their supporting lines are
at least $(P/4)c/S_\sigma=aq(1+d)/S_\sigma$. Both strokes are
horizontally $a$ away. Hence
\[
 a\le\max\{1,S_\sigma/[q(1+d)]\}.
\]

\paragraph{An interior lane.}
Suppose $n=d-\sigma q\ge0$, $s_0=(1+\sigma q)h\le P/2$ and
$(1+\sigma q)d\le S_\sigma+n$. Set $L_0=s_0c/(S_\sigma+n)\le a$.
On the lane $x=a-L_0$ relative to the axis, the adjacent diagonal
lines have phases $\pm(s_0-nL_0/c)$ at the stroke midpoint's phase
origin. Here $0\le s_0-nL_0/c=L_0 S_\sigma/c\le P/2$. All periodic
copies of these lines have distance at least $L_0$ from the test
point, and both columns have horizontal distance at least $L_0$.
Guarding requires $L_0\le1$, or
\[
 a(1+\sigma q)d\le S_\sigma+d-\sigma q.
\]
The center and interior lanes lie between the patrol columns for
every permitted offset.

\subsection{A direct comparison for \texorpdfstring{$0<q\le2/5$}{0 < q at most 2/5}}\label{app:slow-comparison}

Keep the two vertical strokes and add two bends to each diagonal.
Use the transfer in the proof of Theorem~\ref{thm:slow-hooks} with
$e=0$, outer horizontal length one and
\[
 h=\frac{6}{5(1+q)},\qquad
 \tau=\frac{q/5+\sqrt{26/25-q^2}}{1-q^2},\qquad
 v_0=\tfrac15+q\tau .
\]
Write $Y=h-v_0$, $\alpha=\sqrt{1-q^2}$ and $x=\alpha Y/q$.
Take $a=x+1$, $G=4/5$, $\delta=3/5$. After the stroke from $(a,h)$
to $(a,-h)$ the three segments run through $(x,-Y),(-x,Y),(-a,h)$.
Reflect horizontally for the rest of the period. The outer segments
take time $\tau=\sqrt{1+v_0^2}$, and the middle segment has upward
velocity $q$. At the two corners the rising flat levels of the
separator have offsets $+1$ and $-1$, since $Y+q(h+\tau)=1$.
At the stroke midpoint the anchor $(x,0)$ and the wall points
$(a+\delta,\pm G)$ lie in the current disk. The anchor has its
future witness at the first corner. The lower wall point has its
future witness at departure since $qh<1$ and
\[
 \delta^2+(h-G)^2\le(1-qh)^2
 \quad\Longleftrightarrow\quad h(1-q^2)\le2(G-q).
\]
The last inequality holds for $q\le1/2$. At departure the wall
tangency is $\delta^2+G^2=1$. With $Y>0$ these are precisely the
finite-leg witnesses used in the interpolation proof of
Theorem~\ref{thm:slow-hooks}, including its captured reset triangles.
They certify the half-width
\[
 b_L=\frac1q-L(q),\qquad
 L(q)=\frac{q}{1+\alpha}+\alpha(h+\tau)-\frac85 .
\]

We show that every ordinary half-width is at most $b_L-10^{-4}$
throughout $0<q\le2/5$. The following details specify the finite
interval proof in \path{anc/bends/check_bowtie_transition.py}. Put
$l=L(q)+10^{-4}$. The checker verifies $Y>0$, $l<1$,
$1-ql>\sqrt{q^2+1/4}$, and
\[
 l<\frac{2/(1+q)-q}{1+\sqrt{q^2+(1-q)/(1+q)}} .
\]
The third inequality excludes the reverse diagonal branch whenever
$G_c>0$. The last excludes $d\le q$ forward by the farther-wall test
at $d=q$, rationalizing $1/q-\sqrt{1+(c(q)/E_+(q))^2}$, with
$E_+(q)=q(1+q)$.

Suppose $b\ge1/q-l$. Set $v=d/q$ and $t=1/q-a$. The center test
and $\delta_f\le1$ give $-1\le t\le l+1$. Also $0\le v\le5$.
Indeed for $q\le1/5$ we have $qa\ge1-2q>q$, whence
$(1-2q)(1+d)\le S_\sigma\le1+q$ and $v\le3/(1-2q)\le5$.
For $q\ge1/5$, $d\le1$ suffices. It remains to exclude the rectangle
$[0,2/5]\times[0,5]$ in $(q,v)$ for each orientation. The tests use
only expressions regular at zero. For $1-qt>q$ the center test becomes
\[
 t\ge \frac{v-q(1-2\sigma v)/(1+S_\sigma)}{1+qv}.
\]
Put $f_\sigma=2-[1-(2-\sigma)q]v$. If $f_\sigma>0$ the wall
condition requires $G=(1-qt)f_\sigma/c\le1$ and
$t\le l+\sqrt{1-G^2}$. Under the interior-lane hypotheses the further
condition is $(1-qt)v(1+\sigma q)\le S_\sigma+q(v-\sigma)$.
Bisection in $q,v$ and outward interval evaluation of these
inequalities exclude every rectangle. On a rectangle, each new lower
or upper bound on $t$ follows from one of the displayed inequalities.
The forward cover has 1,479 terminal rectangles and the reverse
cover has 8,992. Replay prints the counts for each rejection reason.
Zero enters only as an endpoint of these regular expressions.
The uniform margin $10^{-4}$ gives a strict comparison of the suprema.

\subsection{The active height constraint for \texorpdfstring{$2/5\le q<q_a$}{2/5 through the critical speed}}\label{app:middle-comparison}

Write $c(z)=\sqrt{1-z^2}$. For inclination $z>q$ take
$m=qz+\alpha c(z)$, $\kappa=\alpha z-qc(z)$,
$M=\sqrt{1+2qm+q^2}$, $C=(1+qm)/M$, $S=q\kappa/M$.
Thus $C=\cos\theta_*$, $S=\sin\theta_*$ for the fan in
\eqref{eq:theta-p}. With $p$ from that equation put
\[
 K=\kappa C-mS+(m+q)p,\qquad H(z,q)=\frac{1+K}{m+q}
     =p+\frac{1+\kappa/M}{m+q}.
\]
For an ordinary bowtie set $z=d$ and $A_4=cH/d$.
We record conditions under which guarding requires $a\le A_4$.
On the geometry $a=A_4$, $h=H$, put $T_0=2H(1+1/d)$. Suppose
\[
 0<p<H,\quad A_4>1,\quad
 qmT_0+S-(1+qm)p>H+1,\quad C-q\kappa(T_0-p)>0 .
\]
Use the contact arc $(a-\cos\theta,\sin\theta-t)$ up to $t=p$,
where $t$ is the fan integral through $\theta$. Continue with velocity
$q(\kappa,m)$ until $T_0$, the midpoint of the next left stroke,
then vertically at speed $qm$. The arc has clearance one on the
right stroke and its tangent has clearance at least one on the
remaining stroke. On the diagonal the signed separation along
$(\kappa,m)$ is $H(m+q)-K=1$. On the left stroke the horizontal
distance is at least $2A_4-1>1$, and the vertical tail starts more
than one above every subsequent patrol position. The last displayed
condition keeps the trajectory between the columns. Reflection in
time and height supplies the incoming half. Scaling all positions
and times by $a/A_4>1$ gives a strict escape at the same angle and
speeds, for every corridor containing the columns.

For sufficiency, allow the two strokes to bend inward. Let $\eta\ge0$,
$e=h\eta$, $\beta=1/2$, $j(\eta)=\beta+\sqrt{(1-\beta)^2+\eta^2}$,
$\rho=h j(\eta)$. The right stroke runs through
$(a,h),(a-e,\beta h),(a-e,-\beta h),(a,-h)$. Keep the straight
crossings of duration $2\ell$, and take $G=q(2\ell+\rho)-h$,
$\delta=\sqrt{1-G^2}$. The handoff to the next crossing requires
\[
 mh+q\rho-\kappa e\le1+K.
\]
The construction and its finite witnesses are proved in
Appendix~\ref{app:bent-separator}. At $\eta=0$ its strict auxiliary
conditions reduce to
\[
 0<G<1,\quad 2a+\delta>1,\quad 0<p<h/2,\quad
 q(h-p)<1,\quad 2(G-q)>(1-q^2)(h-p),
\]
together with $d x_j+c y_j>0$ and $0<(-m x_j+\kappa y_j)/\alpha<2\ell$
at $(x_0,y_0)=(-1,h+qp)$ and $(x_1,y_1)=(-C,h-p+S)$.
The wall inequality here is the expanded future-witness distance
at the fan handoff. In an open neighborhood of any geometry with
these strict margins, every centered ordinary patrol with $h\le H$
attains its wall-tight half-width $B(a,d;q)$ of Section~\ref{subsec:tilted}.
The margins also persist for sufficiently small inward bends.

\paragraph{Localizing the best ordinary patrol.}
We give the full domain of the second interval calculation. For each
$2/5\le q\le0.457893$, let $d_T=3q/2-19/200$ and
$b_T=B(A_4(d_T,q),d_T;q)$. The checker verifies the strict
auxiliary conditions at this trial and hence its attainment. It also
checks $b_T>\sqrt{1+(2q)^{-2}}$ and $b_T\sqrt D>1$,
where $D=3q^2+2q-1$. The reverse traversal is thereby excluded.
Every forward competitor at or above $b_T$ must have $d>q$, as
verified using the farther-wall bound on $0\le d\le1$. For $d>q$
write $k=E_+/c$. Since $k^2\ge D$ and $a\le1/k<b_T$, the ellipse
$b_T\le a+\sqrt{1-a^2k^2}$ confines $a$ between
\[
 \frac{b_T-\sqrt{1-k^2(b_T^2-1)}}{1+k^2}
 \quad\text{and}\quad
 \frac{b_T+\sqrt{1-k^2(b_T^2-1)}}{1+k^2}.
\]
An outward lower endpoint for $b_T$ gives the same necessary test.
The center and interior tests further decrease the amplitude upper
bound. For each surviving rectangle in $(q,d)$ the checker verifies
the four strict escape conditions above and imposes $a\le A_4$. It
either excludes the remaining interval of amplitudes or verifies
all strict auxiliary capture conditions throughout it. For the arc
time, the substitution $u=\sin\theta$ gives an increasing integrand
on $[0,S]$, so that
\[
 \frac{S}{1+q}\le p\le
 \frac{S}{\sqrt{1-S^2}\,(\sqrt{1-S^2}+\sqrt{q^2-S^2})}.
\]
These elementary bounds suffice for the localization. Bisection
starts at $[2/5,0.457893]\times[0,1]$. Of its 25,032 terminal
rectangles, 14,907 exclude contenders and 10,125 certify the strict
auxiliary conditions for every remaining amplitude. Both covers use
exact fractions and outward integer square roots.

\paragraph{Strict improvement.}
Fix $2/5\le q<q_a$. If the trial already attains the ordinary
supremum, use it. Otherwise take a maximizing sequence with
half-widths at least $b_T$. By the center test its amplitudes are
bounded. Pass to a convergent subsequence of $(a,d)$. The closed
wall and straight-lane inequalities and the finite cover just
described give $d>q$, $a\le A_4$ and all strict auxiliary margins at
the limit. More explicitly, finitely many closed rectangles suffice,
so one may pass to a subsequence in a single surviving rectangle,
where the bounds hold on the closure. The centered patrol therefore
attains $B(a,d;q)$, at least the limiting half-width. This is a class
maximum, including over offsets.

In a neighborhood of this maximum the smooth constraint $a\le A_4(d,q)$
suffices for attainment of $B$. Its only unconstrained stationary
point is the wall optimum $(a_*,d_*)$, where
$d_*=(1-q)/(2q)$ and $a_*=1/\sqrt{D(1+D)}$.
Indeed $B_a=1-ak^2/\sqrt{1-a^2k^2}$ and
$k'=(2qd-(1-q))/c^3$. Attainment at that point would give width
$F(q)$, excluded by Theorem~\ref{thm:global-gap}. Hence $a=A_4$,
and inward variation of $a$ shows $B_a\ge0$. If $B_a=0$,
differentiating along $a=A_4(d,q)$ gives $B_d=0$, which again
would be the wall optimum. Thus $B_a>0$.

Hold $d$ fixed and impose the bent handoff at equality by setting
\[
 h(\eta)=\frac{1+K}{m+qj(\eta)-\kappa\eta},\qquad
 a(\eta)=c h(\eta)/d,\qquad b(\eta)=a(\eta)+\sqrt{1-G(\eta)^2}.
\]
Since $j'(0)=0$, we have $a'(0)=a\kappa/(m+q)>0$ and
$G(\eta)=a(\eta)k+qh(\eta)(j(\eta)-1)$. Consequently
$b'(0)=B_a a'(0)>0$. The auxiliary inequalities persist, so
sufficiently small positive $\eta$ gives an eight-leg guard strictly
wider than the entire ordinary class. This completes the strict
part of Theorem~\ref{thm:ordinary-transition}.

\section{Bent patrols and the quadratic coefficient}\label{app:variation}

This appendix proves Theorem~\ref{thm:bent-deficit} and establishes the
finite-leg bent separator (Lemma~\ref{lem:bent-separator}) whose straight-crossing
specialization was also used in Appendix~\ref{app:middle-comparison}. Starting
from an ordinary bowtie, each vertical stroke is bent inward into three segments
and each diagonal crossing is broken symmetrically into three segments, giving a
twelve-leg periodic patrol whose second variation at $q_a$ determines both the
exact ordinary limit $C_4$ and the twelve-leg upper coefficient $C_{12}$.

\subsection{A separator on the actual finite legs}\label{app:bent-separator}

Write $c(z):=\sqrt{1-z^2}$ for $z\in(0,1)$. Prescribe a half-crossing duration
$L=x/c(d)$, a mean upward velocity
$d$ and a fraction $\lambda$ of each half-crossing nearest the pole.
The outer pieces have upward velocity $z=d+s$ and the middle piece
has velocity $z_i=d-\lambda s/(1-\lambda)$. Set
\[
 \begin{aligned}
 A&=L[\lambda c(z)+(1-\lambda)c(z_i)], & h&=Ld,\\
 e&=h\eta, & \rho&=h j(\eta),\qquad
 j(\eta)=\beta+\sqrt{(1-\beta)^2+\eta^2}.
 \end{aligned}
\]
With $0<\lambda,\beta<1$ the right stroke goes through
$(A,h),(A-e,\beta h),(A-e,-\beta h),(A,-h)$. The crossing runs
through $P_0,P_1,-P_1,-P_0$, where
$P_0=(A,-h)$ and $P_1=(A-c(z)\tau,-h+z\tau)$, $\tau=\lambda L$.
Horizontal reflection closes the period. Every leg has speed one.
Put $G=q(2L+\rho)-h$, $\delta=\sqrt{1-G^2}$, $b=A+\delta$. At
$s=\eta=0$ this gives the ordinary half-width $B(x,d;q)$.

\begin{lemma}[Finite-leg handoff]\label{lem:bent-separator}
Use the fan parameters $m,\kappa,C,S,p,K$ of
Appendix~\ref{app:middle-comparison} at the outer inclination $z$.
Suppose $0<q<z<1$, $0<z_i<1$, $h>0$, $0\le e<A$, $0<G<1$,
$\delta+e<1$, $0<p<\beta h$ and $2A+\delta-e>1$. Let
$R_e=\sqrt{1-(\delta+e)^2}$, $v=\rho-p$.
The patrol guards half-width $b$ if
\[
 \begin{aligned}
 mh+q\rho-\kappa e &\le1+K,\\
 qv&\le1,\qquad \delta^2+(h-p-R_e)^2\le(1-qv)^2,\\
 |h+q\rho-(z-q)\tau|&\le1,\qquad |-G-(z-q)\tau|\le1 .
 \end{aligned}
\]
In addition, at both $(\xi_0,y_0)=(-1-e,h+qp)$ and
$(\xi_1,y_1)=(-C-e,h-p+S)$ require
\[
 z\xi_j+c(z)y_j\ge0,\qquad
 0\le(-m\xi_j+\kappa y_j)/\alpha\le\tau .
\]
For a single straight crossing ($s=0$), omit the two corner
conditions and replace $\tau$ by the full duration $2L$ in the
witness bound.
\end{lemma}
\begin{proof}
For an affine piece of the patrol on $I$, the points $(t,P)$ with
a witness $w\in I$ satisfying $\|P-g(w)\|+q|w-t|\le1$ form a
convex set. Interpolating the witnesses proves convexity. Any
admissible intrusion meeting the set is caught at the witness time.

In centered stroke time $|\sigma|\le p$ the guard is $(A-e,-\sigma)$.
Attach the fan of Lemma~\ref{lem:attached-fan} translated by $-e$.
For positive $\sigma$ it joins the flat $y=q\sigma$ to its contact
on the current disk. Reflect time and height for negative $\sigma$.
Join the contact to $(b,-\sigma-R_e)$ by a chord in the disk.
The flat endpoint lies in the strip by $2A+\delta-e>1$.

Set physical time zero at departure from $P_0$. At $\sigma=p$,
$t=-v$, use coordinates $(\xi,y)$ relative to $P_0$. If
$N=z\xi+c(z)y\ge0$ and $w=(-m\xi+\kappa y)/\alpha$ belongs to the
first crossing piece, substitution in that piece gives
\[
 \|P_0+(\xi,y)-g(w)\|+q(w+v)=qv+\kappa\xi+my.
\]
From the flat endpoint to the contact, the slopes of the fan lie
in $[-\kappa/m,0]$. The cost and $N$ are nondecreasing and $w$ is
nonincreasing, using $zm-c(z)\kappa=q>0$. The two endpoint tests
therefore place every witness on the finite first piece. The largest
cost, at the contact, is $mh+q\rho-K-\kappa e\le1$. The wall
endpoint has the departure witness by the second line of inequalities.
The fan and its chord can now be replaced by the chord between their
extreme endpoints, sweeping only this same convex capture set.

Let $F_+(t)=q(\rho+t)$ and $F_-(t)=-h-G+qt$. On $F_+$ take $U(t)$
with its transverse coordinate interpolating affinely through the nodes
\[
 \begin{array}{c|rrrr}
 t&-v&\tau&2L-\tau&2L\\ \hline
 U_x(t)&A-e-1&(P_1)_x&-(P_1)_x&-A .
 \end{array}
\]
The first node is the fan's flat endpoint. At the first corner both
flat levels are in the current disk by the two corner conditions,
and reflection gives the other corner. At arrival $U$ has vertical
offset $G$. Thus each interpolated segment lies in the capture set
of its own leg. During the crossing $-U(2L-t)$ lies on $F_-$ in the
same set as $U(t)$. The upper endpoint lies to the left of the
lower one. Their horizontal difference is affine on each piece,
nonpositive at departure, zero at both corners, and nonpositive at
arrival by reflection. A captured chord joins them, with the upper
flat extended to the left wall and the lower flat to the right.
For a single straight crossing simply interpolate between the first
and last nodes in its single capture set. The same horizontal
ordering follows at the two endpoints of the crossing.

Before departure join $U(t)$ to the wall trace interpolating between
$(b,-p-R_e)$ at $t=-v$ and $(b,-h-G)$ at $t=0$. Both nodes are
captured on the first piece, using $\delta^2+G^2=1$ at departure.
At departure only the triangle between $U(0),(A,-h-G),(b,-h-G)$
is swept by the change of chords, in that same set. Time and
position reflection about the crossing midpoint give the arrival
and the reverse handoff to the next fan. At that handoff reset the
wall offset from $+R_e$ to $-R_e$ across a triangle in the current
disk. Repeating by horizontal reflection gives the bounded separator.
On every uncaptured graph $f_t=q\sqrt{1+f_x^2}$, so along any
speed-$q$ path its height relative to the graph cannot increase,
by Cauchy--Schwarz. All resets and intervening chords are captured.
A complete upward transit must pass from below to above this
separator, proving detection.
\end{proof}

\subsection{The second variation}

For the twelve-leg construction fix $\lambda=93/100$ and
$\beta=1/20$. Define the half-width $\mathcal B=A+\sqrt{1-G^2}$
and the handoff ratio
\[
 \mathcal R=h\,\frac{m+qj(\eta)-\kappa\eta}{1+K},\qquad
 (m,\kappa,K)\ \text{evaluated at }z=d+s.
\]
The wall and handoff are tight at $(a_*,d_*,0,0;q_a)$. All the
auxiliary capture tests have strict slack there, apart from the
one-sided requirement $\eta\ge0$. For example the largest first-piece
witness is below $1.256$, while that piece lasts longer than $1.257$.
Holding $\mathcal R=1$ fixes the scale by
\[
 h=\frac{1+K}{m+qj(\eta)-\kappa\eta},\qquad x=h c(d)/d .
\]
For $d$ close to $d_*$, small $s$ and small $\eta\ge0$,
Lemma~\ref{lem:bent-separator} therefore gives a guard.

At the wall optimum $(a_*,d_*,0,0)$ the width gradient in all four
coordinates vanishes. The negative Hessian at $q_a$ is diagonal
with entries
\[
 M_x=\frac{k^2}{\delta^3},\quad
 M_d=\frac{x^2 k k''}{\delta},\quad
 M_s=\frac{x\lambda}{(1-\lambda)c(d)^4},\quad
 M_\eta=\frac{G}{\delta}\frac{qh}{1-\beta}.
\]
The ordinary entries are from Proposition~\ref{prop:tilted}, using
$k'(d_*)=0$. The weighted mean of $c(d+s)$ and $c(z_i)$ has zero
first derivative in $s$ for every $d$ and second derivative
$-\lambda/[(1-\lambda)c(d)^3]$ at zero. On the stroke $j'(0)=0$,
$j''(0)=1/(1-\beta)$. These identities give the other entries and
the vanishing mixed terms.

Let $g=\nabla\mathcal R$ at the critical patrol and
$E_0=-\frac{d}{dq}\mathcal R(a_*(q),d_*(q),0,0;q)|_{q_a}>0$.
Writing $H_z=\partial H/\partial z$, the entries are
\[
 g_x=1/x,\qquad g_s=-H_z/H,\qquad
 g_d=\frac1{d c(d)^2}+g_s,\qquad g_\eta=-\frac{\kappa}{m+q}.
\]
Indeed $m'=-\kappa/c(z)$, $\kappa'=m/c(z)$, and the terminal-angle
terms in differentiating $K$ cancel because the fan integrand at
$\theta_*$ equals $1/M$ and $mC+\kappa S=(m+q)/M$. In particular
$H_z=[mC+\kappa(S+H-p)]/[c(z)(m+q)]>0$.

For allowed coordinates $I$ the least $v^{\mathsf T}M v$ subject
to $g\cdot v=-E_0$ is
\[
 C_I=\frac{E_0^2}{\sum_{i\in I}g_i^2/M_i},\qquad
 v_i=-\frac{E_0\,g_i/M_i}{\sum_{j\in I}g_j^2/M_j}\quad(i\in I).
\]
Set the other entries to zero. For $I=\{x,d\}$ this is
$C_4=C_{\mathrm{tilt}}$. Adding $s$ gives $49.3948159\ldots$;
adding both $s,\eta$ gives $C_{12}=48.4250865\ldots$.
To realize these values take $d=d_*(q)+v_d\varepsilon$,
$s=v_s\varepsilon$, $\eta=v_\eta\varepsilon$ and use the scale
above. The sign of $g_\eta$ gives $v_\eta\ge0$. Differentiating
$\mathcal R=1$ gives $x=a_*(q)+v_x\varepsilon+O(\varepsilon^2)$.
Taylor expansion at the wall optimum at the same speed yields
$F-2\mathcal B=C_I\varepsilon^2+O(\varepsilon^3)$.

For the converse in the ordinary class, take feasible widths within
$o(\varepsilon^2)$ of $W_4$. The preceding construction puts their
loss at $O(\varepsilon^2)$. The reverse direction is excluded by
Appendix~\ref{app:lane-tests}. For forward competitors write
$b_*=F/2$, $E=E_+$. The exact identity
\[
 b_*^2-\bigl(1+c(d)^2/E^2\bigr)=
 \frac{(2q)^2(d-d_*)^2}{DE^2}
\]
gives $d=d_*+O(\varepsilon)>q$. The wall-tight bound $b\le B(x,d;q)$
then applies with $x=a$. Its maximum at fixed $d$ is at
$x_0(d)=1/[k\sqrt{1+k^2}]$. Since $-B_{xx}\ge k^2\ge D$,
$x=x_0(d)+O(\varepsilon)$, and $k^2=D+(2q)^2(d-d_*)^2/c(d)^2$
gives $x_0(d)=a_*+O(\varepsilon^2)$. At each limit $v$ of the
scaled ordinary displacements, the necessary rounded-escape
inequality $x d/[c(d)H(d,q)]\le1$ from Appendix~\ref{app:middle-comparison}
gives $g_xv_x+g_dv_d\le-E_0$. Taylor expansion of $2B$ forces
a loss coefficient at least $C_4$. This proves the limit for $W_4$
and the claimed unrestricted limsup.

\subsection{Uniform estimates and a finite example}\label{app:bent-uniform}

The verifier \path{anc/bends/check_critical_bends.py} covers
$[0.4578,0.457892036652]$ with twelve equal rational intervals,
using enclosure~\eqref{eq:qa-enclosure} for the constant $q_a$.
It reuses the integral with its explicit derivative remainder and
the wall Hessian from the escape and tilted-bowtie checks.

Suppose an ordinary guard has $F-2b\le50\varepsilon^2$. The
reverse bound is smaller than $F-50\varepsilon^2$. The identity
above first gives $|d-d_*|\le7\varepsilon$ using $E\le2q$, and then
$|d-d_*|\le4.5\varepsilon$ on reinserting the angle interval.
Likewise $-B_{xx}\ge D$ gives $|x-x_0|\le10\varepsilon$ and its
value on this interval sharpens the factor to $4.5$. On the
resulting rectangles the necessary handoff ratio $\mathcal R_4=
x d/[c(d)H(d,q)]$ is at most one for every guard. The checker
also tests the complete rounded escape directly there. It takes
$R=h/H$, scales the contact arc and its tangent by $R$, switches
to the vertical at height $21/10$ during the following left stroke,
and checks $a>R$, $21/10-h>R$ and that the tangent stays between
the columns, as in Theorem~\ref{thm:canonical-escape}.

Let $m_x,m_d$ bound the corresponding diagonal entries of
$-\nabla^2B$ from below after subtracting the largest absolute
mixed entry. Let $G_x,G_d$ bound $|\partial_x\mathcal R_4|$ and
$|\partial_d\mathcal R_4|$ from above, and let $e_0$ be a common
lower bound on the negative derivative of the wall-optimal patrol's
ratio over all speed intervals. With $\Delta=(x-a_*,d-d_*)$,
Taylor's theorem and the mean value theorem give
\[
 F-2b\ge m_x\Delta_x^2+m_d\Delta_d^2,\qquad
 e_0\varepsilon\le G_x|\Delta_x|+G_d|\Delta_d|.
\]
Weighted Cauchy--Schwarz gives a lower coefficient
$e_0^2/\max(G_x^2/m_x+G_d^2/m_d)>50.15$, which suffices for
the bound $50$.

For the guarded width use $d=d_*-3.527\varepsilon$,
$s=0.226\varepsilon$, $\eta=1.348\varepsilon$ and impose the
handoff at equality. The smallest auxiliary margin tested in
Lemma~\ref{lem:bent-separator} exceeds $0.00128$. First bound
$F-2B(x,d;q)$ by the wall Hessian on the rectangle joining the
design to $(a_*,d_*)$, using
$|x-a_*|\le\varepsilon\max|(x-a_*)'|$ over the full speed interval.
The additional full-width cost $2(x-A)$ is at most
$L\lambda s^2/[(1-\lambda)c_{\min}^3]$, with $c_{\min}$ the
minimum between the two inclinations, by cancellation of the linear
terms in Taylor's formula for $c$. Since
$G-xk=qh(j(\eta)-1)\le qh\eta^2/[2(1-\beta)]$, the loss from the
smaller inset is at most $(G/\delta)qh\eta^2/(1-\beta)$.
The checker bounds the sum of the three coefficients by $48.689$.
This proves Theorem~\ref{thm:bent-deficit}.

At $q=0.45785$ the same formula with $0.00004203665$ in place of
$\varepsilon$ in the design gives the stated analytic guard by
direct substitution. The independent finite witness
\path{anc/bends/bent_surface_045785.json} has 99 vertices, 72
causal faces and 78 captured faces. Its two central vertical pieces
are each divided at their midpoint, giving fourteen affine time
intervals for twelve geometric legs. Proposition~\ref{prop:surface-cut}
and the exact rational replay give its certified width. The optional
generator replaces each contact arc by four chords using the chord
times of Section~\ref{sec:face-rate}. Verification uses only the
rounded rational patrol, vertices and witnesses in the saved file.

\section{Symmetrizing four arbitrary straight legs}\label{app:four-legs}
We prove Theorem~\ref{thm:four}, showing that on a left neighborhood of $q_a$
no periodic patrol with at most four monotonically traversed straight legs---allowing
arbitrary vertex positions, variable speed $\le 1$, and waits---can exceed the
ordinary symmetric bowtie supremum $W_4(q)$. The proof first localizes any
near-optimal four-leg itinerary to the canonical critical bowtie with vanishing
waiting delay, then shows that the quadratic wall-phase penalty for asymmetry plus
the linear penalty for delay strictly outweighs any amplitude gained from relaxing
the oblique-stroke rounded-escape bound.

\paragraph{Localization.}
In radius-one units, consider any sequence $q\uparrow q_a$ and
competitors of full width $2b\ge W_4(q)$. The ordinary lower construction
and the wall ceiling give $b_*(q)-b\to0$, where
$b_*=\sqrt{1+(3q^2+2q-1)^{-1}}$. Lemma~\ref{lem:gap-crossing},
applied to the original schedule with
five canonical cycles, gives convergence on those cycles to the
critical bowtie, uniformly in physical time after translations. Its
four unit-speed legs meet at genuine corners and its minimal period
is $T_*$. Each sufficiently small fixed time neighborhood of any
interior canonical corner must contain an actual vertex. Otherwise
the three positions at the ends and the middle of that neighborhood
would be collinear in the approximating paths, and hence in the limit.
The five cycles have twenty such neighborhoods. A period at least
$3T_*/2$ could give at most sixteen vertices in this interval for all
sufficiently large indices. Periods are bounded away from zero by the
positive transverse travel and the speed bound. The equalities at
translated times persist throughout a middle full cycle, so any
limiting period is a period of the critical bowtie. It must be
$T_*$. There are now four vertices, one near each critical corner,
with convergent times and positions. Writing $\ell_j$ for their geometric
leg lengths, the total delay $e=T-\sum_{j=1}^4\ell_j\ge0$ tends to zero.
This proves localization even for arbitrarily chosen
periods and clocks.

\paragraph{The wall cost of asymmetry.}
In the resulting neighborhood label the vertices $R_t,R_b,L_t,L_b$ in
itinerary order. Subtract the mean vertical coordinate, and let $V$ be
the resulting vertex tuple. The sign $s_i$ of pole $i$ is $+1$ on the
right and $-1$ on the left. For its endpoints put
$r_j=\sqrt{1-(b-s_i x_j)^2}$. Write $\ell_i$ and $d_i$ for stroke
length and vertical drop, and $T_g=\sum_j\ell_j$. Define the two wall
residuals for the traversal at unit speed,
\[
 L_i=qT_g-q\ell_i-d_i-r_{i,t}-r_{i,b} .
\]
For that traversal, inside the footprint of pole $i$, the maximum of
$qt-y+r$ occurs at its bottom vertex and the minimum of $qt-y-r$ at
its top vertex. At the critical geometry the stroke phase increases
strictly. On the outgoing crossing the phase and the radius both
decrease. On the incoming crossing the phase decreases and the radius
increases. Write $a,h$ for the critical half-spans and
$\ell=\sqrt{a^2+h^2}$. The opposite extrema have strict slack because
the departure phases lie between $(1+q)h$ and $2q\ell-(1-q)h$, both
positive, and the arrival phases are their negatives. With phase zero
at the pole midpoint, the two endpoint phase extrema have magnitude
$(1+q)h+\sqrt{1-(b-a)^2}>1$.
Thus the maximizing and minimizing vertices are unique, with strict
one-sided derivative signs there. These properties persist locally by
continuity, including short oblique strokes. Each footprint comprises
one connected portion of the cycle, disjoint from the other footprint.

An arbitrary clock adds a nondecreasing phase shift equal to $q$ times
the delay accumulated since the start of the portion. The maximum
phase span on that wall can grow by at most $q e_i$, where $e_i$ is
the delay within its footprint. Since every phase on a circle of
length $q(T_g+e)$ must be covered, $L_i\le-q(e-e_i)$, and
$e_R+e_L\le e$. Put $L=(L_R+L_L)/2$, $J=(L_R-L_L)/2$ and
$\Phi=L+J^2$. Hence $L\le-qe/2$, $L\le-|J|$ and, close to the
critical point, $\Phi\le-qe/4$.

Let $S$ be the average of $V$ under horizontal reflection with pole
interchange and under height reflection with time reversal and endpoint
interchange. These commuting involutions act orthogonally. Their common
fixed space consists of ordinary tuples $S=S(a,h)$. Put $A=V-S$.
The function $\Phi$ is invariant. At fixed $b$ its Hessian at the
critical tuple is positive definite on the mean-height-zero space.
Indeed
\[
 2L=q(\ell_R+\ell_L)+2q(\ell_{R_bL_t}+\ell_{L_bR_t})
       -d_R-d_L-\sum_j r_j .
\]
All four norm Hessians are positive semidefinite. Each $-r_j$ has
strictly positive second derivative in its own $x_j$, so a vector in
the kernel has all horizontal entries zero. The two nonvertical crossing
lengths then require $\delta y_{R_t}=\delta y_{L_b}$ and
$\delta y_{R_b}=\delta y_{L_t}$. Modulo vertical translation this
leaves $(\delta y_{R_t},\delta y_{R_b},\delta y_{L_t},\delta y_{L_b})
=(u,-u,-u,u)$. Here $D J=-2(1+q)u$. Since $J=0$ at the symmetric
tuple, the Hessian of $J^2$ removes that last kernel. Positive
definiteness persists on a neighborhood. The first derivative of
$\Phi$ at $S$ vanishes in every antisymmetric direction, by averaging.
Consequently $\Phi(S,b)+c\lVert A\rVert^2\le\Phi(V,b)$ for a fixed
$c>0$. At $S$, the root in $b$ of $\Phi=0$ is
\[
 B(a,h;q)=a+\sqrt{1-G^2},\qquad G=2q\sqrt{a^2+h^2}-(1-q)h .
\]
The derivative $\partial_b\Phi(S,b)$ is positive and bounded above
locally. Combining these estimates gives
\begin{equation}\label{eq:asym-cost}
 b\le B(a,h;q)-c_1(\lVert A\rVert^2+e)\qquad(c_1>0).
\end{equation}

\paragraph{A rounded escape on one oblique stroke.}
For a stroke with downward unit direction $u_i$ let $w_-,w_+$ denote
the unit directions of its preceding and following crossings in forward
time. Put $z_\pm=-u_i\cdot w_\pm$. Near the critical tuple these
inclinations lie in $(q,1)$. Let $H(z,q)$ denote the rounded-contact
height in Appendix~\ref{app:middle-comparison}. For clarity its
formula, with $c=\sqrt{1-z^2}$ and $\alpha=\sqrt{1-q^2}$, is
\[
 m=qz+\alpha c,\quad k=\alpha z-qc,\quad M=\sqrt{1+2qm+q^2},
 \quad S_0=qk/M,\quad C_0=(1+qm)/M,
\]
\[
 p=\int_0^{\arcsin S_0}\frac{d\theta}
 {\cos\theta+\sqrt{q^2-\sin^2\theta}},\qquad
 H=p+\frac{1+k/M}{m+q}.
\]
These functions are smooth here. We claim that every actual guard in
the neighborhood satisfies, for each stroke,
\begin{equation}\label{eq:unequal-rounded}
 E_i:=\ell_i-H(z_-,q)-H(z_+,q)\le C e .
\end{equation}
Use the orthonormal frame with vertical axis $-u_i$ and horizontal axis
pointing outward from the corridor. If $E_i>0$, put
$R=1+E_i/[2(H_-+H_+)]$, $H_\pm=H(z_\pm,q)$. Split the stroke at the
point whose distances to the incoming and outgoing endpoints are
$l_\pm= R H_\pm+E_i/4$. Their sum is exactly $\ell_i$. In the
unit-speed clock, set the splitting point and its time to zero.
The future part of the stroke is $(0,-t)$. Use the radius-$R$ copy
of the contact arc $(-\cos\theta,\sin\theta-t)$, where $t$ is its
integral up to $\theta$. Scaling multiplies all times and positions by
$R$ and preserves speed. At $Rp$ continue on the tangent with velocity
$q(k,m)$ for the outgoing inclination. The stroke clearance is at
least $R$. The signed separation from the outgoing crossing in the
unit normal $(k,m)$ equals
\[
 l_+(m+q)-R[k C_0-m S_0+(m+q)p]\ \ge R .
\]
Here $k C_0-m S_0=k/M$ and the normal component of the crossing
velocity equals $q$, so the separation stays constant. Time reversal
and reflection in the frame's vertical coordinate give the incoming
half using its own inclination and $l_-$. Both halves meet at
$(-R,0)$. This construction handles the two different crossing angles
on a single stroke.

The tail and corridor margins are strict at the critical tuple. More
explicitly, set $T_0=2H(1+1/z)$ at the symmetric tuple. There
$a>1$, $0<p<H$, $qmT_0+S_0-(1+qm)p>H+1$ and
$C_0-qk(T_0-p)>0$. The checker \path{anc/slow/check_slow_widths.py}
verifies them at the enclosed root. Up to the midpoint of the next
opposite stroke, that stroke has horizontal distance at least $2a-1>1$.
The constructed path lies strictly inside the corridor. By that midpoint
the path is more than one above every patrol height. After a small
perturbation these strict conditions persist, including upon rotating
back to physical coordinates. Continue vertically upward in physical
space from that time, and use the downward past tail on the incoming
side. Thus the entire unit-clock trajectory has clearance at least $R$,
shrinking the neighborhood if needed for the uniform tail margins.

On the bounded central time interval the physical clock and the
unit-speed clock, aligned at a vertex, change positions by at most a
fixed multiple of $e$. Indeed in each traversal the arclength delay
is nondecreasing and at most $e$. The interval meets a bounded number
of periods since $T\to T_*$. The comparison follows from the unit speed
of the arclength parametrization. On the two infinite tails the global
extrema of the vertex heights give clearance regardless of clock. If
$E_i>C e$ with $C$ sufficiently large, $R-1$ exceeds the possible
central position error. This would give a strict escape against the
actual clock, proving \eqref{eq:unequal-rounded}.

\paragraph{Comparison with an attained symmetric width.}
Both involutions preserve the smooth function $(E_L+E_R)/2$ because
the inclinations are inner products and reversal interchanges the two
neighbors. Its derivative in direction $A$ vanishes at
$S$. By \eqref{eq:unequal-rounded} and Taylor's theorem,
\[
 h-H(d,q)\le C'(\lVert A\rVert^2+e),\qquad
 d=h/\sqrt{a^2+h^2}.
\]
Scale $(a,h)$ by $\lambda=\min(1,H(d,q)/h)$. This keeps $d$ fixed and
$1-\lambda\le C''(\lVert A\rVert^2+e)$. The ordinary finite-leg
separator in Appendix~\ref{app:bent-separator} guards the full symmetric
width $2B(\lambda a,\lambda h;q)$ because the height condition now holds
and all auxiliary margins at the critical tuple are strict. The same
interval check verifies these margins. The wall functional $B$ has
zero gradient in $(a,h)$ at the critical tuple because that tuple
is its interior maximum. Hence along every competitive sequence
\[
 B(\lambda a,\lambda h;q)\ge B(a,h;q)-o(1)(\lVert A\rVert^2+e).
\]
By \eqref{eq:asym-cost} this attained ordinary half-width is strictly
larger than $b$ whenever $A\ne0$ or $e>0$, sufficiently close to
$q_a$. This contradicts $2b\ge W_4(q)$. If $A=e=0$, the original
guard is already ordinary. A sequence violating the asserted
neighborhood would furnish exactly the competitors just excluded.
Taking suprema and using the eight-leg strict comparison below $q_a$
proves Theorem~\ref{thm:four}.

\section{The five-leg construction}\label{app:five-leg}
Normalize $r=V=1$ and write $\varepsilon=q_a-q$. This appendix proves
Theorem~\ref{thm:five}: starting from a symmetric bowtie near $q_a$, we bend
only the right downward stroke inward by a small offset $e\ge 0$ at its
midpoint while lengthening its half-height to $h_R=H+\gamma e$ ($\gamma=1/10$).
Inscribing a slightly smaller disk centered on an auxiliary vertical motion
relaxes the right-stroke fan handoff at first order in $e$, whereas the extra
stroke duration and wall-offset asymmetry cost only $O(e^2)$ in guarded width,
lowering the quadratic deficit below the four-leg limit $C_4$.

\paragraph{A single bent stroke.}
Write $\alpha=\sqrt{1-q^2}$, $c=\sqrt{1-d^2}$, $m=qd+\alpha c$,
$k=\alpha d-qc$ and $M=\sqrt{1+2qm+q^2}$. Define
\[
 \begin{aligned}
 S&=qk/M,\qquad C=(1+qm)/M,\qquad
 p=\int_0^{\arcsin S}\frac{du}{\cos u+\sqrt{q^2-\sin^2u}},\\
 K&=k/M+(m+q)p,\qquad H=(1+K)/(m+q).
 \end{aligned}
\]
These are the ordinary handoff quantities at inclination $d$.
Set $\gamma=1/10$, choose a small inward displacement $e\ge0$, and put
\[
 h_L=H,\qquad h_R=H+\gamma e,\qquad
 h=(h_L+h_R)/2,\qquad L=h/d,\qquad a=Lc.
\]
The right stroke visits $(a,h_R),(a-e,0),(a,-h_R)$. Cross straight
to $(-a,h_L)$, descend straight to $(-a,-h_L)$, and cross straight
to the starting vertex. Both crossings have duration $2L$. Write
$e_R=e$, $e_L=0$, $\rho_i=\sqrt{h_i^2+e_i^2}$, and for $j\ne i$ put
\[
 G_i=q(2L+\rho_j)-h_i,\qquad \delta_i=\sqrt{1-G_i^2}.
\]
The left and right corridor walls are at $-a-\delta_L$ and
$a+\delta_R$. Their different offsets allow each stroke to use the
full wall phase available at its endpoints.

\begin{lemma}[Single-bend five-leg guard]\label{lem:five-cut}
For $(q,d)$ sufficiently close to the critical pair
$(q_a,(1-q_a)/(2q_a))$ and sufficiently small $e\ge0$, this patrol
guards width $2a+\delta_L+\delta_R$ against every upward transit.
\end{lemma}
\begin{proof}
Use coordinates with positive horizontal direction toward the wall
of stroke $i$, and centered time $\sigma$ on that stroke. Its center
has coordinates $(a-e_i+e_i|\sigma|/\rho_i,-h_i\sigma/\rho_i)$.
Compare the actual disk with the vertical motion $(a-e_i,-\sigma)$.
Define
\[
 u_i=\sqrt{(e_i/\rho_i)^2+(1-h_i/\rho_i)^2},\qquad
 r_i=(1+p u_i)^{-1},\qquad p_i=r_i p.
\]
For $|\sigma|\le p_i$ the full disk of radius $r_i$ about the vertical
motion lies in the actual unit disk because $r_i+p_i u_i=1$.
Use the attached eikonal fan of Lemma~\ref{lem:attached-fan}, scaled by $r_i$ in space and time. Its
inward flat is $y=q\sigma$. Time and height reflection give its
incoming half. Connect the contact by a chord to the wall point
of height $-h_i\sigma/\rho_i-R_i$, with
$R_i=\sqrt{1-(\delta_i+e_i)^2}$. This point is in the actual disk
during both halves. The reflected incoming construction uses the
opposite wall offset until the reset across the captured wall chord.

At the outgoing handoff $\sigma=p_i$, use coordinates $(\xi,y)$
relative to the departure vertex. Write $v_i=\rho_i-p_i$.
On the subsequent straight crossing the witness time
$w=(-m\xi+ky)/\alpha$ gives cost
\[
                 qv_i+k\xi+my
\]
whenever $d\xi+cy\ge0$ and $0\le w\le2L$. This is the finite
straight-leg witness identity of equation~\eqref{eq:obs-witness}.
The two fan endpoints have coordinates
$(-e_i-r_i,h_i+qp_i)$ and $(-e_i-r_iC,h_i+r_i(S-p))$.
The fan slopes belong to $[-k/m,0]$. Consequently the cost and the
normal coordinate are nondecreasing along the fan, and the witness
time is nonincreasing. At the critical pair the normal and finite
witness tests at the endpoints are strict. They persist for this
construction. The maximum cost is $mh_i+q\rho_i-ke_i-r_iK$.
It equals one on the left. On the right the derivative at $e=0+$
of the maximum cost, at fixed $(q,d)$, is
\[
          (m+q)\gamma-k+Kp/H .
\]
It is strictly negative at the critical pair since
$(k-Kp/H)/(m+q)>0.10414$. In a neighborhood the right cost therefore
is at most one, strictly for $e>0$. The expression for $u_R$ is smooth
from the right, as follows also from
$u_R=e\sqrt{2/[\rho_R(\rho_R+h_R)]}$.

The wall point at the handoff has the departure witness, by the strict
inequality
\[
 \delta_i^2+(h_i-h_i p_i/\rho_i-R_i)^2<(1-qv_i)^2,\qquad qv_i<1.
\]
At the critical pair this is the strict ordinary test in
Appendix~\ref{app:variation}, and all quantities are continuous.
The adjoining wall offsets, the fan's position inside the strip,
and $0<p_i<\rho_i$ also have strict slack. The interval checks in \path{anc/bends/check_fifth_leg.py}
enclose these tests at the critical root.

For completeness the two handoffs can be joined despite the unequal
stroke heights. Consider the right-to-left crossing in departure
time $t$ and write $D=2L$. Its two unlit flats are
$F_+(t)=q(\rho_R+t)$ and $F_-(t)=q(t-D-\rho_L)$.
At $t=0$ the second one equals $-h_R-G_R$ and at $t=D$ the
first equals $h_L+G_L$. On $F_+$ interpolate a point from
transverse coordinate $a-e_R-r_R$ at $t=-v_R$ to $-a$ at $t=D$.
On $F_-$ interpolate from $a$ at $t=0$ to $-a+e_L+r_L$ at
$t=D+v_L$. Both segments belong to the convex capture set of this
one straight leg by the endpoint witnesses and the reflected handoff
at the other end. On $[0,D]$ the first point lies to the left of the
second. Their difference is affine and strictly negative at both
endpoints for all sufficiently small perturbations. Connect them
by a captured chord, extending the upper flat to the left and the
lower flat to the right. On the short intervals outside $[0,D]$ use
the chord to the wall point interpolating from the handoff to the
departure wall offset $-G_R$, or from the arrival wall offset $+G_L$
to the incoming handoff. The change of chords at departure or
arrival sweeps a triangle in the same convex capture set. At the
incoming fan the change of sign of the wall offset sweeps a triangle
in the actual disk. These statements use the respective local
quantities at the two ends. Reflect transversely and interchange
$R,L$ for the other crossing. This gives a bounded periodic separator.
Every uncaptured portion is a rising flat or solves the eikonal
equation $f_t=q\sqrt{1+f_x^2}$, so a speed-$q$ path cannot pass
from below to above.
\end{proof}

\paragraph{The width coefficient.}
For a Taylor expansion use independent variables $x=a$, $d$, $e$,
with $h=xd/c$, $h_{R,L}=h\pm\gamma e/2$. Let
$G=q(2L+h)-h$ and $\delta=\sqrt{1-G^2}$. At fixed $x,d$ the mean
change in $G_i$ is $q e^2/(4h)+O(e^3)$, and their first-order
changes are $\mp(1+q)\gamma e/2$. Thus the additional full-width
loss caused by the bend at fixed mean geometry is $M_e e^2+O(e^3)$,
where
\[
          M_e=\frac{Gq}{2\delta h}
                     +\frac{(1+q)^2\gamma^2}{4\delta^3}.
\]
There are no mixed quadratic terms with $e$ because the first
derivative in $e$ vanishes identically at zero. Write $M_x,M_d$ for
the positive diagonal entries of the negative ordinary half-width
Hessian at the critical wall optimum, as in Appendix~\ref{app:variation}.
By Taylor's formula for $2B$, the full-width loss from the moving wall optimum to second order is
$M_x(\Delta x)^2+M_d(\Delta d)^2+M_e e^2$. The design fixes
$(h-\gamma e/2)/H=1$. Let $g_x,g_d$ be the ordinary gradient of
$h/H$ at the critical wall optimum, $g_e=-\gamma/(2H)$, and let
$E_0>0$ be minus the speed derivative of $h/H$ along the wall-optimal
branch there. For $j=x,d,e$ set
\[
         B=\sum_j g_j^2/M_j,\qquad v_j=-E_0g_j/(M_j B).
\]
Take $d=d_*(q)+v_d\varepsilon$, $e=v_e\varepsilon$ and
$x=(H+\gamma e/2)c/d$. Since $v_e>0$ this is the five-leg design.
Differentiating its scale identity gives
$x=x_*(q)+v_x\varepsilon+O(\varepsilon^2)$. Its limiting loss is
$C_5=E_0^2/B$. The four-leg value has the same numerator and uses
only $x,d$ in the denominator. The new term is strictly positive.
The exact interval replay gives the enclosure in Theorem~\ref{thm:five}.
Comparing the two limits and applying Theorem~\ref{thm:four}
proves the strict result against all periodic four-leg guards in a
common left neighborhood of the transition.

At $q=0.45789$ the patrol in Figure~\ref{fig:fifth-leg} uses $t=0.0000020366515$,
$d=d_*(q)-3.702797t$ and $e=1.029609t$. Direct interval checks
of the fan and wall tests give width above $3.36787047176625r$.
Every ordinary bowtie at that speed has width below
$3.36787047176461r$. To obtain the second number, first discard
ordinary competitors with deficit greater than $52r\varepsilon^2$.
The two wall localization identities in Appendix~\ref{app:bent-uniform}
put the remaining angles within $8\varepsilon$ of $d_*$ and their
amplitudes within $10\varepsilon$ of the best wall amplitude at that
angle. On the resulting rectangle the rounded escape forces $h/H\le1$.
The lower Hessian bounds for the wall deficit and upper gradient
bounds for $h/H$ give the stated estimate by weighted Cauchy--Schwarz,
using the wall-optimal $h/H-1$ evaluated directly at this speed.
The arbitrary-four-leg theorem supplies an existential left
neighborhood. The numerical upper bound at this particular speed
uses the ordinary class.

\paragraph{Bending the second stroke.}
The same handoff proof allows both strokes to be bent by $e$ with
equal heights $h_L=h_R=H+\gamma e$. The two wall chords now have
equal insets. The first derivative of each handoff cost is the same
negative number tested above. At fixed $x,d$ the full-width loss
is $[Gq/(\delta h)]e^2+O(e^3)$, since both strokes have duration
$2\sqrt{h^2+e^2}$. The scale constraint is $(h-\gamma e)/H=1$.
Thus replacing $M_e$ by $Gq/(\delta h)$ and $g_e$ by $-\gamma/H$
in the calculation above gives a six-leg construction with
$50.35196081<C_6<50.35196085$. The comparison $C_6<C_5<C_4$
quantifies the benefit of the second bend for these designs.
Determining the optimal five-leg and six-leg coefficients is a
separate variational problem.

\section{The slow-speed ordinary maximum}\label{app:slow-ordinary}
We prove Theorem~\ref{thm:slow-ordinary} in four steps: interior-lane phase
tests rule out reversed bowties and off-scale inclinations, asymptotic wall and
lane estimates localize every competitive ordinary bowtie to centered descending
strokes with inclination $d=q+uq^3$ ($u\in[3/4,5/4]$), the rounded-contact height
gives the exact maximum half-width $B_H(q,u)$ at each such angle, and an analytic
expansion around $u(0)=1$ determines $W_4(q)/r$ through $O(q^4)$.

\subsection{Two crossing lines on the same lane}
In units $r=V=1$, write $a$ for the half-distance between columns and $h$
for the half-height of a stroke. Put $\ell=\sqrt{a^2+h^2}$, $d=h/\ell$,
$c=a/\ell$ and $\sigma=+1$ for downward strokes or $-1$ for upward
strokes. For corridor half-width $b$ and transverse offset $\xi$, put
$b_f=b+|\xi|$ and $\delta=b_f-a$. A necessary condition is $0\le\delta\le1$.
The case $a=0$ gives $b\le1$. Define
\[
 P=4q(h+\ell),\qquad n=d-\sigma q,\qquad S_\sigma=\sqrt{c^2+n^2}.
\]
Time zero is the center of the right stroke and the spatial axis is
midway between the columns.

\begin{lemma}[Interior-lane phase constraint]\label{lem:ordinary-interior}
If $a>1$, guarding requires $aq(1+d)\le S_\sigma$. Whenever
$P/2>S_\sigma/c$, it also requires
\begin{equation}
       aq(1+d)+(a-1)|d-\sigma q|\le S_\sigma .                 \label{eq:ordinary-lane-corr}
\end{equation}
\end{lemma}
\begin{proof}
In the straight-runner coordinates $(x,qt-y)$ the two crossing lines
on the lane at $x$ have phase centers
$\pm(\sigma P/4+nx/c)$ modulo $P$. Indeed their phases at the
right column have magnitudes $(1+\sigma q)h$, and
$(1+\sigma q)h-na/c=\sigma P/4$. On $0\le x<a-1$ neither stroke
can reach the lane. A diagonal's supporting line reaches phases at
vertical distance at most $R=S_\sigma/c$ from its center. The two
phase intervals must cover the entire circle even when the finite
diagonals are extended to their supporting lines.

At $x=0$ the centers are opposite, so $P/4\le R$. If $P/2>R$,
coincident centers would leave an open gap. Starting at zero and
moving continuously along the lane parameter therefore gives
$|n|x/c<P/4$ throughout. The farthest phase then has distance
$P/4+|n|x/c$ from the nearer center. Let $x\uparrow a-1$ and
multiply by $c$. This proves \eqref{eq:ordinary-lane-corr}, including off-center corridors.
\end{proof}

\begin{corollary}[Reversed and two-way ordinary bowties]
The reversed orientation can guard at most $(1/q+4)r$ for $0<q<1$.
Its best width, and the best width of an ordinary patrol that guards
both transit directions, are each asymptotic to $r/q$.
\end{corollary}
\begin{proof}
Take $\sigma=-1$ and $a>1$. If $P/2\le S_-/c$ then
$a\le S_-/[2q(1+d)]\le(1+q)/(2q)$. Otherwise \eqref{eq:ordinary-lane-corr} gives
\[
 a\le\frac{S_-+q+d}{2q+(1+q)d}
 \le\frac{1+2q+d}{2q+(1+q)d}\le1+\frac1{2q}.
\]
Use $b\le a+1$ and restore units. The case $a\le1$ also satisfies
the bound. A two-way ordinary patrol must guard in the unfavorable
direction. The zero-height transverse sweep gives the common lower
bound $r\sqrt{4+q^{-2}}$.
\end{proof}

\begin{lemma}[Unequal straight crossings]\label{lem:ordinary-unequal}
Allow the strokes on $x=\pm a$ to have different heights and vertical
centers. For crossing $i$, write its supporting line in runner coordinates
as $qt-y=\gamma_i+s_i x$, modulo the phase advance $P$ of a period.
Put $R_i=\sqrt{1+s_i^2}$. If $a>1$ and $\max_i R_i<P/2$, guarding
requires distinct center-lane phases. With $g_0\in(0,P)$ representing
$\gamma_2-\gamma_1$ modulo $P$, it requires
\[
 \frac P2+\left|g_0-\frac P2\right|+(a-1)|s_2-s_1|\le R_1+R_2 .
\]
\end{lemma}
\begin{proof}
On $|x|<a-1$ the two crossing lines give intervals of radii $R_1,R_2$
on the phase circle. Their complements are open arcs centered opposite
the original centers, with radii $P/2-R_i$. To cover the circle those
open arcs must be disjoint. If $g$ is the oriented separation of the
line centers in $[0,P)$, this is equivalent to
$\min(g,P-g)\ge P-R_1-R_2>0$. Thus $g$ stays in $(0,P)$ and equals
$g_0+(s_2-s_1)x$ throughout the interval. Equivalently
$P/2+|g-P/2|\le R_1+R_2$. Taking the supremum over $x$ proves the
claim. For ordinary bowties $g_0=P/2$ and $|s_2-s_1|=2|n|/c$,
recovering \eqref{eq:ordinary-lane-corr}.
\end{proof}

\subsection{Localizing every competitive ordinary patrol}
It suffices to consider a patrol at least as wide as the explicit ordinary
guard of Theorem~\ref{thm:slow-limit}, so $b\ge q^{-1}-1+3q/2-O(q^2)$.
Put $t=q^{-1}-a$. Since $b_f=a+\delta$ we have
$t\le1+\delta+O(q)$. The center lane gives $a\le(1+q)/q$ for $a>1$,
so $t$ is bounded, $aq=1+O(q)$ and $d=O(q)$.
Consequently $P/2\to2$, $S_\sigma/c\to1$, and \eqref{eq:ordinary-lane-corr} gives
$|d-\sigma q|=O(q^2)$. Only the forward orientation can compete.
All estimates here are uniform over patrols at or above the benchmark.

Write $y=(d-q)/q^2$, which is bounded. Expanding \eqref{eq:ordinary-lane-corr} to first
order gives $t\ge1+|y|-O(q)$, and hence $|y|\le\delta+O(q)$.
The farther-wall midpoint test of Appendix~\ref{app:lane-tests} uses
\[
 G=\ell[2q-(1-q)d]=1+(1-t-y)q+O(q^2).
\]
For clarity, at that midpoint the adjacent departing diagonal has
squared distance $(\delta+cw)^2+(G-(q-d)w)^2$ at elapsed time $w$.
The arriving diagonal is its reflection and the tested phase is midway
between the two right strokes. For $G>0$ each of these diagonals
ranges in phase between $(1+q)h$ and $G$, both in $[0,P/2]$ on the
departure. Other copies cannot improve the distance, and the opposite
stroke is out of horizontal range. Capture thus requires either
$\delta^2+G^2\le1$, or $d<q$ with $c\delta<(q-d)G$ and distance at
most one to the departing supporting line. In the latter case writing
$e=q-d$ gives $a[2q-(1-q)d]+e\delta\le S_+$ and $\delta<e/S_+$.

In the endpoint case $\delta^2\le2q(t+y-1)+O(q^2)
\le4q\delta+O(q^2)$, giving $\delta=O(q)$. The other case already
gives $\delta=O(q^2)$. Thus $t=1+O(q)$ and $y=O(q)$. Set
\[
      T=(t-1)/q,\qquad u=(d-q)/q^3,\qquad D=\delta/q.
\]
These three variables are bounded. At second order \eqref{eq:ordinary-lane-corr} and the wall
gap become
\[
 T\ge |u|-\tfrac12-O(q),\qquad
 G=1-(T+u+\tfrac12)q^2+O(q^3).
\]
Indeed $aq(1+d)=1-(1+T)q^2+O(q^3)$, $(a-1)|d-q|=|u|q^2+O(q^3)$
and $S_+=1-q^2/2+O(q^4)$. In the supporting-line case $D=O(q^2)$,
contradicting the benchmark $D-T\ge3/2-O(q)$. In the endpoint case
\[
 D^2\le2(T+u+\tfrac12)+O(q)\le4T+2+O(q),\qquad
 D-T\le\tfrac32-(D-2)^2/4+O(q).
\]
The benchmark now forces $D\to2$, $T\to1/2$ and $u\to1$.
In particular every competitive shape eventually has $d>q$ and
$3/4\le u\le5/4$. This step rules out other aspect ratios or traversal
directions that could have changed the class expansion.

\subsection{The exact maximizing branch}
For $u\in[3/4,5/4]$ put $d=q+u q^3$ and $\alpha=\sqrt{1-q^2}$. Use the contact quantities
\[
 m=qd+\alpha c,\quad \kappa=\alpha d-qc,\quad M=\sqrt{1+2qm+q^2},
 \quad S=q\kappa/M,\quad C=(1+qm)/M .
\]
Let $p$ be the integral of
$(\cos\theta+\sqrt{q^2-\sin^2\theta})^{-1}$ from zero to $\arcsin S$.
Set
\[
 H=p+\frac{1+\kappa/M}{m+q},\quad A_H=cH/d,\quad
 k=\frac{2q-(1-q)d}{c},\quad G_H=k A_H,\quad
 B_H=A_H+\sqrt{1-G_H^2}.
\]
We claim that $B_H$ is the \emph{exact maximum half-width at this angle}
in the forward orientation for all sufficiently small $q>0$, uniformly
on the displayed $u$ interval.

For necessity use the homogeneous rounded escape of Appendix~\ref{app:middle-comparison}. Set $T_0=2H(1+1/d)$. Its four conditions at $h=H$ are
\[
 0<p<H,\quad A_H>1,\quad
 qmT_0+S-(1+qm)p>H+1,\quad C-q\kappa(T_0-p)>0.
\]
Here $\kappa=u q^3+O(q^5)$, $p=u q^4+O(q^5)$, $H=1-q+O(q^2)$,
$A_H=q^{-1}-1+O(q)$ and $T_0=2/q+O(q)$. The vertical-tail excess
is $q+O(q^2)$ and the last margin tends to one. Thus the complete
exit stays between the columns. More explicitly, the runner follows
the unit-radius contact arc on the descending stroke, follows the
tangent of velocity $q(\kappa,m)$ to $T_0$, and continues vertically
at speed $qm$. Its signed clearance on the intervening crossing is
$H(m+q)-[\kappa C-m S+(m+q)p]=1$. Time and height reflection give
the incoming half. Multiplying all lengths and times by $h/H>1$
would turn this trajectory into a strict escape. Hence $a\le A_H$.

Conversely the ordinary finite-leg separator at $h=H$ has handoff
at equality and wall inset $\sqrt{1-G_H^2}$. Its remaining tests
hold uniformly. In fact $G_H=1-2u q^2+O(q^3)$, $p>0$, $H/2-p\to1/2$,
$2A_H+\sqrt{1-G_H^2}>1$, $q(H-p)\to0$, and
$2(G_H-q)-(1-q^2)(H-p)\to1$. At both fan endpoints
$(X,Y)=(-1,H+qp),(-C,H-p+S)$ we have $dX+cY\to1$ and
$(-mX+\kappa Y)/\alpha\to1$, strictly inside the crossing time
interval $[0,2H/d]$. These are the single-crossing tests of
Appendix~\ref{app:bent-separator}.

For $b>1$ the farther-wall bound at these angles is
$b_f\le B(a)=a+\sqrt{1-(ak)^2}$.
This concave function has
$B'(A_H)=1-A_Hk^2/\sqrt{1-G_H^2}
=1-1/(2\sqrt u)+O(q)>0$. Therefore its maximum on $[0,A_H]$
is attained only at $A_H$. The centered patrol there guards $B_H$,
proving the claim. Equality in the full width also forces $\xi=0$.

At the class level the localization above gives
$W_4(q)=2\max_{3/4\le u\le5/4}B_H(q,u)$ in the normalized units.
To analyze the apparent singularity use
\[
 p=S\int_0^1
 \frac{dz}{\sqrt{1-S^2 z^2}
 [\sqrt{1-S^2 z^2}+q\sqrt{1-(\kappa/M)^2 z^2}]} .
\]
The integrand is analytic in $q,u$ near $q=0$ with uniformly positive
denominators on the integration interval, and its Taylor series there
converges locally uniformly. Thus $H$ is analytic. The inset equals
$q\sqrt{(1-G_H^2)/q^2}$ for $q>0$, where the radicand tends to $4u>0$.
It follows that
\[
 f(q,u)=\frac{B_H(q,u)-q^{-1}+1}{q}
\]
extends analytically across $q=0$. Its leading term
$f_0(u)=1/2-u+2\sqrt u$ has its unique maximum at $u=1$.
Its second derivative $-1/(2u^{3/2})$ is uniformly negative on
$[3/4,5/4]$, and its derivatives at the two endpoints have opposite
signs. These properties persist for small $q$. The implicit function
theorem at $(0,1)$ gives a unique analytic maximizing $u(q)$.

\subsection{Width through fourth order}
Since $S=O(q^4)$ and $\kappa/M=O(q^3)$ the integrand above differs
from $1/(1+q)$ by $O(q^7)$. This gives
$p=S/(1+q)+O(q^{11})$. Also $m=1+O(q^6)$ and
$\kappa=u q^3+u q^5/2+O(q^6)$. Through order five the height is
\[
 H=\frac{1+u q^3+u q^5/2}{1+q}+O(q^6).
\]
Write $f(q,u)=\sum_{i=0}^3 f_i(u)q^i+O(q^4)$. Expansion of $B_H$ gives
\begin{align*}
 f_0&=\tfrac12-u+2\sqrt u,\\
 f_1&=2u-\tfrac32\sqrt u-\tfrac12,\\
 f_2&=u^2-\tfrac52u-2u^{3/2}+\tfrac{7}{16}\sqrt u+\tfrac38,\\
 f_3&=-2u^2+\tfrac52u+\tfrac72u^{3/2}-\tfrac{59}{64}\sqrt u-\tfrac38 .
\end{align*}
Inserting $u(q)=1+u_1q+u_2 q^2+O(q^3)$ in $f_u=0$ gives
$u_1=-f_1'(1)/f_0''(1)=5/2$ and
$u_2=-(f_2'(1)+u_1 f_1''(1)+u_1^2 f_0'''(1)/2)/f_0''(1)=0$.
Substitution into $2B_H=2/q-2+2qf$ proves the width expansion in
Theorem~\ref{thm:slow-ordinary}. This argument gives an asymptotic neighborhood
of zero with no explicit endpoint.

The file \path{anc/slow/check_slow_bowties.py} checks the series algebra with
exact rational coefficients. It separately integrates the original arc
in $\theta$ with 75-digit arithmetic, locates a stationary point, and
prints the normalized cubic and quartic residuals. Those numerical
comparisons supplement the class localization and analytic proof.

\section{Height and flux of the actual reachable set}\label{app:reach}

This appendix establishes the bounded-variation height and circle-flux balances
for the actual intruder-reachable set, together with a wall-visit count, that
underpin the slow-speed lower bounds in Appendices~\ref{app:oneway-loss},
\ref{app:direction-loss}, and~\ref{app:direction-rate} (for Theorems~\ref{thm:loss}
and~\ref{thm:direction-value} and Corollary~\ref{cor:two-way-height}). Throughout
Appendices~\ref{app:reach}--\ref{app:direction-rate},
we normalize $r=V=1$ and translate the corridor horizontally to $J\times\R$ with
$J=[0,W]$. Let $g=(X,Y)$ be a $T$-periodic polygonal guard of speed at most one
on $[0,W]\times\R$. Put $Q=\sqrt{1+q^2}$, $R=|g'-qe_y|$,
$\delta=2Q-qW$, $K=\delta/q$, and let $\langle\cdot\rangle$ denote
the period average. Choose $H_0$ with $|Y|\le H_0$. We use standard
relative perimeter, BV slicing, traces, the chain rule and coarea
\cite{AFP,Federer}. All perimeter calculations take place in a finite
vertical truncation containing the nonconstant part of the indicators.

\subsection{The space-time boundary and its free speed}
Let $E$ consist of the endpoints reachable from some earlier point
with $y<-H_0-2$, by a finite path with Lipschitz constant strictly
less than $q$ avoiding the closed unit disk. The set is open relative
to the hallway: an avoiding finite path has positive clearance and
speed slack, so its endpoint and terminal time can be perturbed on a
short final subinterval. Projection onto the hallway preserves the speed
bound. Waiting gives a complete lower half-strip. Successful guarding
excludes a complete upper half-strip, since a surviving path reaching
$y>H_0+2$ extends to a transit with safe vertical tails. The set is
periodic because all starting times are admitted.

Here are details establishing finite relative perimeter at the original
radius. Off the disk world tube the indicator obeys
$(\partial_t+v\cdot\nabla)\mathbf1_E\ge0$ in distributions for every
$|v|<q$, by propagation along the corresponding straight lines.
Reflect across the two walls and extend periodically in $x$ with
period $2W$. Folding a straight segment gives an admissible path,
so propagation holds off the reflected obstacles. In the original
hallway its own center is the nearest reflected center. Convolve
in time and space with a nonnegative smooth kernel of radius
$\epsilon$. At distance greater than $1+2\epsilon$ from the original
center the convolution $u$ satisfies $u_t\ge q|\nabla u|$.
Take $r'>1+2\epsilon$ and a regular value $c\in[1/3,2/3]$, also
regular for restrictions to the wall and the finitely many affine
circle cylinders. Intersect $\{u>c\}$ with the exterior of this
radius-$r'$ tube. Its free outward normal satisfies
$-N_t\ge q|N_{\rm sp}|$. On the active cylinder the outward vector
surface element is $(g'\cdot n,-n)\,dt\,ds$, where $n$ is the outward
disk normal. Integrating the constant time field over one period
bounds the free integral of $|N_t|$ by $2r'T$. The cone inequality
then bounds that of $|N_{\rm sp}|$ by $2r'T/q$. The whole cylinder
has area at most $2\pi r'\sqrt2\,T$. These estimates give a uniform
relative perimeter bound as $\epsilon\to0$ and $r'\downarrow1$.
The superlevel indicators converge in $L^1$ to $\mathbf1_E$ since
$|\mathbf1_{\{u>c\}}-\mathbf1_E|\le3|u-\mathbf1_E|$ almost everywhere.
Cutting out the enlarged disk adds at most $\pi((r')^2-1)T$ of error.
Lower semicontinuity proves finite perimeter. On the original tube
the interior trace is zero, and the exterior indicator trace is zero
or one almost everywhere. The trace-jump formula gives precisely the
same vector surface element on the trace-one part, denoted $\Gamma$.
The cylinder seams have zero surface area.

On the free boundary the exact normal speed is $q$. Indeed in a small
convex space-time box away from the walls and sensor, a neighborhood of
each free topological boundary point has $E=\{(t,z):t>f(z)\}$.
The strict forward cone from the point is reachable by approximation
and concatenation, and the strict backward cone is unreachable.
Waiting and the same cones between nearby fibers give
$|f(z)-f(w)|\le |z-w|/q$. At a differentiability point $z$ take a
path reaching $(f(z)+h^2,z)$ and its predecessor $z_h$ a time $h$
earlier. For small $h$ the predecessor exists: all $y<-H_0-1$
are reachable by a safe vertical extension, so the free boundary
has a positive spatial buffer from the source $y<-H_0-2$.
Then $|z_h-z|\le qh$ and $f(z_h)\le f(z)+h^2-h$. Differentiability
implies $q|\nabla f|\ge1$. The reverse inequality was the Lipschitz
bound. The graph area formula now gives $-N_t=q|N_{\rm sp}|$ almost
everywhere on the reduced free boundary.

\subsection{The two height derivatives}\label{app:height-balances}
Use the lower orientation and define
$H(x,t)=\int(\mathbf1_E(t,x,y)-\mathbf1_{\{y<0\}})\,dy$.
At regular times write $\Sigma_t$ for the free spatial boundary with
outward normal $\nu$, and $\Gamma_t$ for the active circle trace.
Set $j=(g'-qe_y)\cdot n$. Write $\varepsilon_t,\beta_t$ for the
transverse projections of $(1-\nu_y)\,ds,\nu_x\,ds$ on $\Sigma_t$,
and $\gamma_t,\kappa_t$ for those of $j\,ds,n_x\,ds$ on $\Gamma_t$.
Put $\eta_t=\varepsilon_t(J)$. For a transverse Borel set $I$,
\[
 \mu_t(I)=\int_{\mathbb S^1}\mathbf1_{\{X+n_x\in I\}}j_+\,ds,
 \qquad m_t=2R-\gamma_t(J),\qquad b(t)=q\eta_t+m_t .
\]
The measure $\mu_t-\gamma_t$ is nonnegative on the whole line when
the active trace is extended by zero. Its mass is $m_t$.
The space-time derivative and the derivative on a time slice obey
\begin{equation}\label{eq:reach-time-new}
 D_tH=q\,dx\,dt+q\varepsilon_t\,dt-\gamma_t\,dt,\qquad
 D_xH(\cdot,t)=-\beta_t+\kappa_t .
\end{equation}
To see this, slice the free boundary in time using its normal speed.
The time derivative of the spatial indicator contributes $q\,ds\,dt$
there and $-g'\cdot n\,ds\,dt$ on the circle. Integrating the vertical
derivative up a full lane gives $-1$, hence
$\operatorname{pr}_{x\#}(\nu_y\,ds|_{\Sigma_t})
=dx+\operatorname{pr}_{x\#}(n_y\,ds|_{\Gamma_t})$.
Splitting $1=\nu_y+(1-\nu_y)$ proves the time formula. The horizontal
derivative follows by direct projection. On each lane the projected
active $n_y$ density is between $-1$ and $1$ by pairing the two
circle points. Thus for every interval $I$ of length $d$,
\begin{equation}\label{eq:reach-space-new}
 \int_{\Sigma_t\cap\{x\in I\}}|\nu_x|\,ds
 \le\sqrt{\eta_t(4d+\eta_t)} .
\end{equation}
This is Cauchy--Schwarz using
$\nu_x^2=(1-\nu_y)(1+\nu_y)$ and
$\int_{\Sigma_t\cap\{x\in I\}}\nu_y\,ds\le2d$.
The additional horizontal variation from the active circle is at
most four. Time-direction BV slicing disintegrates $D_tH$ over $dx$.
Since $\mu_t$ has a lane density, the positive measure
$q\varepsilon_t\,dt+(\mu_t-\gamma_t)\,dt$ on the open hallway
also disintegrates over $dx$. Its mass per unit time averages at most
$\langle b\rangle$. This justifies using the time balance both on
bands and on almost every individual lane.

Integrating over the period yields
\[
 \langle b\rangle=2\langle R\rangle-qW\le\delta,\qquad
 \langle\eta_t\rangle\le K,\qquad
 \langle(R-Q)^2\rangle\le Q\delta .
\]
Here $\langle R^2\rangle\le Q^2$ follows from periodicity of $Y$,
and $2\langle R\rangle\ge qW$ from nonnegativity. More locally
$\langle\mu\rangle-q\mathbf1_J dx$ is a nonnegative measure on the
line of mass at most $\delta$. If $B\subset J$ is a union of intervals,
all flux at times with $X\in B$ lies in $B+[-1,1]$. For
$p=\Pr(X\in B)$ and $d=|J\cap(B+[-1,1])|$ we obtain
$2\langle R\mathbf1_B(X)\rangle\le qd+\delta$.
The left side is at least $2Qp-2\sqrt{Q\delta p}\ge Qp-\delta$,
so $p\le(qd+2\delta)/Q$.

In particular, if $q\downarrow0$ and $qW\to2$, the law of $2X/W$
converges to the uniform law on $[0,2]$. For a bounded Lipschitz test
function replace its argument $2(X+n_x)/W$ in the circle integral
by $2X/W$ at an error at most $2\operatorname{Lip}(f)/W$ times
the total mass. The local measure comparison and $R\to1$ in mean
then give $2\langle f(2X/W)\rangle\to\int_0^2 f(s)\,ds$.
For the upper orientation reflect physical height, including the
guard velocity. The corresponding norm is $|g'+qe_y|$, and all
estimates hold with the same $K$.

\subsection{Counting visits to the wall layers}\label{app:wall-access}
Write $P=qT$ and suppose $W>6$. At the left wall let $N$ count components
of $\{X\le2\}$ reaching $\{X\le3/2\}$ in one period, cutting at a time
with $X>2$. The right wall requires such a time. We prove
\begin{equation}\label{eq:wall-access-new}
                         N/P\le 2^{14}K .
\end{equation}
For $1\le a\le3/2$ let $N(a)$ count the components of $\{X\le a\}$
reaching $\{X\le1\}$. Set $V(a)=\int_{\{X\le a\}}|q-Y'|\,dt$ and
$d(a)=(1-V(a)/P)_+$. Denote the maximum wall-footprint half-length
on each such component by $r_j$. Its detected phases have span at most
the variation of $qt-Y$ plus $2r_j$. Wall coverage gives
$P\le V(a)+2\sum_j r_j$. The area of its largest disk cap outside
the left wall is $\arcsin r_j-r_j\sqrt{1-r_j^2}\ge2r_j^3/3$,
by differentiation. The cap begins at zero, and the outside positive
circle flux dominates its area increase, using the disk with center
$(X,Y-qt)$. Total outside flux is at most $\delta T=KP$, so
Hölder's inequality gives $d(a)^3\le12K(N(a)/P)^2$.

For any lane interval $I$, pairing the two circle points at
$z=x-X$ gives their flux density as in \eqref{eq:both-kernel}.
Integrated over a period restricted by any Borel set of center
positions, the signed $X'$ term vanishes. After integration in $z$ it
is $X'$ times a bounded function of $X$, and its enclosed primitive
returns to the starting value. This also works for the outside
half-line $I=(-\infty,0)$. The remaining maximum there is at least
$(1-X)_+|Y'-q|$ after integration in $z$, hence $V(1/2)\le2KP$.
For $I=(1/2,3/2)$ all centers in $(1/2,a]$ project onto the entire
interval. They supply at least $V(a)-2KP$ flux. The total on $I$
is at most $P(1+K)$ by Appendix~\ref{app:reach}. The residual from
centers above $a$ is therefore at most $P[d(a)+3K]$.

For centers in $(1,3/2)$, the transverse coefficient
$\int_{I-X}|z|/\sqrt{1-z^2}\,dz$ is at least
$a_0=2-\sqrt3$ (the minimum is on the centered interval).
Every counted component at level $x$ contributes two crossings of that
level for almost every $x$. Coarea and the residual bound give
$2a_0\int_a^{3/2}N(x)\,dx/P\le d(a)+3K$. Consequently
\[
 \int_a^{3/2}d(x)^{3/2}\,dx\le C_0\sqrt K\,[d(a)+3K],
 \qquad C_0=\sqrt{12}/(2a_0)=3+2\sqrt3 .
\]
For $K=0$ the cap bound already gives $d=0$. For $K>0$ write
$z(a)=d(a)+3K$ and $F(a)=\int_a^{3/2}z(x)^{3/2}\,dx$. Since
$(u+v)^{3/2}\le\sqrt2(u^{3/2}+v^{3/2})$ and $z\ge3K$,
$F(a)\le C_1\sqrt K\,z(a)$ with $C_1=\sqrt2 C_0+\sqrt{3/2}<11$.
Integrating $F'=-z^{3/2}\le-F^{3/2}/(C_1^{3/2}K^{3/4})$ yields
$F(a)\le4C_1^3 K^{3/2}/(a-1)^2$. Monotonicity of $z$ gives
$\tfrac14 z(3/2)^{3/2}\le F(5/4)$, hence
$z(3/2)\le4\,32^{2/3}C_1^2K<8192K$.

Each component counted by $N$ crosses $(3/2,2)$ twice. Its
contribution to the transverse coefficient for the same lane interval
is at least $2\int_{3/2}^2\sqrt{1-(x-3/2)^2}\,dx
=\sqrt3/4+\pi/6>1/2$. Comparing the summed crossing costs with
the residual $Pz(3/2)$ proves \eqref{eq:wall-access-new}, including
threshold tangencies and the case $K=0$. Horizontal reflection gives
the identical right-wall estimate.

\section{A positive one-way additive loss}\label{app:oneway-loss}
We prove Theorem~\ref{thm:loss}, showing that $D_1(q)=2/q-W_1(q)/r\ge 2^{-21}$
for every $0<q\le 2^{-24}$. Under a smaller deficit $K$, comparing the mean
reachable heights of two adjacent interior bands prevents deep transverse turns
and forces net progress exceeding $5/2$ across nearly every length-$16$ subinterval
of a window of duration $16\lceil 2/q\rceil$, driving the guard across a distance
wider than the corridor itself.

\subsection{Reachable heights and wall visits}
Normalize $r=V=1$ and first take a periodic polygonal
guard $g=(X,Y)$ of period $T$ covering $[0,W]$. Write
\[
 Q=\sqrt{1+q^2},\qquad R=|g'-q e_y|,\qquad
 \delta=2Q-qW,\qquad K=\delta/q,\qquad
 \langle f\rangle=T^{-1}\int_0^T f(t)\,dt .
\]
We use the estimates proved in Appendix~\ref{app:reach}. Let $\mathcal R$
be the space-time set reachable from far below the guard by an
avoiding path of speed strictly less than $q$. Its renormalized
height is $H(x,t)=\int(\mathbf1_{\mathcal R}-\mathbf1_{\{y<0\}})\,dy$.
On its free spatial boundary let $\nu$ be the outer unit normal and
$\eta_t=\int(1-\nu_y)\,ds$. For the moving unit disk write $n$ for
the outward normal, $\Gamma_t$ for its active reachable-boundary trace,
and $j=(g'-qe_y)\cdot n$. The positive circle flux into transverse lanes
$I$ and the trace error are
\[
 \mu_I(t)=\int_{\mathbb S^1}\mathbf1_{\{X+n_x\in I\}}j_+\,ds,\qquad
 m_t=2R-\int_{\Gamma_t}j\,ds\ \ge 0,\qquad b(t)=q\eta_t+m_t .
\]
The reachable-boundary and localized occupation estimates give
\begin{equation}\label{eq:loss-budgets}
 \langle b\rangle\le qK,\qquad \langle\eta_t\rangle\le K,\qquad
 \langle(R-Q)^2\rangle\le Q\delta,\qquad
 \Pr\{\min(X,W-X)\le2\}\le(6q+2\delta)/Q
\end{equation}
for $W>6$, with uniform time modulo $T$. In fact
$\langle b\rangle=2\langle R\rangle-qW$ and
$\langle R^2\rangle\le Q^2$. Equation~\eqref{eq:wall-access-new}
gives at most $2^{14}qTK$ components of $\{X\le2\}$ reaching
$X\le3/2$ per period, and the same bound at the right wall. Its proof
counts crossings of the layer $(3/2,2)$ and charges their disk flux
to the lane $(1/2,3/2)$. For a uniform starting phase $s$, it follows
that a window $[s,s+L]$ reaches either $3/2$ wall layer with probability
at most
\begin{equation}\label{eq:loss-wall}
             (6q+2\delta)/Q+2^{15}qLK .
\end{equation}
Indeed, it either starts in a width-two wall layer or a counted component
starts within the window. Counting component starts with multiplicity
works for $L>T$ as well.

We will use two ordered, disjoint open bands $I_0,I_1$ of common width
$w$, whose convex hull has length $\ell$. Let $S(t)$ be the difference
of their mean reachable heights. When $I_0$ stays outside the disk
footprint on $[u,v]$, the storage inequalities are
\begin{align}
 \int_u^v\mu_{I_1}(t)\,dt
     &\le w(|S(u)|+|S(v)|)+\int_u^v b(t)\,dt, \label{eq:loss-storage}\\
 |S(t)|&\le 2+\sqrt{\eta_t(4\ell+\eta_t)}. \label{eq:loss-spatial}
\end{align}
At a regular time when the entire hull misses the footprint, the
initial $2$ in \eqref{eq:loss-spatial} can be dropped. To see the
source of these estimates, the rate of $\int_I H\,dx$ is $q|I|$
plus the free $q(1-\nu_y)$ integral minus the signed active $j$
integral on $I$. Subtracting the equal-width identities cancels the
first term. The error upon replacing active by positive flux on the
two disjoint bands is at most $m_t$ in total. Integration gives
\eqref{eq:loss-storage}. The spatial derivative of $H$ is the
projection of $-\nu_x\,ds$ on the free boundary plus that of $n_x\,ds$
on the active disk. The free part on the hull has variation at most
$\sqrt{\eta_t(4\ell+\eta_t)}$ by \eqref{eq:reach-space-new}.
The difference of band means weights the active contribution by a
function in $[0,1]$ up to an overall sign, and each sign of $n_x$
has total magnitude two on a full circle. This proves
\eqref{eq:loss-spatial}, including the improvement off the disk.
For each band pattern below one can fix a countable dense collection
of generic translates, including the reflected patterns. The time
identities then hold for their absolutely continuous representatives,
with the spatial bounds simultaneous at a common full-measure set of
times. All band placements below have slack for such a translate.

\subsection{The cap at an interior turn}
\begin{figure}[ht]
\centering
\includegraphics[width=.76\linewidth]{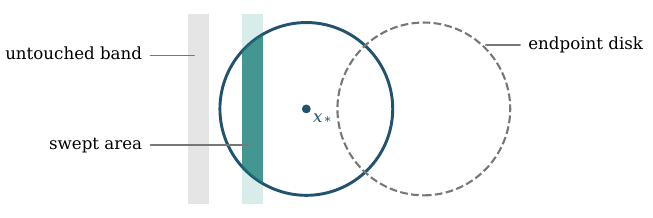}
\caption{Both endpoint disks miss the two bands. At the intervening
minimum a piece of disk enters the inner band while the outer band
stays untouched. Vertical translation of either disk leaves the band
areas unchanged.}\label{fig:turn-bands}
\end{figure}

Call a time good if it is regular for these identities and
$\eta_t\le a_0=1/4096$. Suppose a window stays in
$3/2<X<W-3/2$ and $\int b\le1/16$ on the window. For any good
$u<v$ therein, the entire intervening range of $X$ lies in
\begin{equation}\label{eq:loss-between}
 [\min(X(u),X(v))-1,\ \max(X(u),X(v))+1].
\end{equation}
For the lower end, if the minimum $x_*=\min_{[u,v]}X$ is at least
one below both endpoint values, choose $|p-x_*|<1/64$, and use
\[
 I_0=(p-11/8,p-9/8),\qquad I_1=(p-3/4,p-1/2).
\]
Both bands have width $1/4$ and their hull has length $7/8$ (Figure~\ref{fig:turn-bands}). The
outer band $I_0$ is untouched throughout. At $u$ and $v$ the whole
hull misses the disk, so $|S|<1/16$ by \eqref{eq:loss-spatial}.
Equation \eqref{eq:loss-storage} gives integrated flux into $I_1$
at most $3/32$. Yet the disk area in $I_1$ starts at zero and, at
the minimum, is at least
\[
             \tfrac14\,2\sqrt{1-(49/64)^2}>\tfrac14.
\]
This area increase is bounded by integrated positive circle flux.
Indeed, translating the disk with center $(X(t),Y(t)-qt)$ leaves
the area in a vertical band identical to that of the physical disk,
and its signed area derivative there is the circle integral of $j$.
This is a contradiction. Reflection proves the upper end of
\eqref{eq:loss-between}.

If the window endpoints $s_0,s_M$ are good, orient the $x$-axis so
that $X(s_0)\le X(s_M)$. For every pair of good times $u<v$ one has
\begin{equation}\label{eq:loss-order}
                        X(u)-X(v)\le2 .
\end{equation}
Otherwise, if $X(s_0)\le X(v)+1$, the overshoot at $u$ contradicts
\eqref{eq:loss-between} on $[s_0,v]$. If $X(s_0)>X(v)+1$, the
undershoot at $v$ contradicts it on $[u,s_M]$. Endpoint coincidences
obey the same inequalities.

\paragraph{Transverse progress on a short interval.}
Consider a length-16 interval $[u,v]$ in the same window whose
endpoints are good. Suppose $\int_u^v2R\,dt\ge30$. Then
\[
                         |X(v)-X(u)|>5/2 .
\]
Otherwise \eqref{eq:loss-between} confines the centers to an interval
of length at most $9/2$ and the full disk footprint to length at most
$13/2$. Enclose this footprint in an open band $I_1$ of width $27/4$
in the corridor. Place a disjoint band $I_0$ of the same width at a
gap of $1/16$ on the side with more room. Both choices are possible
when $W>30$, and their hull has length less than $14$. Since
$\sqrt{a_0(56+a_0)}<1/8$, \eqref{eq:loss-storage}--\eqref{eq:loss-spatial}
give
\[
 \int_u^v2R\,dt=\int_u^v\mu_{I_1}\,dt
      \le 2(27/4)(2+1/8)+1/16=115/4<30 ,
\]
a contradiction. By \eqref{eq:loss-order}, the progress is positive
in the common orientation of the window.

\paragraph{A long window and the finite constant.}
Suppose $0<q\le q_0=2^{-24}$ and $2/q-W\le k=2^{-21}$. Then $W>30$,
$1\le Q\le3/2$, $0\le K\le k+q_0$, and $\delta=qK$. Here we used
$K=2/q-W+2q/(Q+1)$. Set
\[
       M=\lceil2/q\rceil,\qquad L=16M,\qquad s_j=s+16j
                \quad(0\le j\le M)
\]
with $s$ uniform modulo $T$. The endpoint failure probability is at
most $8192K$ by \eqref{eq:loss-budgets}. The probability of window
budget $\int_s^{s+L}b>1/16$ is at most $16qLK$. Equation
\eqref{eq:loss-wall} bounds contact with a wall layer.

Call a grid interval deficient when either endpoint is not good or
its integrated $2R$ is below $30$. In the latter case its mean
$R$ is below $15/16$, at distance at least $1/16$ from $Q$. Jensen's
inequality on the interval, followed by averaging over $s$ in
\eqref{eq:loss-budgets}, bounds that probability by $256Q\delta$.
The expected deficient fraction is therefore at most
$8192K+256Q\delta$. Markov's inequality bounds the chance of more
than $1/32$ of the grid being deficient by $32(8192K+256Q\delta)$.
Summing the four failure probabilities gives at most
\[
 (6q+2\delta)/Q+2^{15}qLK+8192K+16qLK
                      +32(8192K+256Q\delta)<3/4 .
\]
For the last strict inequality substitute $K\le k+q_0$,
$\delta\le q_0(k+q_0)$ and $qL\le32+16q_0$. These substitutions give
a rational upper bound below $0.708$. This argument uses only first moments
and requires no independence.

Choose a phase avoiding all four failures and orient its window as in
\eqref{eq:loss-order}. At least $31/32$ of the intervals make progress
greater than $5/2$. Each remaining interval loses at most $16$ by the
guard speed bound. Total transverse displacement is consequently at least
\[
      M\left(\frac{31}{32}\frac52-\frac{16}{32}\right)
                 =\frac{123}{64}M\ge\frac{123}{32q}>\frac{2Q}{q}.
\]
This exceeds the entire corridor width by the area ceiling. It proves
$2/q-W>k$ for periodic polygonal guards in the stated speed range.
At each fixed $q$, Theorem~\ref{thm:periodic-approximation} supplies a
periodic polygonal guard on widths arbitrarily close from below to
each width guarded by a bounded-excursion schedule. Passing to the
supremum and restoring units proves Theorem~\ref{thm:loss}.

\section{The unbounded additive value of direction information}
\label{app:direction-loss}

This appendix proves the qualitative divergence $D_2(q)=2/q-W_2(q)/r\to+\infty$
as $q\downarrow 0$ in Theorem~\ref{thm:direction-value}, together with the
fixed-cap obstruction $\limsup_{q\downarrow 0}qW_2(q,hr)/r<2$ in
Corollary~\ref{cor:two-way-height}; Appendix~\ref{app:direction-rate} upgrades
the divergence to the sharp $\Theta(q^{-1/2})$ rate. Work in normalized units
$r=V=1$ on $[0,W]\times\R$. If $D_2(q)$ did not tend to $+\infty$,
Theorem~\ref{thm:periodic-approximation} would supply a sequence of successful
periodic polygonal two-way guards with $q\downarrow0$ and
\[
 W=2/q+O(1),\qquad 0\le K=2\sqrt{1+q^2}/q-W\le K_0
\]
for some fixed $K_0$. All limits below refer to this sequence. Draw a starting
phase uniformly from the guard's own period and denote expectation over this
phase by $\mathbb E$; neither the period nor the longitudinal excursion is
assumed uniformly bounded along the sequence.

\subsection{Comparing the two sets of reachable points}
At physical time $t$ write $E_b(t)$ and $E_a(t)$ for the spatial sets
reachable from below and above in the past. Use Appendix~\ref{app:reach}
in the physical and height-reflected coordinates, respectively.
Thus $H_b(x,t)$ is the renormalized measure of the lower set,
and $H_a(x,t)$ that of the upper set in reflected coordinates. Their
equal-measure rays in physical height are $y<H_b$ and $y>-H_a$.
Write $G=H_b+H_a$. If both slices are rays, positive $G$ is the
length of their overlap and negative $G$ measures the intervening gap.

\begin{lemma}[Approximation by vertical rays]\label{lem:ray-approx}
At a regular time of either reachable set with $\eta_t\le M$, if
$W>2M+20$, the area of its symmetric difference from its equal-measure
rays is at most $C(M+1)^2$ for an absolute $C$. On a transverse
interval of length $d$ the essential oscillation of its height is at
most $\sqrt{M(4d+M)}+4$.
\end{lemma}
\begin{proof}
Use the lower orientation. Almost every vertical section ends in zero
at the top and one at the bottom. If it is not a single ray up to null
sets, its indicator has at least one upward jump. By vertical slicing,
the total measure of these exceptional lanes is at most the integral
of the negative vertical outward normal over the spatial boundary.
The free part contributes at most $\eta_t$ and the circle at most two.
Every typical exceptional lane therefore has a typical single-ray lane
within distance $M+4$, also at either wall. For two typical lanes the
$L^1(dy)$ distance of their indicators is at most the integral of
the absolute horizontal normal across the intervening strip. This
follows by horizontal slicing and the one-dimensional variation
inequality. By \eqref{eq:reach-space-new} this integral is at most
$\sqrt{M(4d+M)}+4$. Comparing with the nearby ray and then moving
its endpoint to the equal-measure height costs at most twice that
amount in $L^1(dy)$. Integrating on the exceptional lanes proves the
area bound. The same horizontal variation estimate applies to the
renormalized heights themselves. Reflection gives the upper case.
\end{proof}

For a physical window $[u,v]$ of length $L$ let $S$ be the spatial
union of the closed sensing disks during the window, clipped to the
corridor. Its transverse hull $J_0$ has length at most $L+2$, and
$|S|\le\pi+2L$. The latter follows by adding one cylindrical strip
of area at most twice the length of each polygonal leg to the initial
disk. Waiting outside $S$ preserves reachability, hence both removed
sets $E_i(u)\setminus E_i(v)$ lie in $S$, for $i=b,a$. Integrating
\eqref{eq:reach-time-new} on $J_0$, with the full active trace in this
band, gives
\begin{align}
 0\le |S\setminus(E_i(u)\setminus E_i(v))|
 &\le \pi+2\int_u^v |1-R_i(t)|\,dt
       +q(L+2)L+\int_u^v b_i(t)\,dt .                 \label{eq:shared-sweep}
\end{align}
Here $R_b=|g'-qe_y|$, $R_a=|g'+qe_y|$, and $b_i=q\eta_i+m_i$ is
the nonnegative error from Appendix~\ref{app:reach}. In particular
$\mathbb E b_i\le qK_0$ and $\mathbb E|1-R_i|\to0$. The interval identity
can be used with the path-dependent $J_0$: at typical endpoint times
it holds as an equality of measures in $x$, by first testing a
countable dense set of continuous functions.

\begin{lemma}[Tightness of the height difference]\label{lem:gap-tight}
At a uniform independent lane and phase, the sequence of laws of $G$
is tight.
\end{lemma}
\begin{proof}
First consider the lanes inside the disk footprint at a uniform phase
$u$. Given $\epsilon>0$, choose $M$ so that both $\eta_i$ at $u$ and
$u+L$ are at most $M$ with probability at least $1-\epsilon$,
uniformly in the deterministic $L$. This uses $\mathbb E\eta_i\le K_0$.
Choose a fixed $L$ much larger than $C(M+1)^2+10$, where $C$ is as
in Lemma~\ref{lem:ray-approx}. For this fixed $L$, the mean bounds
on $b_i$ and $|1-R_i|$ and Markov's inequality give each right side
of \eqref{eq:shared-sweep} at most $\pi+L/16$, with additional failure
probability at most $\epsilon$ for sufficiently small $q$.
The integral height balance gives
\[
 |E_i(u)\setminus E_i(u+L)|\ge
 2L-2\int_u^{u+L}|1-R_i|\,dt-q(L+2)L-\int_u^{u+L}b_i\,dt.
\]
Both removed sets thus have area at least $2L-L/16$. Their common
part has area at least $2L-L/8-\pi$ since both lie in $S$.

Replace the four endpoint sets by their ray families. The sum of their
area errors is $A\le4C(M+1)^2$. The two ray decrements therefore
have common area in $S$ greater than $A$, while the sum of their
areas outside $S$ is at most $A$. To see what this implies, enclose
$S$ between horizontal levels $l,h$ of separation $d\le L+2$.
On a lane with a common decrement in $S$, both initial rays contain
some point of $[l,h]$. If $H_b(u)>h+d$, the lower ray decrement
extends through an entire length $d$ above $h$. If $-H_a(u)<l-d$,
the upper ray decrement extends the same distance below $l$. In
either case the decrement outside $S$ pays for the whole common
decrement inside that lane. Since the total outside area is too small,
there is a typical lane in $J_0$ where $0\le G(u)\le3d$. The spatial
bound of Lemma~\ref{lem:ray-approx}, applied twice on $J_0$, now gives
$|G(x,u)|\le B_\epsilon$ throughout the footprint up to null sets,
with $B_\epsilon$ independent of $q$, except on phases of probability
at most $2\epsilon$.

To pass from the footprint to a uniform lane, use the periodic
one-dimensional height balance at almost every $x$:
\begin{equation}
 D_tG=(2q-\mu_b(x,t)-\mu_a(x,t))\,dt+de_x,\qquad e_x\ge0 .
                                                        \label{eq:G-evolution}
\end{equation}
The densities $\mu_i$ are full positive circle flux. For a $C^1$
nondecreasing function whose derivative lies in $[0,1]$, vanishes
on $[-B_\epsilon,B_\epsilon]$ and equals one for $|z|\ge B_\epsilon+1$,
integrate the BV chain rule for its composition with $G$ over a
period. The singular increments in \eqref{eq:G-evolution} are
nonnegative, including their contribution at jumps. Consequently
\[
 2q\,\mathbb E\int_0^W {\bf1}_{\{|G|>B_\epsilon+1\}}\,dx
 \le \mathbb E\int_0^W {\bf1}_{\{|G|>B_\epsilon\}}(\mu_b+\mu_a)\,dx
 \le 8(1+q)\epsilon.
\]
The final estimate uses support in the footprint and total flux
at most $4(1+q)$. Since $qW\to2$, this proves tightness.
\end{proof}

\subsection{A stationary limit at a randomly chosen lane}
Choose $s$ uniformly on $(0,2)$, independently of the starting phase
$u$, and sample the physical lane $x=sW/2$. For macroscopic time
$\tau$ write $t(\tau)=u+\tau/q$. Set
\[
 Z_q(\tau)=2X(t(\tau))/W,\qquad G_q(\tau)=G(x,t(\tau)),
 \qquad d\lambda_q=\frac{\mu_b(x,t(\tau))+\mu_a(x,t(\tau))}{2q}\,d\tau .
\]
The joint process is stationary in $\tau$. The paths $Z_q$ have
Lipschitz constants $2/(qW)\to1$. On any fixed time interval of
length $l$, $\mathbb E\lambda_q\le 2(1+q)l/(qW)$. On that interval
each height increment is $l-\lambda_q$ up to an error whose expected
absolute value tends to zero. Indeed the error measure for either
height has expected mass at most $lK_0/W$. The difference of the two
circle fluxes has total variation at most $8q$ per unit physical time,
using $|a_+-b_+|\le|a-b|$ and $\int|n_y|\,ds=4$.

By compactness of equi-Lipschitz paths, tightness of locally finite
measures under the displayed local mass bounds, and
Lemma~\ref{lem:gap-tight}, pass to a subsequential joint weak limit
$(s,Z,\lambda,G(0))$. Here paths converge locally uniformly and
measures vaguely. One may use the Skorokhod representation on this
product of separable complete metric spaces \cite{Billingsley}.
The limiting measure has no atom at any fixed deterministic time,
by the local mean mass bound on open intervals. Define a right
continuous process by its increments
\begin{equation}
             G(v)-G(u)=2(v-u)-2\lambda((u,v]) .
                                                        \label{eq:macro-G}
\end{equation}
Here and below $u,v$ in limit-process formulas denote macroscopic
times. Convergence of the finite-dimensional laws at deterministic
times follows from the vanishing increment errors. In particular
$(s,Z,\lambda,G)$ is stationary. The lane $s$ remains independent of
$Z$, and $Z(0)$ is uniform on $[0,2]$ by the localized flux estimate
in Appendix~\ref{app:reach}.

The support of $\lambda$ lies in $\{t:Z(t)=s\}$ because the original
flux needs $|Z_q-s|\le2/W$. For a Lipschitz real path, almost every
level has finitely many preimages in each compact interval, all
differentiability points with nonzero derivative. These statements
follow from the one-dimensional area formula, the null image of
the nondifferentiability set, and the area formula on the
zero-derivative set. They apply almost surely to the independent $s$.
Thus the level hits are isolated, with a change of side at each hit.

Each hit carries an atom of $\lambda$ of mass at least two. In fact,
bracket a transversal hit by rational times whose endpoint values of
$Z$ are on opposite sides of $s$. Eventually the physical centers
cross the whole lane footprint in this bracket. Pairing the two
circle points with $n_x=z=x-X$ gives the lane flux density for either
direction, for $|z|<1$,
\begin{equation}
 \max\left\{|X'|\frac{|z|}{\sqrt{1-z^2}},\,|Y'\mp q|\right\}
                         +X'\frac{z}{\sqrt{1-z^2}} .      \label{eq:both-kernel}
\end{equation}
The signed term integrates to zero between endpoints on opposite
sides outside the footprint, being a potential derivative. The
first term integrates to at least
$\int_{-1}^1 |z|/\sqrt{1-z^2}\,dz=2$ by one-dimensional coarea.
These facts also follow by truncating the integrable endpoint
singularities. Passing to the limit on the rational bracket and
shrinking it proves the claim. Write $a_i=\lambda\{t_i\}$ at the
successive hits.

\subsection{Every crossing divides its height change equally}
We show that the values at a hit satisfy
\begin{equation}
                   G(t_i-)=a_i,\qquad G(t_i)=-a_i .
                                                        \label{eq:macro-swap}
\end{equation}
Fix $A>0$ and partition $[-A,A]$ into dyadic intervals of common
length $l$ (successively halved). In an interval with endpoints
$u,v$, the physical duration is $L=l/q$. Sum
\eqref{eq:shared-sweep} over the mesh, for both orientations, and
divide by $W$. The expectation of the sum of the normalized areas
of the symmetric differences of the two actual decrements is at most
$C_A l+o(1)$ as $q\to0$ with $l$ fixed. Indeed the sum of
$qL^2/W$ is $2Al/(qW)$, the constants per interval disappear,
and the mean speed and budget errors tend to zero after normalization.

We can replace the endpoint sets by rays for this estimate using
the uniform random lane. Take $M=W^{1/4}$, so the probability of any
endpoint failing $\eta_i\le M$ tends to zero on a fixed mesh. On the
good event the normalized sum of ray-area errors is $O_{A,l}(M^2/W)$.
On the exceptional event truncate each lane symmetric-difference
length at one, which bounds its contribution to the expectation by
the number of mesh intervals times the exceptional probability.
Consequently the expected sum of these truncated lane discrepancies
for the ray decrements is still bounded by $C_A l+o(1)$.

On a lane translate the initial bottom of the upper ray to zero.
The two ray-decrement intervals then have endpoints
$[G_q(u)-d_b,G_q(u)]$ and $[0,d_a]$, where
$d_i=(H_i(u)-H_i(v))_+$ in macroscopic time. Empty decrements
have length zero. Both $d_i$ converge in joint law on the fixed mesh
to $d=(\lambda((u,v])-l)_+$. The symmetric difference of two
equal-length intervals of length $d$ equals twice the smaller of
$d$ and the distance between their left endpoints. By bounded
continuity on this mesh, the limiting estimate is
\[
  \mathbb E\sum_{[u,v]}
   \min\{1,\,2d,\,2|G(u)-d|\}\ \le C_A l .
\]
Almost surely no dyadic endpoint is a hit, and there are finitely
many hits on the compact interval. As the mesh size decreases, the
cell containing an interior hit eventually contains only that hit.
There $d\to a_i$ and $G(u)\to G(t_i-)$. Fatou's lemma gives
$\min\{1,2a_i,2|G(t_i-)-a_i|\}=0$ at every such hit. Since
$a_i\ge2$, this proves the first equality in \eqref{eq:macro-swap}.
The second follows from \eqref{eq:macro-G}.

\subsection{The incompatible occupation times}
Almost surely hits occur arbitrarily far in both time directions.
On the shift-invariant event of no hits, $G(t)=G(0)+2t$, which
contradicts stationarity of a finite random variable conditional on
that event if it has positive probability. A stationary locally
finite point set cannot have a first or a last point with positive
probability, since its location would have a translation-invariant
probability distribution on $\mathbb R$.

Between consecutive hits \eqref{eq:macro-G} and \eqref{eq:macro-swap}
give $2(t_{i+1}-t_i)=a_i+a_{i+1}$. The intervals
\[
                        [t_i-a_i/2,\ t_i+a_i/2]
\]
therefore tile the time line with disjoint interiors. Exactly half
of each such interval has $Z>s$ and half has $Z<s$, since a hit
changes side and there are no intervening hits. For clarity, average
these indicators over $[0,B]$. Only the tiles containing the two
endpoints can contribute an imbalance. After division by $B$, their
expected contribution tends to zero: it is bounded by twice
$\mathbb E\min(1,L_0/B)$ where $L_0$ is the almost surely finite
length of the tile at a uniform time, by stationarity. Endpoint ties
can be assigned to the tile to their right. The same estimate holds
after multiplication by any indicator depending only on $s$.
Consequently $\Pr(Z(0)>s\mid s)=1/2$ for almost every $s$.
Independence and the uniform marginal instead give $1-s/2$. This
contradiction excludes every sequence of bounded two-way deficits
and proves the divergence in Theorem~\ref{thm:direction-value}.

\subsection{Keeping the guard in a fixed longitudinal strip}\label{app:fixed-two-direction}
We prove the necessity in Corollary~\ref{cor:two-way-height}. Work in
unit-radius coordinates. If a sequence with height bounded by
$h<\infty$ attained $qW\to2$, Theorem~\ref{thm:periodic-approximation}
would give a periodic polygonal sequence with the same cap and limit.
Now $\delta=2\sqrt{1+q^2}-qW\to0$ and $K/W=\delta/(qW)\to0$.
The height sum $G$ is bounded deterministically by $2(h+2)$ since
each reachable section agrees with its oriented half-strip outside
$[-h-2,h+2]$. At a regular time the area discrepancy for either
ray family is at most $2(h+2)(\eta_t+2)$. Indeed the exceptional
lanes have measure at most $\eta_t+2$ by vertical slicing as in
Lemma~\ref{lem:ray-approx}, and both indicators agree outside that
fixed vertical interval. Dividing by $W$ and averaging gives a bound
tending to zero, even if $K$ grows.

The stationary limit and mesh comparison now apply with these two
estimates. In detail, on each fixed macroscopic window of length $l$
the expected lane increment error is at most $lK/W+O(l/W)$.
In the sweep comparison the mean physical budget is at most
$\delta\to0$, and $R_i\to1$ in mean by Appendix~\ref{app:reach}.
At the finitely many endpoints of each fixed mesh, the expected ray
discrepancies divided by $W$ tend to zero by the fixed-strip bound
above. The limits \eqref{eq:macro-G} and \eqref{eq:macro-swap}, and
the uniform transverse marginal, follow exactly as before. They force
the incompatible conditional probabilities $1/2$ and $1-s/2$.
No such sequence exists, which gives the strict limsup for each finite
cap by the area ceiling.

\section{The rate of the directional advantage}\label{app:direction-rate}
This appendix proves the quantitative bounds in Theorem~\ref{thm:direction-value}:
we first establish a single-window version of the two-orientation occupation
argument of Appendix~\ref{app:direction-loss} to prove $D_2(q)\ge c_0/\sqrt q$
for small $q$, and then construct the reclining-stroke scan ($4\sqrt 2$) and its
diagonal-return refinement ($2\sqrt 6$) in Subsections~I.1--I.2.

Work first in normalized units $r=V=1$ on $[0,W]\times\R$ with $0<q\le 1/2$.
The estimates of Appendices~\ref{app:reach} and~\ref{app:direction-loss} used
below hold at a uniform phase of any successful periodic polygonal two-way guard,
with no advance bound on its period or longitudinal excursion.
Let $W$ be its width, $Q=\sqrt{1+q^2}$,
$\delta=2Q-qW\ge 0$ (by~\eqref{eq:area-bound}), $K=\delta/q\ge 0$, and $k=1+K\ge 1$.
We may restrict to $1\le qW\le3$ by~\eqref{eq:area-bound} and the immediate
$1/q$ deficit when $qW<1$. All implicit constants in this section are absolute.

Write $E_b,E_a$ for the spatial sets reachable in the past from
below and above and $H_b,H_a$ for their renormalized heights in
their respective orientations. The equal-measure rays in physical
height are $y<H_b$ and $y>-H_a$. Put $G=H_b+H_a$.
The spatial excess $\eta_i$ has mean at most $K$. On each lane
the time balance is $D_tH_i=q\,dt-\mu_i\,dt+de_i$ with nonnegative
error. The mean total error, integrated over lanes per unit
physical time, is at most $\delta$. Here $\mu_i$ is the full positive
circle flux density. Its difference between orientations has total
variation at most $8q$ per physical time. Writing
$R_{b,a}=|g'\mp q e_y|$ gives $\mathbb E|1-R_i|\lesssim\sqrt\delta+q$.
A typical vertical section is an exact ray except on lanes of total
measure at most $\eta_i+2$.

\begin{lemma}[Swept height on one lane]\label{lem:lane-sweep}
On a physical time window of duration $L$ put $l=qL$ and
$a_i(x)=\int\mu_i(x,t)\,dt$, $a=(a_b+a_a)/2$. At almost every lane
whose sensing slice is empty at the starting time, the length of
its physical swept set is at most $a+l$. Each $|a_i-a|\le l$.
If the center crosses the whole footprint then $a_i\ge2$ for both
orientations.
\end{lemma}
\begin{proof}
Set $z=x-X(t)$ and $R=\sqrt{1-z^2}$ when $|z|<1$. The moving
interval on this lane has endpoints $Y-R,Y+R$. New swept length
is at most the total positive variation of the upper endpoint
plus the negative variation of the lower endpoint. An interval
which appears at a later footprint entry starts with zero length.
Summing over the components gives the integral of
\[
  (Y'+X'z/R)_++(-Y'+X'z/R)_+
     =\max\{|Y'|,|X'z/R|\}+X'z/R .
\]
This is the paired-circle density at $q=0$. The densities for the
two orientations replace $Y'$ by $Y'\mp q$ inside the maximum,
changing it by at most $q$. Integrating proves the first two claims.
The integrable endpoint singularities can be handled on each affine
piece by truncation. The lanes exactly tangent to constant-transverse
pieces form a null set. During a complete passage the signed term
integrates to zero between the exterior endpoints, as it is the
derivative of the half-length $R$ extended by zero. The remaining
term integrates to at least $\int_{-1}^1 |z|/\sqrt{1-z^2}\,dz=2$
by one-dimensional coarea. All densities are nonnegative outside
that crossing interval as well.
\end{proof}

\begin{lemma}[Linear scale bound on $G$]\label{lem:bounded-end}
For each $p>0$ there are constants $A>1$ and $\eta>0$ such that
if $qk^2<\eta$, then at an independent uniform lane and phase
$\mathbb P(|G|>Ak)<p$. Also, except on phases of probability less
than $p$, $|G|\le Ak$ on almost every lane in the current disk
footprint.
\end{lemma}
\begin{proof}
Here are the quantitative choices in the physical-window proof
of Lemma~\ref{lem:gap-tight}. Fix a sufficiently small absolute
failure probability $p_0$. At both endpoints of a physical window
of length $L$, both spatial excesses are at most $M=4k/p_0$
except with probability $p_0$. With $qk^2$ small,
$W>2M+20$. The symmetric-difference area between each endpoint
set and its rays is at most $C(M+1)^2$ by Lemma~\ref{lem:ray-approx}.
Take $L=C'(M+1)^2$ with $C'$ sufficiently large. Write $u,v$ for
the endpoints. The physical sweep has area at most $2L+\pi$.
Each actual removed set lies in it and has area at least $2L-L/16$
except with additional probability $p_0$. The mean subtracted
errors divided by $L$ are at most
$C(\sqrt\delta+q+\delta)+q(L+2)$. At these good endpoints the
two ray decrements have a common part of positive area. Indeed
the actual decrements have common area at least $2L-L/8-\pi$,
and replacing their endpoints by rays costs at most $4C(M+1)^2$.
Outside the transverse sweep hull the decrements vanish since
waiting preserves reachability there.

The ray decrement lengths $d_i(x)=(H_i(x,u)-H_i(x,v))_+$ are
bounded by $C_1 k$ at almost every lane. Their integrals are at
most $2L+\pi$, since a net height decrease costs at least that much
actual removal on the lane. In an interval of length $M+1$ adjacent
to a typical $x$, select a lane with $d_i\le(2L+\pi)/(M+1)$.
The two endpoint oscillation bounds of Lemma~\ref{lem:ray-approx} on this interval
give $d_i(x)\le(2L+\pi)/(M+1)+2\sqrt{M(4M+4+M)}+8$.
Where the two physical ray decrements $[H_b(x,v),H_b(x,u)]$ and
$[-H_a(x,u),-H_a(x,v)]$ overlap, we have $G(x,v)\le 0\le G(x,u)$, yielding
$0\le G(x,u)=G(x,v)+d_b(x)+d_a(x)\le d_b(x)+d_a(x)\le2C_1k$.

We can find such an overlap within $C_2 k$ of the current footprint.
The expected total positive error in the lane time balances over
the window is at most $2L\delta=O_{p_0}(k)$ under $qk^2<\eta$.
By Markov's inequality this total is at most $C_3 k$ except
with probability $p_0$. Choose $\eta$ small enough that $qL<1/16$.
Choose a constant $D_0$ sufficiently large compared with $C_3$
and $1/p_0$. If the guard stays within $2D_0k+4$ of $X(u)$
transversely, the common ray lane already found is close enough.
Otherwise, on the side it reaches, an interval of length $D_0k$
within $2D_0k+4$ of $X(u)$ has every lane fully crossed from
one exterior side of its footprint to the other. On each such lane
each orientation's integrated circle flux is at least two by the
paired-circle kernel in~\eqref{eq:both-kernel}. Exclude the non-ray slices
at the four endpoints, of total lane measure at most $4(M+2)$.
Also exclude lanes where the summed time balance error exceeds
$1/4$, of measure at most $4C_3k$. Take the crossed interval
at least two units from $X(u)$, so it starts outside the disk.
There is a remaining lane. Its two decrement lengths sum to at
least $2a-2qL-1/4$ by the two exact height balances. Both removed
sets lie in the physical swept set, of length at most $a+qL$ by
Lemma~\ref{lem:lane-sweep}. Their common length is therefore at
least $a-3qL-1/4>0$.
The endpoint rays make this an overlap of the two ray decrements.
In both cases the spatial oscillation bound over distance $C_2k$
at time $u$ now gives $|G|\le A_0k$ throughout the current footprint
up to null sets.

Finally apply the periodic lane chain rule to a monotone truncation
whose derivative vanishes on $[-A_0 k,A_0 k]$ and equals one
outside its unit enlargement. Positive singular error increments
can be discarded. As in Lemma~\ref{lem:gap-tight} this bounds
$2qW\mathbb P(|G|>A_0 k+1)$ by the mean footprint flux on bad
phases, at most $12(1+q)p_0$. Since $qW\ge1$, choosing $p_0$
small enough in terms of $p$ proves both claims.
\end{proof}

\paragraph{An occupation identity on a single window.}
Take the starting phase uniform in the guard period, and independently
take the physical lane $x$ uniform in $[0,W]$. Use a deterministic
window of physical duration $L=l/q$, with an absolute $0<l<1/16$
chosen below. Use subscripts $-,+$ for its endpoints, $\Delta$ for
the latter value minus the former and a bar for their arithmetic
mean. Set $w(s)=1-s$ on $[0,1]$, extended constantly on either side,
and $\psi(v)=\max(-1,\min(1,v/b))$ for a fixed smoothing width $b>1$
in physical units. At the sampled lane put
$S(t)=\psi(x-X(t))$ and use $w=w(x/W)$ as a weight.

Choose $p$ in Lemma~\ref{lem:bounded-end} and $M=Ak$ as there. Let
$T$ be an odd nondecreasing $C^1$ truncation with $0\le T'\le1$,
identical to its argument on $[-M,M]$ and constant beyond its
unit enlargement. Stationarity and the product identity give
\[
      0=\mathbb E\,w\Delta(ST(G))
       =\mathbb E\,w\overline S\,\Delta T(G)+\mathbb E\,w\Delta S\,\overline {T(G)} .
\]
The lane time balance is $D_tG=2q\,dt-(\mu_b+\mu_a)\,dt+de$,
with a nonnegative error measure. By the BV chain rule and
Lemma~\ref{lem:bounded-end},
\[
              \mathbb E|\Delta T(G)-2l+2a|\le C l(p+qk).
\]
Indeed the drift discarded by truncation costs at most $2pl$ at a
uniform lane. The discarded flux vanishes on the good footprint
phases, and has full spatial mass at most $4(1+q)$ per physical
time otherwise. Its expected normalized cost is at most
$4(1+q)pL/W$. The error measure contributes at most its full
mass, also at jumps, whose expected lane value is at most
$2L\delta/W$. Here and below $1\le qW\le3$.

\paragraph{Cancellation of the flux-weighted side.}
We claim
\[
 |\mathbb E\,w\overline S a|
 \le C l\bigl(l+b^{-1}+W^{-1}+\sqrt\delta+q+\delta\bigr).
 \tag{*}
\]
To see this, integrate over the random lane first. At time $t$
the circle flux is supported within physical distance one of $X(t)$.
Replace $w(x/W)\overline S(x)$ by
\[
      w(X(t)/W)\,
      \frac{\psi(X(t)-X_-)+\psi(X(t)-X_+)}2
\]
at cost $O(W^{-1}+b^{-1})$ per unit mass. The full spatial mass of
the mean of the two circle fluxes is $R_b+R_a$. Mass outside the
corridor has mean at most $\delta$ per orientation, since the
active trace is confined to the corridor and the missing flux is
part of the nonnegative budget of Appendix~\ref{app:reach}. Thus we may replace
the spatial mass by two with mean error $O(\sqrt\delta+q+\delta)$.
All these errors have the factor $L/W\le l$ after lane averaging.

In the remaining expression use stationarity and the oddness of
$\psi$. For a uniform starting phase, writing $X(t)$ relative to
that phase, the unweighted-in-time average of the preceding display
is exactly
\[
   \frac1{2L}\int_0^L \mathbb E\bigl[
     (w(X(t)/W)-w(X(0)/W))\psi(X(t)-X(0))\bigr]\,dt .
\]
The term involving the ending phase becomes the second term of
this difference by using the lag $L-t$ and shifting the starting
phase. The display is bounded by $L/(2W)\le l/2$; multiplying by the window mass factor $2L/W\le 2l$ proves $(*)$.
This cancellation includes time windows with arbitrary transverse
reversals and windows spent entirely near a single lane.

\paragraph{The change of side at a crossing.}
The second product term satisfies
\[
    |\mathbb E\,w\Delta S\,\overline{T(G)}|
    \le C l^2+C lqk+C(M+1)\bigl((b+1)/W+qk+lqk\bigr).
 \tag{**}
\]
First remove lanes within $b+1$ of either endpoint center, of
probability at most $4(b+1)/W$. Outside these strips a change in
$S$ requires the two endpoints to straddle the full lane footprint.
Such lanes have probability at most $L/W\le l$. Endpoint slices
of both orientations are exact rays except with probability
$4(K+2)/W$. Write the separate balances as
$\Delta H_i=l-a+e_i$. The circle-flux difference and the positive
time error give $\mathbb E e_*\le C lqk$, where $e_*=|e_b|+|e_a|$.
Indeed the time error integrates to at most $2L\delta/W$ in
expectation and the orientation difference to at most $8qL/W$.
Remove also the event $e_*>1/4$. At every removed lane the
absolute integrand is at most $2(M+1)$.

On a remaining straddled lane put $d=a-l$. By
Lemma~\ref{lem:lane-sweep}, $d\ge2-l$ and the actual physical
swept set has length at most $a+l$. The two decrement intervals
of the exact rays have lengths $d-e_i>0$ and both lie in this
set, since a point outside the sweep may wait safely. Their
symmetric difference therefore has length at most $4l+e_*$.
Translate the initial boundary of the upper ray to zero. Replacing
the decrement lengths by $d$ gives the two intervals
$[G_--d,G_-]$ and $[0,d]$, at additional cost at most $e_*$.
Their symmetric difference is $2\min(d,|G_--d|)$. Since
$4l+2e_*<2d$, we obtain
\[
 |\overline G|\le2l+\tfrac32 e_*,
 \qquad |\overline{T(G)}|\le|\overline G|,
\]
since $T(G_-)+T(G_+)=T(G_-)-T(-G_+)$ and $T$ is odd and $1$-Lipschitz. Integrating over the
straddled lanes gives the first two terms of $(**)$.

\paragraph{The transverse marginal.}
At a stationary phase $X/W$ is within $O(\sqrt\delta+q+\delta)$
of uniform measure on $[0,1]$ when tested against a fixed bounded
Lipschitz function. Indeed the period mean circle flux dominates
$q\,dx$ on the corridor, with total surplus at most $\delta$ on
the full line, by Appendix~\ref{app:height-balances}. Replacing the test point
$(X+n_x)/W$ by $X/W$ costs $O(1/W)$, and the full mass $2R_i$
has mean absolute deviation $O(\sqrt\delta+q)$ from two.
Also $qW=2\sqrt{1+q^2}-\delta$. The independent lane thus gives
\[
  \mathbb E\,w\overline S=-1/6+O(\sqrt\delta+q+\delta+b/W).
\]
For an exactly uniform center and a hard sign, the integral is
$\int_0^1(1-s)(2s-1)\,ds=-1/6$. Smoothing changes an interval
of lane length at most $2b/W$ in normalized coordinates.

Suppose $qk^2$ is sufficiently small. The product identity bounds
$2l|\mathbb E w\overline S|$ by $2|\mathbb E w\overline S a|$ plus the
right side of $(**)$ and $C l(p+qk)$. After division by $l$,
choose $p$ small, $b$ large and $l$ small to make the terms
in $p,b^{-1},l$ arbitrarily small. These terms have absolute
constants independent of $M$. Now $M=Ak$ with $A$ fixed by $p$,
and the remaining terms tend uniformly to zero as $qk^2\to0$,
including $(M+1)((b+1)/W+qk+lqk)/l$. The nonzero transverse
marginal contradicts this inequality. Consequently an absolute
$c_0>0$ satisfies $k\ge c_0/\sqrt q$ for all sufficiently small
$q$ and every successful polygon.

Finally $k=1+(2/q-W)+2q/(\sqrt{1+q^2}+1)$. Periodic polygonal
approximation at the same speed transfers the bound to the
two-way supremum. The one-way deficit $2/q-W_1/r$ is bounded
by Theorem~\ref{thm:slow-hooks}. Absorbing this bounded term proves
the lower estimate in Theorem~\ref{thm:direction-value}.

\subsection{A scan with reclining strokes}
In each
almost vertical stroke give the guard transverse velocity $-q$.
The succeeding rightward transfer shifts $X+qt$ by exactly two
capture radii. This makes the times at which adjacent strokes reach
the same height fit together across the full scan.

Set $v=\sqrt{1-q^2}$, choose a stroke half-height $h$, and take
\[
 t_0=2h/v,\qquad b_0=2/(1+q),\qquad s=b_0-qt_0.
\]
Use an even number $n$ of strokes, alternately up and down, with
velocity $(-q,\pm v)$. Stroke $j$, $0\le j<n$, starts at transverse
position $qt_0+js$ and ends at $js$. Transfer to the next start
at unit transverse speed over distance $b_0$. From the end of the
last stroke return along the bottom to the original start. The
guarded corridor will be $[0,(n-1)s+qt_0]$. For the parameters
below the final return has length $(n-1)s-qt_0>0$, and the period
$P$ is at most $n(t_0+b_0+s)$.

For each fixed $y\in[-h,h]$ let $u_j(y)\in[0,t_0]$ be its time
within stroke $j$. The time and transverse position of the corresponding
witness, in one pass, are
\[
 t_j(y)=j(t_0+b_0)+u_j(y),\qquad
 x_j(y)=qt_0+js-q u_j(y).
\]
Adjacent witnesses have $\Delta t>0$, $\Delta x\ge b_0-2qt_0>0$
and $\Delta x+q\Delta t=2$. Define $f(x,y)$ by linear interpolation
of the times in their ordered transverse coordinates, using the
outermost time constantly beyond each end. Every point $(x,y)$ at
time $f(x,y)$ inside the corridor is captured. Between two witnesses,
at interpolation fraction $z$, its distance plus time cost from the
first is $2z$ and from the second $2(1-z)$. One is at most one.
In an end extension the transverse distance at the witness time
is at most $qt_0<1$.

The graph $f$ is continuous in $x,y$. Set
\[
 J=\frac{2t_0+b_0}{b_0-2qt_0},\qquad B_0=t_0(1+qJ).
\]
It is $J$-Lipschitz in $x$, and $|f(x,y)-f(x,y')|\le B_0$ at
the same $x$. For the latter bound use the same two witness indices
and the same interpolation fraction at the new height. The
interpolated position changes by at most $qt_0$ and its time by
at most $t_0$. Finish by the transverse Lipschitz bound. The
constant end extensions satisfy the same estimate.

On every pass any admissible path meets the graph
$t=kP+f(x,\operatorname{clip}(y,[-h,h]))$, by continuity from the
start of the pass to its end. An uncaptured meeting must lie
outside $[-h,h]$ in physical height. For consecutive meetings the
elapsed time is at most $P+B_0+J|\Delta x|$. If $qJ<1$, the
Euclidean speed constraint therefore gives
\[
 |\Delta y|^2\le
  q^2(P+B_0+J|\Delta x|)^2-|\Delta x|^2
  \le \frac{q^2(P+B_0)^2}{1-q^2J^2}.
\]
Thus $q(P+B_0)\le2h\sqrt{1-q^2J^2}$ prevents a change between
the two exterior height ranges in either direction.

For fixed $c>0$ take $h=c/\sqrt q$ and
\[
 n=2\left\lfloor
 \frac{2h(1-q^2J^2)/q-B_0}{2(t_0+b_0+s)}
 \right\rfloor .
\]
All strict hypotheses hold for sufficiently small $q$ and this
rounding ensures the capture condition. Now
$s=2-2c\sqrt q+O_c(q)$ and $qn=1-2\sqrt q/c+O_c(q)$.
The resulting width is
\[
 q[(n-1)s+qt_0]=2-(2c+4/c)\sqrt q+O_c(q).
\]
Taking $c=\sqrt2$ gives an upper coefficient $4\sqrt2$ using
$W_1(q)/r\le2\sqrt{1+q^{-2}}$.
For example $q=10^{-4}$, $h=141$ and $n=9850$ guard a width between
$19418.316795571r$ and $19418.316795572r$. The exact period leaves
more than $0.1699$ in the inequality
$2h\sqrt{1-q^2J^2}-q(P+B_0)>0$.

\subsection{Using the return trip to block a transit}

We can improve the preceding coefficient by sweeping on a diagonal
during the return. Fix $0<D<1$ and put $a=D/2$. Use the same even
number $n$ of reclining strokes, stroke time $t_0$, transfer $b_0$,
step $s$ and drift $d_0=qt_0$ as above. Write
\[
 W=(n-1)s+d_0,\quad T_s=(n-1)(t_0+b_0)+t_0,\quad J_0=T_s/W .
\]
On the even pass the scan covers the band centered at $-a$, starting
and ending at its upper height $h-a$. On the odd pass it covers the
band centered at $+a$, starting and ending at its lower height
$-h+a$. At the end of the even scan connect straight to $(W,h+a)$,
return on the diagonal to $(0,-h-a)$ and connect straight to the
start of the odd scan at $(d_0,-h+a)$. After the odd scan reflect
the return heights. Each pass and its return have total duration
\[
 P=T_s+2c_0+\sqrt{W^2+(2h+D)^2},\qquad c_0=\sqrt{D^2+d_0^2}.
\]

\begin{lemma}[Diagonal-return two-way scan]\label{lem:diagonal-return}
Suppose $n\ge2$ is even, $b_0>2d_0$, $d_0\le1$ and $qJ_0<1$.
With $Q_0=\sqrt{1-q^2J_0^2}$, the guard captures every transit
in either direction if
\[
             q(P+2t_0)/Q_0\le2h+D,\qquad
             q(t_0+c_0)/Q_0\le1-D .
\]
\end{lemma}
\begin{proof}
Use the captured time graph $f$ of the preceding construction on
each main scan, with the appropriate center and starting orientation.
Throughout its band $|f(x,y)-J_0x|\le t_0$. In fact, at witness $j$
with local time $v\in[0,t_0]$ the difference is
\[
    j(t_0+b_0-J_0s)+(1+qJ_0)v-J_0d_0,\qquad
    (n-1)(t_0+b_0-J_0s)=J_0d_0-t_0 .
\]
Its four extremal values are bounded in magnitude by $t_0$ since
$J_0d_0<t_0$. Interpolation gives the same bound inside. At the
two constant end extensions also test $x=0,W$. There the local
time differs from $0$ or $T_s$ by at most $t_0$.

Choose one meeting with the clipped main graph on every scan, as
before, and one transverse meeting on each diagonal return, which
runs all the way from the right wall to the left wall. Put $u=x/W$
at that return meeting and let $R$ be the diagonal length. Its
time in the pass is $T_s+c_0+(1-u)R$. Subtracting $J_0x=uT_s$
gives $(1-u)P+(2u-1)c_0$. Thus the time from the preceding main
meeting to the return is at most $(1-u)P+c_0+t_0+J_0|\Delta x|$.
From the return to the next main meeting the bound is
$uP+c_0+t_0+J_0|\Delta x|$. Between successive main meetings the
bound is $P+2t_0+J_0|\Delta x|$. Whenever elapsed time is at most
$A+J_0|\Delta x|$, the Euclidean speed bound implies
$|\Delta y|\le qA/Q_0$ by maximizing the same quadratic as before.

On an even pass a switch from below the main band to above the next
one would need more than $2h+D$ vertical motion. For the opposite
switch, the return height at $u$ is $Y_R=(h+a)(2u-1)$. An uncaptured
return meeting is either below $Y_R-1$ or above $Y_R+1$. From above
the even band to below this observation the required motion is
greater than
\[
             h-a-Y_R+1=(1-u)(2h+D)+(1-D).
\]
From above this observation to below the odd band it is greater
than $Y_R+1+h-a=u(2h+D)+(1-D)$. Both possibilities are excluded
by the hypotheses, including at $u=0,1$. On an odd pass reflect
all heights and interchange the directions. Hence the exterior
side at every main meeting is constant along a surviving path.
Either complete transit would require a change of side.
\end{proof}

Take $h=c/\sqrt q$ and
\[
 n=2\left\lfloor
      \frac{2h+D-q^{1/4}}{2q(t_0+b_0+s)}
                \right\rfloor,\qquad c>0 .
\]
Here $n\asymp1/q$, $P=n(t_0+b_0+s)+O_c(1)$ because
$R-W\le(2h+D)^2/(2W)$, and $qJ_0=O_c(\sqrt q)$. The first
hypothesis on time in Lemma~\ref{lem:diagonal-return} has margin
$q^{1/4}-O_c(\sqrt q)$, and the left side of the second one tends
to zero. Thus every sufficiently small speed admits the guard.
Since $qn=1-(4-D)/(2h)+o(\sqrt q)$ and
$s=2-2c\sqrt q+O_c(q)$, its width satisfies
\[
       qW=2-\bigl(2c+(4-D)/c\bigr)\sqrt q+o(\sqrt q).
\]
Set $c=\sqrt{(4-D)/2}$ and then let $D\uparrow1$. The one-way area
ceiling gives $\limsup\sqrt q(W_1-W_2)/r\le2\sqrt6$.

For a finite example take $q=10^{-4}$, $h=122$, $D=97/100$ and
$n=9876$. The margins in the two time tests exceed $0.00193$ and
$0.00550$, respectively. The guarded width exceeds $19507.09959627r$.


\begin{thebibliography}{99}
\sloppy

\bibitem{Koopman1946}
B.~O. Koopman, \emph{Search and Screening}, Operations Evaluation Group Report No.~56, Office of the Chief of Naval Operations, Navy Department, Washington, D.C., 1946.

\bibitem{Balinski1965}
M.~L. Balinski, W.~Lucas, G.~Owen, and F.~M. Sand, \emph{A Study of Optimal Patrol and Transit Strategies in a Rectangular Barrier Zone Using Mathematical Games}, Final Report, Mathematica, Princeton, NJ, for the Office of Naval Research, Contract Nonr-4937(00), November 4, 1965. DTIC accession no.~AD0475962.

\bibitem{Washburn1976}
A.~R. Washburn, \emph{Patrolling a Channel Revisited}, Technical Report NPS55Ws75121, Naval Postgraduate School, Monterey, CA, January 1976.

\bibitem{Koopman1980}
B.~O. Koopman, \emph{Search and Screening: General Principles with Historical Applications}, Pergamon Press, New York, 1980.

\bibitem{Washburn1982}
A.~R. Washburn, ``On patrolling a channel,'' \emph{Naval Research Logistics Quarterly}, vol.~29, no.~4, pp.~609--615, 1982, \href{https://doi.org/10.1002/nav.3800290408}{\nolinkurl{doi:10.1002/nav.3800290408}}.

\bibitem{AndersenBouchetJackson1996}
L.~D. Andersen, A.~Bouchet, and B.~Jackson, ``Orthogonal A-trails of 4-regular graphs embedded in surfaces of low genus,'' \emph{Journal of Combinatorial Theory, Series B}, vol.~66, no.~2, pp.~232--246, 1996, \href{https://doi.org/10.1006/jctb.1996.0017}{\nolinkurl{doi:10.1006/jctb.1996.0017}}.

\bibitem{Washburn2009}
A.~R. Washburn, \emph{Barrier Games}, Technical Report NPS-OR-09-001, Naval Postgraduate School, Monterey, CA, June 2009.

\bibitem{Washburn2010}
A.~Washburn, ``Barrier games,'' \emph{Military Operations Research}, vol.~15, no.~3, pp.~31--41, 2010, \href{https://doi.org/10.5711/morj.15.3.31}{\nolinkurl{doi:10.5711/morj.15.3.31}}. Errata: vol.~21, no.~4, p.~83, 2016.

\bibitem{GalyaevMaslov2011}
A.~A. Galyaev and E.~P. Maslov, ``On the border patrolling problem,'' \emph{Journal of Computer and Systems Sciences International}, vol.~50, no.~5, pp.~837--846, 2011, \href{https://doi.org/10.1134/S1064230711050091}{\nolinkurl{doi:10.1134/S1064230711050091}}.

\bibitem{Hew2015}
P.~C. Hew, ``Linear and concentric arc patrols against smart evaders,'' \emph{Military Operations Research}, vol.~20, no.~3, pp.~39--48, 2015.

\bibitem{Hew2017}
P.~C. Hew, ``Barriers against smart intruders: Screening the limiting approach angle,'' \emph{Military Operations Research}, vol.~22, no.~1, pp.~5--18, 2017.

\bibitem{Rowling2019}
S.~J. Rowling, ``Effective sweep-width for barrier missions against an evasive target,'' \emph{Journal of Defense Modeling and Simulation}, vol.~18, no.~4, pp.~417--428, 2021 (first published online 2019), \href{https://doi.org/10.1177/1548512919875524}{\nolinkurl{doi:10.1177/1548512919875524}}.

\bibitem{AFP}
L.~Ambrosio, N.~Fusco and D.~Pallara, \emph{Functions of Bounded Variation and Free Discontinuity Problems}, Oxford University Press, Oxford, 2000.

\bibitem{Federer}
H.~Federer, \emph{Geometric Measure Theory}, Springer, Berlin, 1969.

\bibitem{Billingsley}
P.~Billingsley, \emph{Convergence of Probability Measures}, 2nd ed., Wiley, New York, 1999, \href{https://doi.org/10.1002/9780470316962}{\nolinkurl{doi:10.1002/9780470316962}}.

\end{thebibliography}
\end{document}